\pdfoutput=1
\documentclass[11pt]{amsart}
\usepackage[margin=2.5cm]{geometry}
\usepackage{amsmath,amssymb,amsthm,mathtools,enumitem}
\usepackage{tikz}
\usetikzlibrary{arrows.meta}
\usepackage[hyphens]{url}
\usepackage[colorlinks=true,linkcolor=blue,citecolor=blue,urlcolor=blue]{hyperref}

\newtheorem{conjecture}{Conjecture}[section]
\newtheorem{lemma}{Lemma}[section]
\newtheorem{proposition}[lemma]{Proposition}
\newtheorem{corollary}[lemma]{Corollary}
\newcommand{\T}{\mathcal T}
\newcommand{\C}{\mathbb C}
\newcommand{\Z}{\mathbb Z}
\newcommand{\Pb}{\mathbf P}
\newcommand{\Eb}{\mathbf E}
\newcommand{\cB}{\mathcal B}
\newcommand{\cH}{\mathcal H}
\newcommand{\one}{\mathbf 1}
\newcommand{\dbar}{\bar\partial}
\newcommand{\Arm}{\mathrm{Arm}}
\newcommand{\pc}{\xleftrightarrow{\,+\,}}
\newcommand{\sface}{\mathcal{F}}
\newcommand{\fce}{\mathsf{f}}
\newcommand{\ii}{\mathsf{i}}

\begin{document}
\sloppy

\title[Cardy's formula for the high-temperature Ising model]{A heuristic derivation of Cardy's formula for the
high-temperature Ising model}
\author[Yu Feng]{Yu Feng$^{\dagger}$}
\thanks{$^{\dagger}$yuproba@umich.edu, yufeng\_proba@163.com}
\makeatletter
\g@addto@macro\@setauthors{\par\vspace{4pt}{\centering\normalfont\small University of Michigan, Ann Arbor, USA\par}}
\makeatother

\begin{abstract}
The spin clusters of the Ising model on the triangular lattice at high temperature are
expected to have the same large-scale behavior as the clusters of critical percolation, and in particular
to satisfy Cardy’s formula. This conjecture is natural, since, by the Edwards--Sokal coupling, the spin
configuration is critical percolation on the clusters of its Fortuin--Kasteleyn representation,
which are typically small. We explain, in intuitive terms, a heuristic strategy to derive Cardy's
formula, and we postpone a detailed proof to a forthcoming paper. There is no exact color-switching
lemma for Smirnov's observable, but, if the probabilities that define it are computed conditionally
on the Fortuin--Kasteleyn clusters, it satisfies Smirnov's relations exactly on the graph obtained
by contracting every cluster to a point. Placing every cluster at its barycenter rotates, dilates,
and shears the faces of the lattice. By Green's formula, the holomorphicity of the scaling limit of
the observable reduces to a property of a single random interface, namely, that the expected contour
integrals of test functions along it vanish. We derive this property by homogenization along the
interfaces, using the key symmetry of the triangular lattice under the rotation
$z\mapsto e^{2\pi\ii/3}z$.
\end{abstract}


\date{\today}
\renewcommand*{\shortauthors}{}
\maketitle

\setcounter{tocdepth}{1}
\tableofcontents

\section{Introduction}\label{sec:intro}

\subsection{Motivation}\label{sec:motivation}

Two-dimensional critical percolation is expected to have a conformally invariant scaling limit,
which does not depend on the microscopic structures of the lattice. Let $\Omega$ be a Jordan domain
with boundary points $x_1,b,x_2,x_3$ in counterclockwise order, and let $\varphi$ be the conformal
map from $\Omega$ onto the upper half-plane that sends $(x_1,x_2,x_3)$ to $(0,1,\infty)$, so that
$\varphi(b)\in(0,1)$. For critical site percolation on the triangular lattice $\delta\T$, Cardy
predicted~\cite{Car92} and Smirnov proved~\cite{Smi01} that
\begin{equation}\label{eq:cardy-perc}
\lim_{\delta\to0}P_\delta\bigl([x_1,b]\leftrightarrow[x_2,x_3]\bigr)=\frac{\Gamma(2/3)}{\Gamma(1/3)^2}\int_0^{\varphi(b)}
t^{-2/3}(1-t)^{-2/3}\,dt ,
\end{equation}
where $[x_1,b]$ and $[x_2,x_3]$ are counterclockwise boundary arcs, and
$[x_1,b]\leftrightarrow[x_2,x_3]$ is the event that an open path joins them. The key step of
Smirnov's proof is an exact identity at the lattice level, obtained by a color-switching argument
that uses the independence of the site colors and the symmetry of the triangular lattice. This
identity shows that the discrete contour integrals of an observable vanish, up to boundary terms
that tend to zero as $\delta\to0$; by Morera's theorem, the scaling limits of the observable are
therefore holomorphic.

A natural test case for universality is the Ising model on $\delta\T$ at an inverse temperature
$\beta\in(0,\beta_c)$, where $\beta_c=\frac14\log3$. In the Edwards--Sokal coupling~\cite{ES88} of
the spins $\sigma$ with the Fortuin--Kasteleyn configuration $\omega$~\cite{FK72}, the spins are
obtained by coloring the clusters of $\omega$:
\begin{equation}\label{eq:ES-intro}
\Pb(\sigma\mid\omega)=2^{-k(\omega)}\,\one_{\{\sigma\text{ is constant on every cluster of }\omega\}},
\end{equation}
where $k(\omega)$ is the number of clusters, and these clusters are typically small for
$\beta<\beta_c$. Thus $\sigma$ is critical percolation on the small clusters of $\omega$, and the
spin clusters, that is, the maximal connected sets of sites with the same spin, are expected to
satisfy~\eqref{eq:cardy-perc}. In the language of the loop $O(n)$ model, the convergence of the
interfaces to those of critical percolation was conjectured by Smirnov~\cite[Section~2.2,
Conjecture~3]{Smi06}; see Section~\ref{sec:background}.

This conjecture connects two of the most studied lattice models. B\'alint, Camia, and
Meester~\cite{BCM10} used the representation~\eqref{eq:ES-intro} to show that the spin clusters are
critical at every $\beta<\beta_c$. To the best of our knowledge, Cardy's formula would be the first
conformally invariant quantity shown to be shared by the two models. In the scaling limit, the spin
clusters of an interacting model would then have the same crossing probabilities as those of
independent percolation, at every $\beta<\beta_c$. In renormalization-group terms, the Ising
interaction is irrelevant at the percolation fixed point, as far as the geometry of the clusters is
concerned (Section~\ref{sec:background}). Our strategy builds on this connection. Since the spins
are fair coins on the clusters of $\omega$, color switching survives at the level of clusters, which
is the starting point of our argument.

The goal of this paper is to explain a heuristic strategy to derive this conformal invariance; we
will give the details in a forthcoming paper~\cite{CardyLong}. In Sections~\ref{sec:ES}
and~\ref{sec:smirnov-beta0}, we will recall the model and Smirnov's argument, and in
Section~\ref{sec:ideas}, we will describe the main ideas of our strategy.

\subsection{Background and connections to related work}\label{sec:background}

The Cardy--Smirnov formula~\eqref{eq:cardy-perc} identifies the scaling limit of the percolation
interfaces as the Schramm--Loewner evolution SLE$_6$~\cite{Sch00,Smi01,CN07}, and it also leads to
the identification of the full scaling limit~\cite{CN06}. The universality of Cardy's formula was
tested numerically by Langlands, Pouliot, and Saint-Aubin~\cite{LPSA94}. The formula is known for
certain correlated models of triangular type~\cite{CL07,BCL10a,BCL10b}, where the proofs rest on
color-switching identities that hold exactly at the lattice level, and for models obtained from
Bernoulli percolation by local dynamics~\cite{CNS02,CNS04}, whose macroscopic cluster boundaries can
be coupled to those of Bernoulli percolation. Khristoforov and Smirnov~\cite{KS21} emphasized that
the conformal invariance of percolation remains open on other graphs. For the high-temperature Ising
model, neither an identity of this type nor such a coupling is known.

Beffara studied geometric extensions of Smirnov's argument, examining in~\cite{Bef08} the role of
the embedding of the lattice in the plane and the local defects of analyticity that it creates, and
in~\cite{Bef13} periodic graphs obtained by piecing together large triangular regions of the
triangular lattice. The observables of~\cite{KS21,KSS25}, built from loop configurations with
disorders and from the components of the complement of an interface, are closely related to the
observable of Section~\ref{sec:ideas-switching}. Homogenization of random Beltrami equations also
leads to deterministic effective complex structures~\cite{ARST22,IM23}. In our setting, the
microscopic distortion may collapse or reverse faces, and it is probed along interfaces that depend
on the same random medium; our strategy combines color switching on clusters, orthogonality
identities, and an energy identity to control this dependence.

The conjecture originates in the loop $O(n)$ model on the hexagonal lattice~\cite{DMNS81}, in which
a configuration of disjoint loops has a weight proportional to
$n^{\#\text{loops}}x^{\text{total length}}$. For $n\in[0,2]$, Nienhuis~\cite{Nie82,Nie84} predicted
the critical value $x_c(n)=(2+\sqrt{2-n})^{-1/2}$. For every $x>x_c(n)$, the model is expected to
have the scaling limit of its dense phase. For $n=1$, the loops are the interfaces of the Ising
model on $\T$ with $x=e^{-2\beta}$, so that $x_c(1)=1/\sqrt3$ corresponds to $\beta_c$, and the
dense phase contains critical percolation, which is the case $x=1$. Accordingly, Smirnov
conjectured~\cite[Section~2.2, Conjecture~3]{Smi06} that, for $n=1$ and every $x>x_c(1)$, the
interfaces converge to SLE$_6$, as for critical percolation. This includes the Ising model on $\T$
at every $\beta\in(0,\beta_c)$, for which $x\in(x_c(1),1)$.

The conjecture is also supported by a renormalization-group argument. Let $\nu=\frac43$ be the
correlation-length exponent of percolation. By the criterion of Harris~\cite{Har74}, extended to
correlated percolation by Weinrib~\cite{Wei84}, correlations that decay faster than
$r^{-2/\nu}=r^{-3/2}$ are irrelevant at the percolation fixed point. In our setting, the
correlations decay exponentially, so the model should flow to this fixed point, whose central charge
is $c=0$. As $\beta\to\beta_c$, the correlation length $\xi(\beta)$ tends to infinity. At scales
smaller than $\xi(\beta)$, the spin clusters look like those of the critical Ising model, which are
described by the tricritical point of the $q\to1$ Potts model~\cite{SV89} (see also~\cite{CK80}) and
whose interfaces converge to SLE$_3$~\cite{CDCHKS14}; at larger scales, they should look like
percolation clusters. Delfino~\cite{Del09} described this crossover, from the correlated fixed point
to that of percolation, by an integrable field theory. The mechanism of Section~\ref{sec:ideas}
concerns the scales larger than $\xi(\beta)$, and we do not know how the constant $\eta$ of
Section~\ref{sec:ideas-selfavg} behaves as $\beta\to\beta_c$.

For the high-temperature Ising model on $\T$, B\'alint, Camia, and Meester~\cite{BCM10} showed that,
if the clusters of $\omega$ are colored plus with probability $\theta$ instead of $1/2$
in~\eqref{eq:ES-intro}, a phase transition occurs exactly at $\theta=1/2$. The triangular lattice
matters here because it is self-matching, so that a quadrilateral is crossed by plus sites unless
the other pair of sides is joined by minus sites, and, together with the symmetry between the two
colors, this makes the spin clusters critical at every $\beta<\beta_c$. Recently, Tassion and
Vanneuville~\cite{TV25} proved that, for every $\beta<\beta_c$, the event that a box is crossed by a
path of plus spins is noise sensitive under Glauber dynamics; they also observed that, for this
model, the Russo--Seymour--Welsh (RSW) estimates, which go back to~\cite{Rus78,SW78}, follow
from~\cite{KT23}. Strong RSW estimates, uniform in the boundary conditions, will appear
in~\cite{SRSW}; for the critical random-cluster model, estimates of this type were proved
in~\cite{CDH16} for $q=2$ and in~\cite{DMT21} for $1\le q\le4$. The exponential decay of the
Fortuin--Kasteleyn connectivities follows from~\cite{ABF87,DRT19}, and, through the Edwards--Sokal
coupling, it gives the exponential spatial mixing that replaces independence
(Section~\ref{sec:ideas-selfavg}); see also the correlation inequality of Ding, Song, and
Sun~\cite{DSS23}. The arm bounds of Section~\ref{sec:returns} will appear in~\cite{SRSW}. To the
best of our knowledge, no conformal invariance result is known for this model at any
$\beta\in(0,\beta_c)$.

\subsection{The high-temperature Ising model and Cardy's formula}\label{sec:ES}

The sites of $\T$ are the points $m+ne^{\ii\pi/3}$, $m,n\in\Z$, and its faces are equilateral
triangles. We write $\cH$ for the hexagonal lattice dual to $\T$, whose vertices are the faces of
$\T$; every edge $e$ of $\T$ is crossed by one edge $e^*$ of $\cH$, and similarly for $\delta\T$ and
$\delta\cH$. A discrete domain $\Omega_\delta$ is a finite union of faces of $\delta\T$ whose
boundary is a simple closed curve; thus $\Omega_\delta$ is a closed topological disk. Its sites are
the vertices of its faces, and the Ising model on $\Omega_\delta$ with free boundary conditions is
the probability measure on $\{-1,+1\}^{\Omega_\delta\cap\delta\T}$ given by
\begin{equation*}
\Pb_{\Omega_\delta,\beta}(\sigma)=\frac1Z\exp\Bigl(\beta\sum_{zz'}\sigma_z\sigma_{z'}\Bigr),
\end{equation*}
where the sum runs over the edges $zz'$ of $\Omega_\delta$. Here and throughout, $zz'$ denotes the
edge with endpoints $z$ and $z'$, and $zz'z''$ the face with vertices $z$, $z'$, and $z''$. In the
Edwards--Sokal coupling~\cite{ES88}, given $\sigma$, each edge $zz'$ is an open bond with
probability $(1-e^{-2\beta})\one_{\{\sigma_z=\sigma_{z'}\}}$, independently of the other edges, and,
given the configuration $\omega$ of open bonds, $\sigma$ is given by~\eqref{eq:ES-intro}. We write
$P^\omega$ and $E^\omega$ for the conditional probability and expectation given $\omega$. For
$\beta<\beta_c$, by~\cite[Theorem~1.2]{DRT19}, applied to the random-cluster model with $q=2$, the
probability that the cluster of a site reaches distance $n\delta$ is at most $e^{-cn}$, even with
wired boundary conditions at this distance, for some $c=c(\beta)>0$ (see also~\cite{ABF87} for the
Ising model). In particular,
\begin{equation*}
\Pb\bigl(z\text{ and }z'\text{ lie in one cluster of }\omega\bigr)\le C\,e^{-c|z-z'|/\delta},
\end{equation*}
so that all the clusters of $\omega$ that meet a bounded region have diameter at most
\[\ell:=\delta\log^2(1/\delta),\]
except on an event of probability at most
\[C\delta^{-2}\log^2(1/\delta)\,e^{-\frac c2\log^2(1/\delta)}=\delta^{\frac c2\log(1/\delta)-2+o(1)},\]
by a union bound over the $O(\delta^{-2})$ sites of the region and the sites at distance at least
$\ell/2$ from each of them. This bound is smaller than any power of $\delta$. The interfaces of
$\sigma$ are the curves formed by the edges $e^*$ of $\delta\cH$ for which $e$ joins sites of
different spins; in the plane, and on the sphere of Section~\ref{sec:averaged}, they are circuits.
Since an open bond joins two sites of the same cluster, which have the same spin, interfaces never
cross open bonds.

Let
\[\tau=e^{2\pi\ii/3}\quad\text{and}\quad w_k=\tau^k.\]
Let $\Delta$ be the triangle with vertices $w_1,w_2,w_3$, and let $\Phi_\Omega$ be the conformal map
from $\Omega$ onto the interior of $\Delta$ with $\Phi_\Omega(x_k)=w_k$, $k\in\{1,2,3\}$. Let
$\Omega_\delta$ be discrete domains with boundary sites $x^\delta_1,b^\delta,x^\delta_2,x^\delta_3$
in counterclockwise order, which approximate $(\Omega;x_1,b,x_2,x_3)$ as $\delta\to0$ in the
following sense: $\partial\Omega_\delta$ converges to $\partial\Omega$ uniformly as a parametrized
curve, and $x^\delta_1,b^\delta,x^\delta_2,x^\delta_3$ converge to $x_1,b,x_2,x_3$. We write $[x,y]$
for the counterclockwise boundary arc from $x$ to $y$, and we will often drop the superscript
$\delta$ of the marked sites $x^\delta_1,b^\delta,x^\delta_2,x^\delta_3$.

\begin{conjecture}\label{conj:cardy}
For every $\beta\in(0,\beta_c)$,
\begin{equation}\label{eq:cardy}
\lim_{\delta\to0}\Pb_{\Omega_\delta,\beta}\bigl([x^\delta_1,b^\delta]\pc[x^\delta_2,x^\delta_3]\bigr)
=\frac{|\Phi_\Omega(b)-w_1|}{|w_2-w_1|},
\end{equation}
where $[x^\delta_1,b^\delta]\pc[x^\delta_2,x^\delta_3]$ denotes the event that a path of plus spins
in $\Omega_\delta$ joins $[x^\delta_1,b^\delta]$ to $[x^\delta_2,x^\delta_3]$.
\end{conjecture}

By the Schwarz--Christoffel formula, the right-hand side of~\eqref{eq:cardy} is that
of~\eqref{eq:cardy-perc}, in Carleson's form. It does not depend on $\beta$.
Conjecture~\ref{conj:cardy} is the counterpart, for crossing probabilities, of Smirnov's
conjecture~\cite[Section~2.2, Conjecture~3]{Smi06} that the interfaces converge to SLE$_6$
(Section~\ref{sec:background}).

The goal of this paper is to explain intuitively a strategy to prove Conjecture~\ref{conj:cardy}.
Except for the exact identities, which we will prove in Section~\ref{sec:switching} and
Appendix~\ref{app:exact}, and for the results that we cite, the arguments are heuristic; the details
will be given in~\cite{CardyLong}.

\subsection{Smirnov's argument at infinite temperature}\label{sec:smirnov-beta0}

\begin{figure}[t]
\centering
\begin{tikzpicture}[scale=0.82,
  bd/.style={line width=0.45pt},
  pl/.style={line width=1.4pt},
  mi/.style={line width=1.4pt,dash pattern=on 3.2pt off 2.2pt},
  hx/.style={line width=0.35pt,black!55},
  ar/.style={line width=0.6pt,-{Stealth[length=4.2pt,width=3.2pt]}}]
\def\domain{plot[smooth cycle,tension=0.62] coordinates {(4,0.3) (4.15,1.55) (4.5,2.35) (3.85,2.4)
  (3.1,1.95) (2.35,2.4) (1.45,2.8) (0.45,2.35) (0.05,1.1) (0.35,-0.25) (1.25,-1) (2.55,-1.3)
  (3.6,-0.85)}}
\def\base{%
  \draw[bd] \domain;
  \fill (4,0.3) circle (1.7pt) node[right=1pt] {$x_1$};
  \fill (0.45,2.35) circle (1.7pt) node[above left=-1pt] {$x_2$};
  \fill (1.25,-1) circle (1.7pt) node[below left=-1pt] {$x_3$};
  \foreach \c in {(1.69,1.087),(2.62,0.55),(1.69,0.013)}
    \draw[hx] \c ++(0.62,0) -- ++(120:0.62) -- ++(180:0.62) -- ++(240:0.62) -- ++(300:0.62) --
      ++(0:0.62) -- cycle;
  \draw[ar] (2,0.55) -- (2.31,1.087) node[above right=-2.5pt] {\small $\hat e$};
  \draw[ar] (2,0.55) -- (1.38,0.55) node[above left=-3pt] {\small $\tau\hat e$};
  \draw[ar] (2,0.55) -- (2.31,0.013) node[below right=-2.5pt] {\small $\tau^2\hat e$};
  \node at (1.98,0.33) {\small $z$};}
\def\toparm#1{\draw[#1] (2.3,2.7) .. controls (2.55,2) and (1.45,2.05) .. (1.85,1.6) .. controls
  (2.1,1.35) and (1.75,1.2) .. (1.69,1.087);}
\def\rightarm#1{\draw[#1] (2.62,0.55) .. controls (3.25,0.4) and (3.55,1.25) .. (3.5,0.5) ..
  controls (3.45,-0.15) and (2.65,-0.5) .. (3,-0.85) .. controls (3.15,-1) and (3.35,-1.1) ..
  (3.5,-1.3);}
\def\leftarm#1{\draw[#1] (1.69,0.013) .. controls (1.55,-0.4) and (1.05,-0.6) .. (0.75,-0.3) ..
  controls (0.5,-0.05) and (0.65,0.45) .. (0.3,0.55) .. controls (0.1,0.6) and (0,0.5) ..
  (-0.2,0.5);}
\def\panel#1#2#3#4#5{%
  \base
  \begin{scope}\clip \domain; \toparm{#1} \draw[#2] (1.69,1.087) -- #5; \rightarm{#3} \leftarm{#4}
  \end{scope}}
\def\pdot#1{\fill #1 circle (1.6pt);}
\def\mdot#1{\filldraw[fill=white,line width=0.6pt] #1 circle (1.6pt);}
\begin{scope}
\panel{pl}{pl}{pl}{mi}{(2.62,0.55)}
\pdot{(1.69,1.087)} \pdot{(2.62,0.55)} \mdot{(1.69,0.013)}
\node[above] at (1.45,2.8) {$A_1$};
\node[left] at (0.02,1.5) {$A_2$};
\node[right] at (3.92,-0.5) {$A_3$};
\node at (2.1,-1.95) {(a)};
\end{scope}
\begin{scope}[xshift=5.6cm]
\panel{mi}{mi}{pl}{mi}{(1.69,0.013)}
\mdot{(1.69,1.087)} \pdot{(2.62,0.55)} \mdot{(1.69,0.013)}
\node at (2.1,-1.95) {(b)};
\end{scope}
\begin{scope}[xshift=11.2cm]
\panel{pl}{pl}{mi}{pl}{(1.69,0.013)}
\pdot{(1.69,1.087)} \mdot{(2.62,0.55)} \pdot{(1.69,0.013)}
\node at (2.1,-1.95) {(c)};
\end{scope}
\end{tikzpicture}
\caption{Color switching at $\beta=0$. Solid lines are plus paths, dashed lines minus paths, and the
dots in the hexagons around $z$ indicate the colors of these hexagons. (a)
$E_1(z+\varepsilon\hat e)\setminus E_1(z)$. (b) After exploring the extremal arms to $A_2$ and
$A_3$, closest to $x_3$, and reversing the colors of the unexplored sites. (c) The global color flip
of (b), which is $E_2(z+\varepsilon\tau\hat e)\setminus E_2(z)$.}\label{fig:smirnov}
\end{figure}
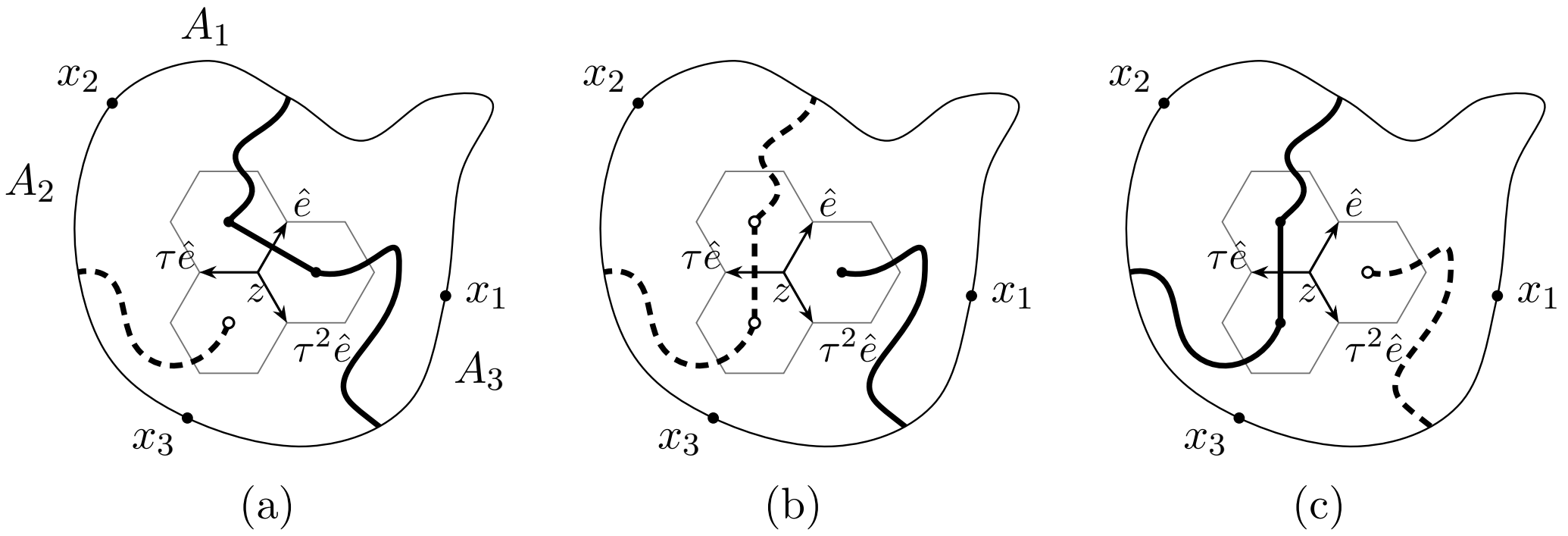

Let $\beta=0$, and let $A_k=[x_k,x_{k+1}]$, with indices modulo $3$. For a vertex $z$ of
$\delta\cH$, let $E_k(z)$ be the event that a simple plus path joins $A_{k-1}$ to $A_k$ and
separates $z$ from $A_{k+1}$, and let
\begin{equation}\label{eq:H}
H(z)=\sum_{k=1}^3w_k\,P\bigl(E_k(z)\bigr).
\end{equation}
Neighboring vertices of $\delta\cH$ are at distance $\varepsilon=\delta/\sqrt3$. Smirnov's key
observation~\cite{Smi01} is that, for every vertex $z$ and every unit vector $\hat e$ such that
$z+\varepsilon\hat e$ is a neighboring vertex,
\begin{equation}\label{eq:smirnov}
P\bigl(E_1(z+\varepsilon\hat e)\setminus E_1(z)\bigr)=P\bigl(E_2(z+\varepsilon\tau\hat e)\setminus E_2(z)\bigr)
=P\bigl(E_3(z+\varepsilon\tau^2\hat e)\setminus E_3(z)\bigr).
\end{equation}
The first event means that plus arms join the two hexagons adjacent to the edge $zz'$, where
$z'=z+\varepsilon\hat e$, to $A_1$ and $A_3$, and a minus arm joins the third hexagon to $A_2$
(Figure~\ref{fig:smirnov}(a)). Choose the arms to $A_2$ and $A_3$ to be the extremal ones, closest
to the common endpoint $x_3$ of these arcs. These arms can be found by an exploration started at
$x_3$, which does not reveal the colors of the sites beyond them, so that these colors remain
independent and fair. Reversing these colors therefore preserves the measure, and it turns the arm
to $A_1$ into a minus arm (b); a global color flip then gives the second event (c).
By~\eqref{eq:smirnov}, $H$ is approximately discrete holomorphic, in the sense that its discrete
contour integrals vanish up to boundary terms, which are bounded by the probabilities of three arms.
Next to a site $u$ of $A_k$, only $E_k$ and $E_{k+1}$ can occur, and
\begin{equation}\label{eq:H-boundary}
H=w_k+(w_{k+1}-w_k)\,P\bigl([x_k,u]\pc A_{k+1}\bigr).
\end{equation}
Every subsequential limit of $H$ is therefore holomorphic, by Morera's theorem, and maps $A_k$ into
the side $[w_k,w_{k+1}]$ of $\Delta$, so that it is equal to $\Phi_\Omega$,
and~\eqref{eq:H-boundary} at $u=b$ gives~\eqref{eq:cardy} at $\beta=0$~\cite{Smi01}; see
also~\cite{Bef07}. The relation~\eqref{eq:smirnov} is exact at the lattice level.

\section{Description of the main ideas}\label{sec:ideas}

This section describes the whole strategy step by step, without detailed proofs;
Sections~\ref{sec:switching} to~\ref{sec:cardy} will add derivations and remarks. Here and
throughout, $X\approx Y$ denotes a heuristic approximation, in the sense that $X-Y$ should tend to
zero, in probability if $X$ and $Y$ are random, as $\delta\to0$ and then as the auxiliary scales
introduced below tend to zero. Moreover, $X\lesssim Y$ means that $X\le CY$ for a constant $C$, and
$X\asymp Y$ means that $X\lesssim Y$ and $Y\lesssim X$. We never use an exact exponent of the model
that relies on conformal invariance; the only exact exponents that we use are universal ones, such
as the five-arm exponent, which do not. The exponent $d$ that we will introduce in
Section~\ref{sec:ideas-switching} only serves the exposition; in the arguments, RSW upper bounds can
replace it.

\subsection{A randomly distorted lattice}\label{sec:ideas-distorted}

By~\eqref{eq:ES-intro}, given $\omega$, all the sites of a cluster of $\omega$ have the same spin,
and distinct clusters have independent fair spins. Thus $\sigma$ is critical site percolation on the
contracted graph, obtained from $\delta\T$ by contracting every cluster of $\omega$ to a point. Its
vertices are the clusters, and every edge of $\delta\T$ that joins two different clusters gives an
edge between them; the open bonds are contracted, and a closed edge whose endpoints lie in one
cluster becomes a loop. This graph is not necessarily a triangulation:
\begin{itemize}[leftmargin=2em]
\item two neighboring clusters are usually joined by several parallel edges;
\item a face with exactly one open side becomes a face with two sides, and a face whose three sides
are open disappears;
\item a cluster that contains a circuit of open bonds separates the sites inside this circuit from
the rest of the graph.
\end{itemize}
No face has more than three sides, however, since identifying the sites of a cluster can only make
sides collapse. These features will not matter, since the face relation~\eqref{eq:face}, the defect
formula~\eqref{eq:defect}, and the energy identity~\eqref{eq:energy} below concern each face of
$\delta\T$ separately.

\emph{The barycenter map.} We draw the contracted graph in the plane by moving every site $z$ to the
barycenter of its cluster $K_z$,
\begin{equation*}
\cB(z)=\frac1{|K_z|}\sum_{z'\in K_z}z',
\end{equation*}
and by extending $\cB$ affinely to every face. Since the clusters have diameter at most
$\ell=\delta\log^2(1/\delta)$ with high probability (Section~\ref{sec:ES}), we have
$|\cB(z)-z|\le\ell$, so that $\cB$ is close to the identity at large scales. At the lattice level,
it distorts the faces near the clusters with at least two sites, and it maps every open bond to a
point (Figure~\ref{fig:barycenter}). The image of a face $\fce$ is the convex hull of the images of
its vertices; it is
\begin{itemize}[leftmargin=2em]
\item a triangle, if these three images are not collinear;
\item a segment, if exactly two vertices of $\fce$ lie in one cluster, or, in rare cases, if the
three barycenters are collinear;
\item a point, if the three vertices of $\fce$ lie in one cluster.
\end{itemize}
The edges between two neighboring clusters are all mapped to the segment between their barycenters
(Figure~\ref{fig:contraction}), and the images of distinct faces may overlap; the argument uses only
the images of the individual faces. Any rule that places every cluster at a point near it, depends
only on the cluster, and is covariant under the symmetries of $\T$ would do; the constant $\eta$ of
Section~\ref{sec:ideas-selfavg} depends on this choice, but the conclusion does not.

\begin{figure}[t]
\centering
\begin{tikzpicture}[scale=1.35,
  tri/.style={line width=0.7pt},
  cir/.style={line width=0.45pt}]
\def\site#1{\fill #1 circle (1.2pt);}
\begin{scope}
\draw[tri] (0,0) -- (1.2,0) -- (0.6,1.039) -- cycle;
\draw[cir] (0.6,0.346) circle (0.346);
\foreach \p in {(0,0),(1.2,0),(0.6,1.039)} {\site{\p}}
\node at (0.6,-0.42) {(a)};
\end{scope}
\begin{scope}[shift={(2.25,0)}]
\begin{scope}[rotate=20,scale=0.85]
\draw[tri] (0,0) -- (1.2,0) -- (0.6,1.039) -- cycle;
\draw[cir] (0.6,0.346) circle (0.346);
\foreach \p in {(0,0),(1.2,0),(0.6,1.039)} {\site{\p}}
\end{scope}
\node at (0.52,-0.42) {(b)};
\end{scope}
\begin{scope}[shift={(4.15,0)}]
\draw[tri] (0,0) -- (1.56,0) -- (0.78,0.727) -- cycle;
\draw[cir] (0.78,0.242) ellipse (0.45 and 0.242);
\foreach \p in {(0,0),(1.56,0),(0.78,0.727)} {\site{\p}}
\node at (0.78,-0.42) {(c)};
\end{scope}
\end{tikzpicture}
\caption{(a) An equilateral face with its inscribed circle. (b) Its image under a rotation and a
dilation $z\mapsto Az$: the face remains equilateral, the circle remains a circle, and Smirnov's
relation is preserved. (c) Its image under the shear $z\mapsto z+\mu\bar z$ with $\mu=0.3$, which
stretches the face horizontally by the factor $1.3$ and compresses it vertically by the factor
$0.7$: the face is no longer equilateral, the circle becomes an ellipse, and Smirnov's relation
acquires the defect in~\eqref{eq:defect}.}\label{fig:shear}
\end{figure}
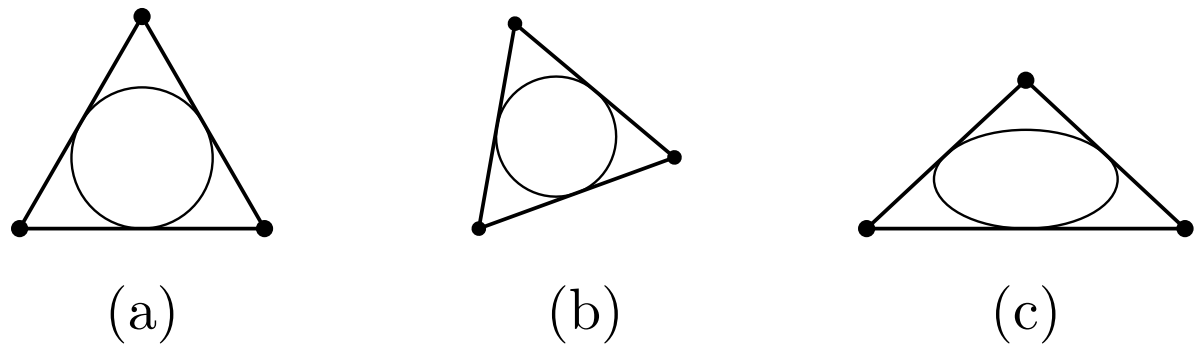

\emph{Rotations, dilations, and shears.} On a face $\fce$, the map $\cB$ is affine, of the form
$z\mapsto A_{\fce}z+A'_{\fce}\bar z+\mathrm{const}$. If $A_{\fce}\neq0$, it is, up to a translation,
\begin{equation*}
z\mapsto A_{\fce}\,(z+\mu_{\fce}\bar z),\qquad\mu_{\fce}=A'_{\fce}/A_{\fce},
\end{equation*}
namely a shear $z\mapsto z+\mu_{\fce}\bar z$, followed by a rotation and a dilation
$z\mapsto A_{\fce}z$, which keep the face equilateral (Figures~\ref{fig:shear} and~\ref{fig:rds}).
Writing $\mu_{\fce}=|\mu_{\fce}|e^{2\ii\theta}$, we have, for real $a$ and $b$,
\begin{equation*}
\begin{aligned}
v&=a\,e^{\ii\theta}+b\,\ii e^{\ii\theta},\\
v+\mu_{\fce}\bar v&=\bigl(1+|\mu_{\fce}|\bigr)\,a\,e^{\ii\theta}+\bigl(1-|\mu_{\fce}|\bigr)\,b\,\ii e^{\ii\theta},
\end{aligned}
\end{equation*}
so that, in these orthonormal coordinates, the shear is the linear map
$(a,b)\mapsto\bigl((1+|\mu_{\fce}|)\,a,(1-|\mu_{\fce}|)\,b\bigr)$. Hence the image of $\fce$ is
degenerate if $|\mu_{\fce}|=1$, and it is inverted, with the orientation of its vertices reversed,
if $|\mu_{\fce}|>1$, since the factor $1-|\mu_{\fce}|$ is then negative and the shear also reflects
the face. Inverted images do occur, for instance when the clusters of two vertices of $\fce$ extend
in opposite directions, far enough that their barycenters lie in the reverse order. The argument,
however, never uses the ellipticity of the distortion.

\emph{The two parts of the distortion.} Let $e_0(\fce),e_1(\fce),e_2(\fce)$ be the sides of $\fce$
in counterclockwise order, starting from any side; another starting side only multiplies the
relations below by a power of $\tau$. The index $j$ records the direction of the side. More
precisely, the edge of $\delta\cH$ that crosses $e_j(\fce)$ leaves the vertex $\fce$ of $\delta\cH$
in the direction $\tau^j\hat e$, where $\hat e$ is the direction of the edge that crosses
$e_0(\fce)$, as in~\eqref{eq:smirnov}. The vectors $d_j(\fce)$ that join the images of the endpoints
of $e_j(\fce)$, oriented clockwise around $\fce$, add up to zero, so that they split uniquely as
\begin{equation}\label{eq:pq}
d_j(\fce)=p_{\fce}\,\tau^j+q_{\fce}\,\bar\tau^j,\qquad j\in\{0,1,2\},
\end{equation}
with $p_{\fce}=\frac13\sum_j\bar\tau^jd_j(\fce)$ and $q_{\fce}=\frac13\sum_j\tau^jd_j(\fce)$. Let
$p^0_{\fce}$ be the value of $p_{\fce}$ for the undistorted face, for which $q_{\fce}=0$ and
$|p^0_{\fce}|=\delta$. As the affine map sends every side vector $v$ of $\fce$ to
$A_{\fce}v+A'_{\fce}\bar v$,
\begin{equation*}
\begin{aligned}
p_{\fce}&=A_{\fce}\,p^0_{\fce},\\
q_{\fce}&=A'_{\fce}\,\overline{p^0_{\fce}},\\
\tfrac{\sqrt3}4\bigl(|p_{\fce}|^2-|q_{\fce}|^2\bigr)&=\text{signed area of the image of }\fce,
\end{aligned}
\end{equation*}
since $|A_{\fce}|^2-|A'_{\fce}|^2$ is the Jacobian of the map and $\frac{\sqrt3}4\delta^2$ is the
area of $\fce$. Thus $p_{\fce}$ describes the rotation and the dilation of $\fce$, and $q_{\fce}$
its shear. In particular, $q_{\fce}=0$ if and only if the image of $\fce$ is an equilateral triangle
with the orientation of $\fce$, or a point, and the signed area is negative exactly when the image
is inverted. At $\beta=0$, there are no open bonds, $\cB$ is the identity, and $q_{\fce}=0$ for
every face.

\begin{figure}[t]
\centering
\begin{tikzpicture}[scale=1.25,
  bond/.style={line width=1.5pt,line cap=round},
  img/.style={line width=0.9pt,dash pattern=on 2.6pt off 1.6pt},
  face/.style={line width=0.5pt}]
\def\site#1{\fill #1 circle (1.3pt);}
\def\bary#1{\filldraw[fill=white,line width=0.6pt] #1 circle (1.9pt);}
\begin{scope}
\draw[face] (0,0) -- (1.2,0) -- (0.6,1.039) -- cycle;
\draw[bond] (0,0) -- (0.6,-1.039);
\draw[bond] (0.6,1.039) -- (-0.6,1.039);
\foreach \p in {(0,0),(1.2,0),(0.6,1.039),(0.6,-1.039),(-0.6,1.039)} {\site{\p}}
\draw[img] (0.3,-0.52) -- (1.2,0) -- (0,1.039) -- cycle;
\bary{(0.3,-0.52)} \bary{(1.2,0)} \bary{(0,1.039)}
\node[font=\small] at (0.62,0.3) {$\fce$};
\node at (0.3,-1.6) {(a)};
\end{scope}
\begin{scope}[xshift=4cm,yshift=-0.5cm]
\draw[face] (0,0) -- (1.2,0) -- (0.6,1.039) -- cycle;
\draw[bond] (0,0) -- (1.2,0);
\draw[bond] (0.6,1.039) -- (1.8,1.039);
\foreach \p in {(0,0),(1.2,0),(0.6,1.039),(1.8,1.039)} {\site{\p}}
\draw[img] (0.6,0) -- (1.2,1.039);
\bary{(0.6,0)} \bary{(1.2,1.039)}
\node[font=\small] at (0.4,0.35) {$\fce$};
\node at (0.7,-1.1) {(b)};
\end{scope}
\end{tikzpicture}
\caption{(a) A face $\fce$ whose vertices lie in three clusters; open bonds are thick, and one
cluster is a single site. The map $\cB$ sends every site to the barycenter of its cluster (circles),
and $\fce$ to the dashed triangle. (b) If a side of $\fce$ is open, the image of $\fce$ is
degenerate.}\label{fig:barycenter}
\end{figure}
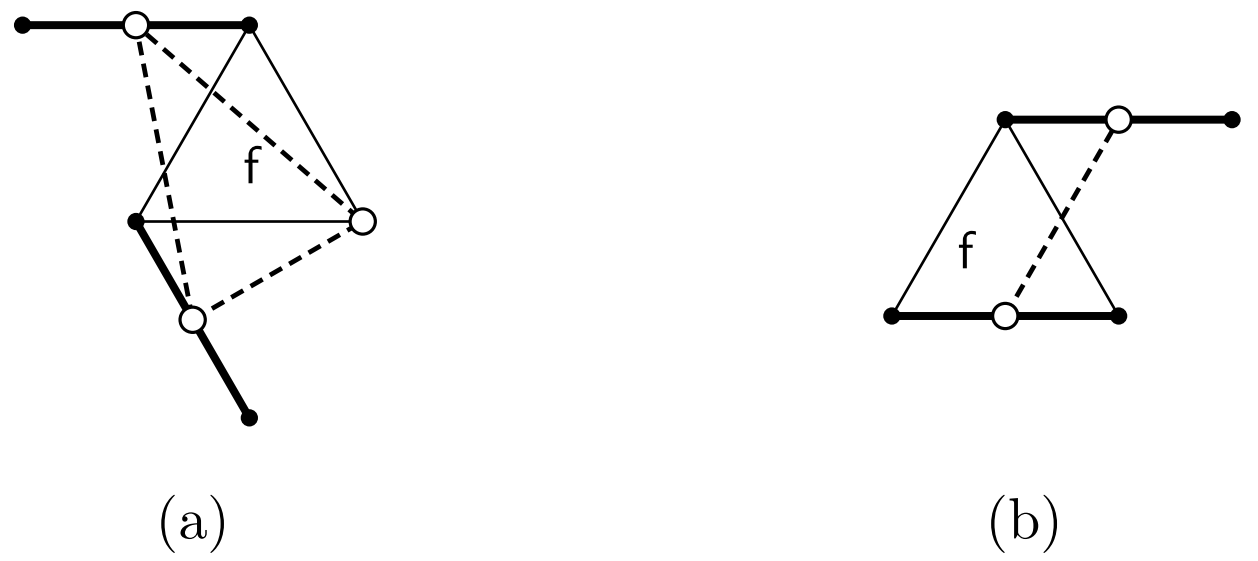

\begin{figure}[t]
\centering
\begin{tikzpicture}[scale=1.15,
  lat/.style={line width=0.35pt},
  bond/.style={line width=1.5pt,line cap=round},
  hl/.style={line width=0.9pt,dash pattern=on 2.6pt off 1.6pt},
  img/.style={line width=0.4pt}]
\def\site#1{\fill #1 circle (1.3pt);}
\def\bary#1{\filldraw[fill=white,line width=0.6pt] #1 circle (1.9pt);}
\begin{scope}
\draw[lat] (2,0) -- (4,0);
\draw[lat] (1.5,0.866) -- (2.5,0.866);
\draw[lat] (3.5,0.866) -- (4.5,0.866);
\draw[lat] (1,1.732) -- (2,1.732);
\draw[lat] (3,1.732) -- (5,1.732);
\draw[lat] (1.5,2.598) -- (4.5,2.598);
\draw[lat] (2,3.464) -- (4,3.464);
\draw[lat] (2,0) -- (2.5,0.866);
\draw[lat] (3,0) -- (4.5,2.598);
\draw[lat] (4,0) -- (5,1.732);
\draw[lat] (1.5,0.866) -- (2,1.732);
\draw[lat] (1,1.732) -- (2,3.464);
\draw[lat] (3,1.732) -- (4,3.464);
\draw[lat] (2.5,2.598) -- (3,3.464);
\draw[lat] (2,0) -- (1,1.732);
\draw[lat] (3,0) -- (1.5,2.598);
\draw[lat] (4,0) -- (3.5,0.866);
\draw[lat] (4.5,0.866) -- (3,3.464);
\draw[lat] (5,1.732) -- (4,3.464);
\draw[lat] (2.5,2.598) -- (2,3.464);
\draw[hl] (2,1.732) -- (2.5,2.598);
\draw[hl] (2.5,2.598) -- (3,1.732);
\draw[hl] (2.5,0.866) -- (3,1.732);
\draw[hl] (2.5,0.866) -- (3.5,0.866);
\draw[hl] (3,1.732) -- (3.5,0.866);
\draw[bond] (2,1.732) -- (3,1.732);
\draw[bond] (3,1.732) -- (3.5,2.598);
\draw[bond] (3,0) -- (3.5,0.866);
\draw[bond] (4,1.732) -- (5,1.732);
\draw[bond] (4,1.732) -- (4.5,2.598);
\draw[bond] (1.5,2.598) -- (2.5,2.598);
\site{(2,0)}
\site{(3,0)}
\site{(4,0)}
\site{(1.5,0.866)}
\site{(2.5,0.866)}
\site{(3.5,0.866)}
\site{(4.5,0.866)}
\site{(1,1.732)}
\site{(2,1.732)}
\site{(3,1.732)}
\site{(4,1.732)}
\site{(5,1.732)}
\site{(1.5,2.598)}
\site{(2.5,2.598)}
\site{(3.5,2.598)}
\site{(4.5,2.598)}
\site{(2,3.464)}
\site{(3,3.464)}
\site{(4,3.464)}
\bary{(3.25,0.433)}
\bary{(2.833,2.021)}
\bary{(4.5,2.021)}
\bary{(2,2.598)}
\node[font=\small] at (3,1.125) {$\fce$};
\node[font=\small] at (2.35,2.051) {$\fce'$};
\node at (3,-0.45) {(a)};
\end{scope}
\begin{scope}[xshift=5.7cm]
\draw[img] (1,1.732) -- (2,2.598);
\draw[img] (1,1.732) -- (1.5,0.866);
\draw[img] (1,1.732) -- (2.833,2.021);
\draw[img] (2,2.598) -- (2,3.464);
\draw[hl] (2,2.598) -- (2.833,2.021);
\draw[img] (2,3.464) -- (3,3.464);
\draw[img] (1.5,0.866) -- (2.833,2.021);
\draw[img] (1.5,0.866) -- (2,0);
\draw[img] (1.5,0.866) -- (2.5,0.866);
\draw[hl] (2.833,2.021) -- (2.5,0.866);
\draw[img] (2,2.598) -- (3,3.464);
\draw[img] (3,3.464) -- (2.833,2.021);
\draw[img] (3,3.464) -- (4,3.464);
\draw[img] (2,0) -- (2.5,0.866);
\draw[img] (2,0) -- (3.25,0.433);
\draw[hl] (2.5,0.866) -- (3.25,0.433);
\draw[hl] (2.833,2.021) -- (3.25,0.433);
\draw[img] (2.833,2.021) -- (4.5,2.021);
\draw[img] (2.833,2.021) -- (4,3.464);
\draw[img] (4,3.464) -- (4.5,2.021);
\draw[img] (3.25,0.433) -- (4,0);
\draw[img] (3.25,0.433) -- (4.5,2.021);
\draw[img] (3.25,0.433) -- (4.5,0.866);
\draw[img] (4.5,2.021) -- (4.5,0.866);
\draw[img] (4,0) -- (4.5,0.866);
\site{(2,0)}
\bary{(3.25,0.433)}
\site{(4,0)}
\site{(1.5,0.866)}
\site{(2.5,0.866)}
\site{(4.5,0.866)}
\site{(1,1.732)}
\bary{(2.833,2.021)}
\bary{(4.5,2.021)}
\bary{(2,2.598)}
\site{(2,3.464)}
\site{(3,3.464)}
\site{(4,3.464)}
\node at (3,-0.45) {(b)};
\end{scope}
\end{tikzpicture}
\caption{The map $\cB$ on a patch of $\delta\T$. (a) Open bonds are thick. They form four clusters
with several sites, whose barycenters are marked by circles; every other site is a cluster by
itself. The vertices of the face $\fce$ lie in three different clusters, and the face $\fce'$ has an
open side. (b) The image under $\cB$: every site is moved to the barycenter of its cluster, and the
map is affine on every face. The face $\fce$ becomes a distorted triangle, and $\fce'$ becomes a
segment (both dashed). The four edges that join the two clusters containing the vertices of $\fce'$
are all mapped to this segment, and the face whose three vertices form one cluster is mapped to a
point.}\label{fig:contraction}
\end{figure}
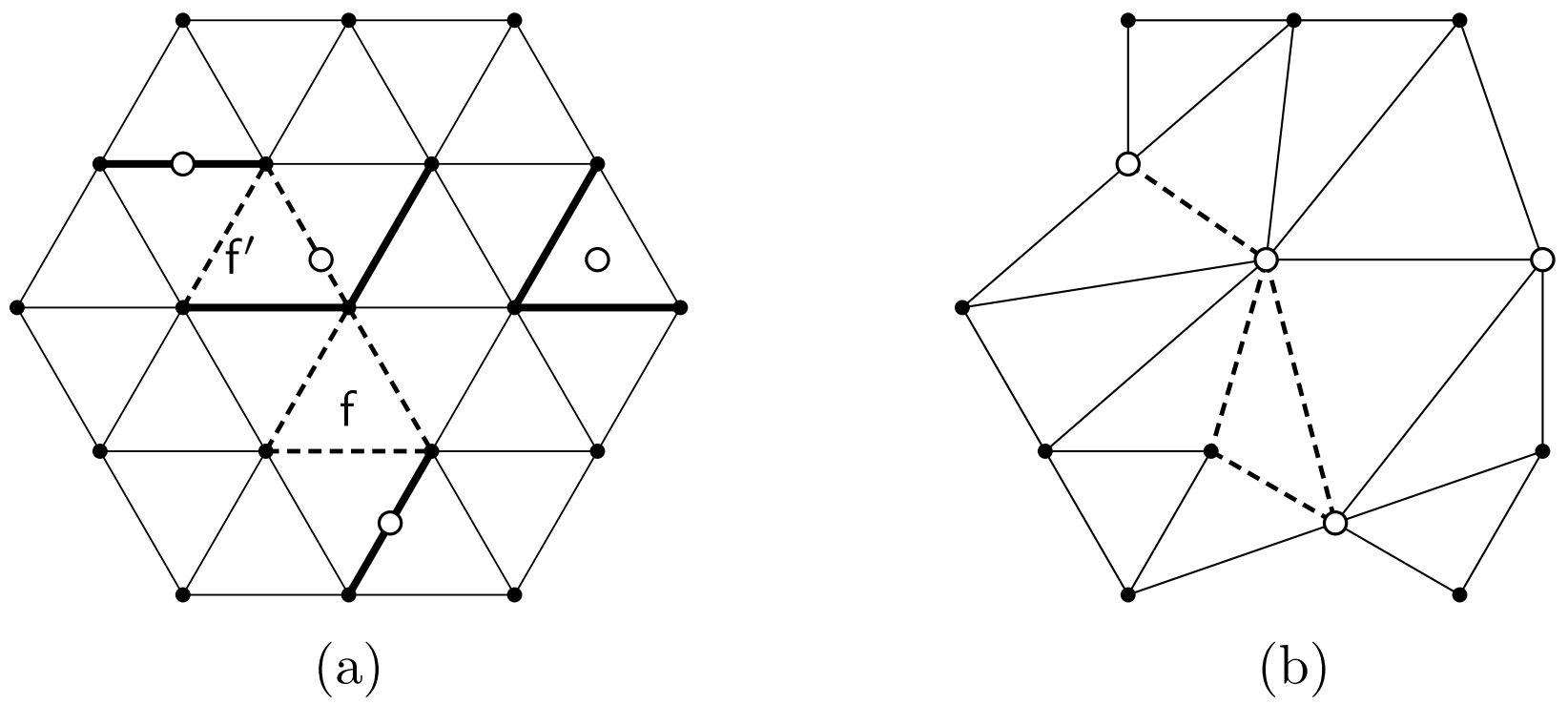

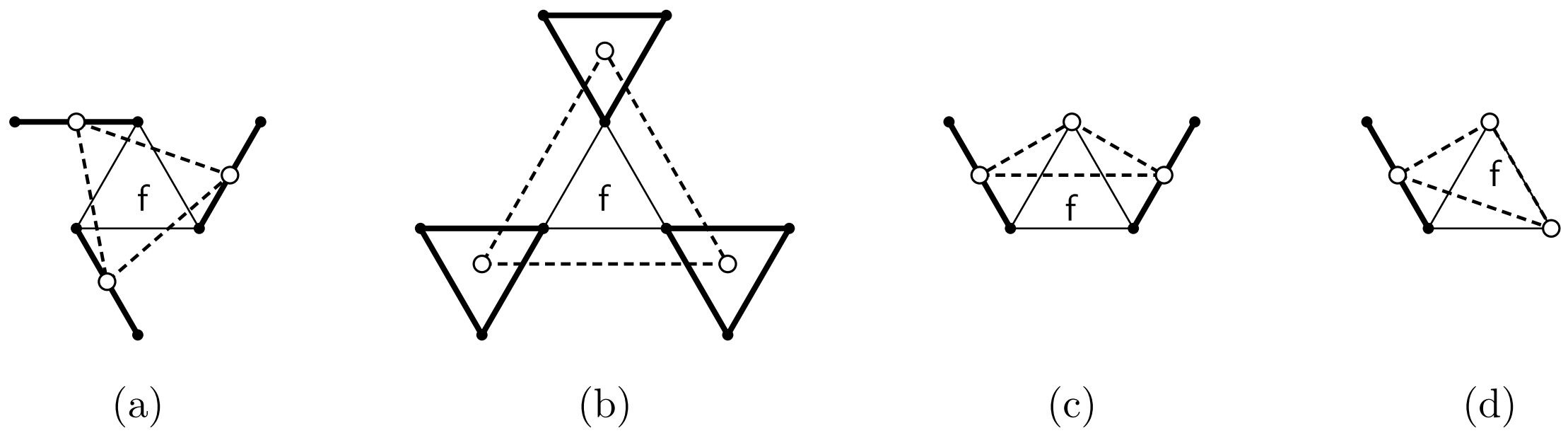
\begin{figure}[t]
\centering
\begin{tikzpicture}[scale=1.15,
  bond/.style={line width=1.5pt,line cap=round},
  img/.style={line width=0.9pt,dash pattern=on 2.6pt off 1.6pt},
  face/.style={line width=0.5pt}]
\def\site#1{\fill #1 circle (1.3pt);}
\def\bary#1{\filldraw[fill=white,line width=0.6pt] #1 circle (1.9pt);}
\begin{scope}
\draw[face] (0,0) -- (1,0) -- (0.5,0.866) -- cycle;
\draw[bond] (0,0) -- (0.5,-0.866);
\draw[bond] (1,0) -- (1.5,0.866);
\draw[bond] (0.5,0.866) -- (-0.5,0.866);
{\site{(0,0)}} {\site{(1,0)}} {\site{(0.5,0.866)}}
{\site{(0.5,-0.866)}} {\site{(1.5,0.866)}} {\site{(-0.5,0.866)}}
\draw[img] (0.25,-0.433) -- (1.25,0.433) -- (0,0.866) -- cycle;
\bary{(0.25,-0.433)} \bary{(1.25,0.433)} \bary{(0,0.866)}
\node[font=\small] at (0.55,0.25) {$\fce$};
\node at (0.5,-1.45) {(a)};
\end{scope}
\begin{scope}[xshift=3.8cm]
\draw[face] (0,0) -- (1,0) -- (0.5,0.866) -- cycle;
\draw[bond,line join=round] (-1,0) -- (0,0) -- (-0.5,-0.866) -- cycle;
\draw[bond,line join=round] (2,0) -- (1,0) -- (1.5,-0.866) -- cycle;
\draw[bond,line join=round] (0,1.732) -- (0.5,0.866) -- (1,1.732) -- cycle;
{\site{(0,0)}} {\site{(1,0)}} {\site{(0.5,0.866)}}
{\site{(-1,0)}} {\site{(-0.5,-0.866)}} {\site{(2,0)}} {\site{(1.5,-0.866)}} {\site{(0,1.732)}}
{\site{(1,1.732)}}
\draw[img] (-0.5,-0.289) -- (1.5,-0.289) -- (0.5,1.443) -- cycle;
\bary{(-0.5,-0.289)} \bary{(1.5,-0.289)} \bary{(0.5,1.443)}
\node[font=\small] at (0.5,0.25) {$\fce$};
\node at (0.5,-1.45) {(b)};
\end{scope}
\begin{scope}[xshift=7.6cm]
\draw[face] (0,0) -- (1,0) -- (0.5,0.866) -- cycle;
\draw[bond] (0,0) -- (-0.5,0.866);
\draw[bond] (1,0) -- (1.5,0.866);
{\site{(0,0)}} {\site{(1,0)}} {\site{(0.5,0.866)}}
{\site{(-0.5,0.866)}} {\site{(1.5,0.866)}}
\draw[img] (-0.25,0.433) -- (1.25,0.433) -- (0.5,0.866) -- cycle;
\bary{(-0.25,0.433)} \bary{(1.25,0.433)} \bary{(0.5,0.866)}
\node[font=\small] at (0.5,0.17) {$\fce$};
\node at (0.5,-1.45) {(c)};
\end{scope}
\begin{scope}[xshift=11cm]
\draw[face] (0,0) -- (1,0) -- (0.5,0.866) -- cycle;
\draw[bond] (0,0) -- (-0.5,0.866);
{\site{(0,0)}} {\site{(1,0)}} {\site{(0.5,0.866)}} {\site{(-0.5,0.866)}}
\draw[img] (-0.25,0.433) -- (1,0) -- (0.5,0.866) -- cycle;
\bary{(-0.25,0.433)} \bary{(1,0)} \bary{(0.5,0.866)}
\node[font=\small] at (0.55,0.45) {$\fce$};
\node at (0.5,-1.45) {(d)};
\end{scope}
\end{tikzpicture}
\caption{Clusters at the vertices of a face $\fce$, and the image of $\fce$ under $\cB$ (dashed);
open bonds are thick, and the barycenters are marked by circles. In (a), (b), and (c), the image is,
respectively, a rotated, a dilated, and a sheared copy of $\fce$. (a) Each vertex of $\fce$ is
joined by an open bond to the outer neighbor in the counterclockwise direction around $\fce$: the
image is $\fce$ rotated about its center by $\arctan(\sqrt3/2)\approx41^\circ$, together with a
dilation by the factor $\sqrt7/2$; a pure rotation requires larger clusters. (b) Each vertex forms a
cluster with its two outer neighbors: the image is $\fce$ dilated about its center by the factor
$2$. (c) The two lower vertices are joined to their upper outer neighbors: up to a translation, the
image is obtained from $\fce$ by the shear $z\mapsto z+\frac12\bar z$, which stretches $\fce$
horizontally by the factor $\frac32$ and compresses it vertically by the factor $\frac12$. (d) Only
the left bond of (c) is kept. This is the general case: the distortion has all three components, a
rotation by $\arctan\bigl(1/(2\sqrt3)\bigr)\approx16^\circ$ clockwise, a dilation by the factor
$\sqrt{13/12}\approx1.04$, and a shear with $|\mu_{\fce}|=1/\sqrt{13}\approx0.28$. Each panel shows
a single face of $\delta\T$; for the contracted graph, which is not necessarily a triangulation, see
Figure~\ref{fig:contraction}.}\label{fig:rds}
\end{figure}

\subsection{Color switching on clusters}\label{sec:ideas-switching}

For $\beta>0$, colors can no longer be switched site by site, since all the sites of a cluster of
$\omega$ have the same spin. What survives is color switching cluster by cluster, or equivalently
site by site on the contracted graph. Given $\omega$, the colors of the clusters are independent
fair coins, by~\eqref{eq:ES-intro}, so that reversing the colors of whole clusters preserves the
conditional law.

To use this, we write Smirnov's observable in terms of a single interface
(Section~\ref{sec:averaged}). We add three ghost vertices $s_1,s_2,s_3$ outside the domain, joined
to each other, with $s_k$ joined to every site of $A_k$, and we call ${\sface}=s_1s_2s_3$ the added
face. The new edges are never open bonds, so that the spins of $s_1,s_2,s_3$ are independent fair
coins, independent of the model inside the domain. Let
\begin{equation}\label{eq:wstar}
W=\sum_{k=1}^3w_k\,\one_{\{\sigma_{s_{k-1}}=\sigma_{s_k}\neq\sigma_{s_{k+1}}\}},
\end{equation}
which takes each of the values $w_1,w_2,w_3$ with probability $\frac14$, and the value $0$ if the
three spins are equal.

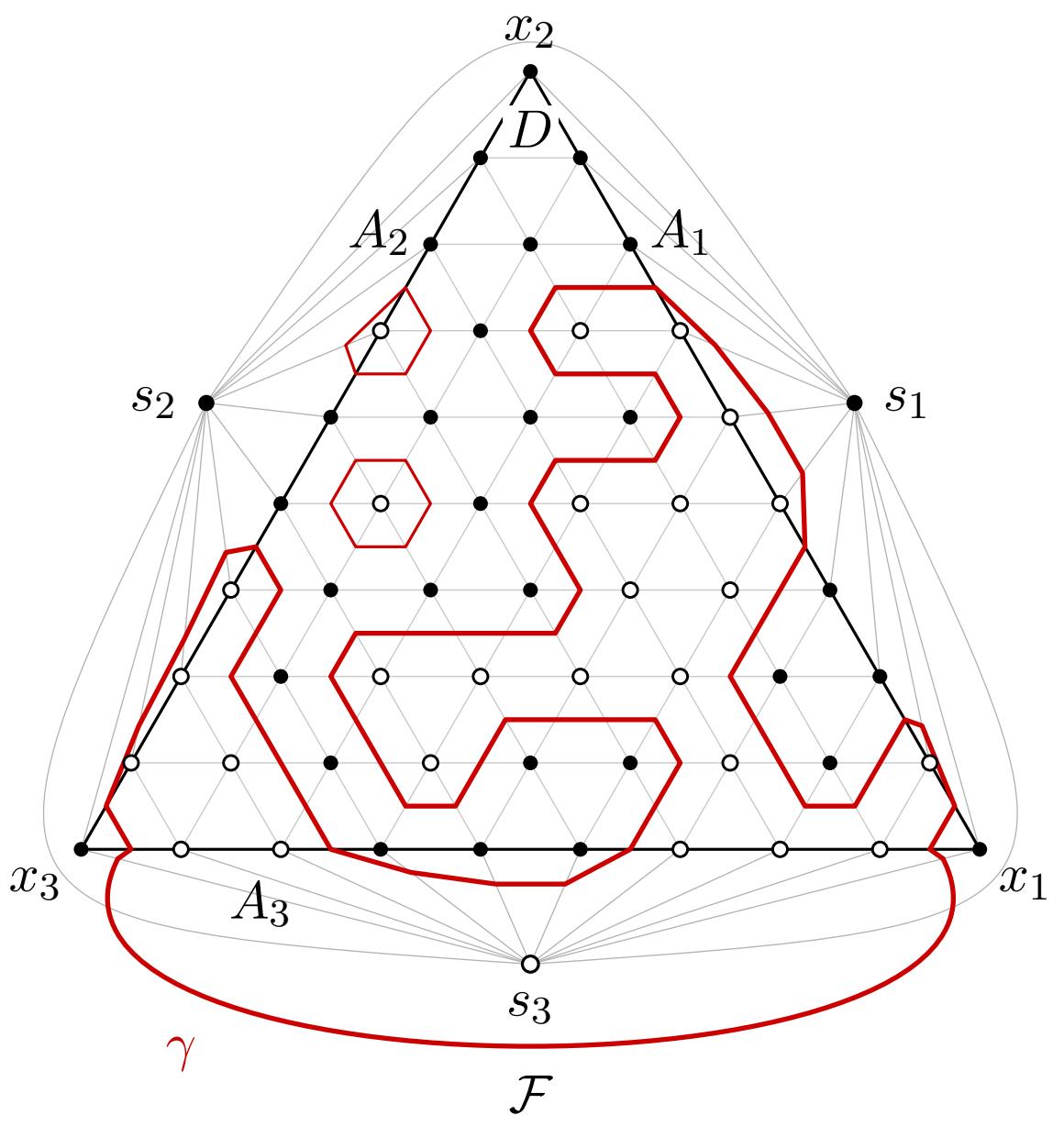
\begin{figure}[tp]
\centering
\begin{tikzpicture}[scale=0.72,
  lat/.style={line width=0.25pt,black!22},
  fan/.style={line width=0.25pt,black!30},
  loop/.style={line width=0.6pt,red!80!black,line join=round},
  Lstar/.style={line width=1pt,red!80!black,line join=round,line cap=round}]
\def\plus#1{\fill #1 circle (2.1pt);}
\def\minus#1{\filldraw[fill=white,line width=0.55pt] #1 circle (2.1pt);}
\coordinate (s1) at (7.746,4.472);
\coordinate (s2) at (1.254,4.472);
\coordinate (s3) at (4.5,-1.15);
\begin{scope}[x={(1,0)},y={(0.5,0.866)}]
\foreach \k in {0,...,8} {
  \draw[lat] (0,\k) -- (9-\k,\k);
  \draw[lat] (\k,0) -- (\k,9-\k);
  \draw[lat] (\k+1,0) -- (0,\k+1);
}
\end{scope}
\draw[line width=0.6pt] (0,0) -- (9,0) -- (4.5,7.794) -- cycle;
\begin{scope}[x={(1,0)},y={(0.5,0.866)}]
\foreach \k in {0,...,9} {
  \draw[fan] (s3) -- (\k,0);
  \draw[fan] (s1) -- (9-\k,\k);
  \draw[fan] (s2) -- (0,\k);
}
\end{scope}
\draw[fan] (1.254,4.472) .. controls (-1.299,-0.75) .. (4.5,-1.15);
\draw[fan] (4.5,-1.15) .. controls (10.299,-0.75) .. (7.746,4.472);
\draw[fan] (1.254,4.472) .. controls (4.5,9.294) .. (7.746,4.472);
\draw[loop] (2.75,3.031) -- (3.25,3.031) -- (3.5,3.464) -- (3.25,3.897) -- (2.75,3.897) --
  (2.5,3.464) -- cycle;
\draw[loop] (3.25,4.763) -- (3.5,5.196) -- (3.25,5.629) -- (2.649,5.051) -- (2.75,4.763) -- cycle;
\draw[Lstar] (8.632,-0.094) -- (8.5,0) -- (8.75,0.433) -- (8.422,1.238) -- (8.25,1.299) --
  (8,0.866) -- (7.75,0.433) -- (7.25,0.433) -- (7,0.866) -- (6.75,1.299) -- (6.5,1.732) --
  (6.75,2.165) -- (7,2.598) -- (7.25,3.031) -- (7.226,3.77) -- (6.878,4.373) -- (6.351,5.051) --
  (5.75,5.629) -- (5.25,5.629) -- (4.75,5.629) -- (4.5,5.196) -- (4.75,4.763) -- (5.25,4.763) --
  (5.75,4.763) -- (6,4.33) -- (5.75,3.897) -- (5.25,3.897) -- (4.75,3.897) -- (4.5,3.464) --
  (4.75,3.031) -- (5,2.598) -- (4.75,2.165) -- (4.25,2.165) -- (3.75,2.165) -- (3.25,2.165) --
  (2.75,2.165) -- (2.5,1.732) -- (2.75,1.299) -- (3,0.866) -- (3.25,0.433) -- (3.75,0.433) --
  (4,0.866) -- (4.25,1.299) -- (4.75,1.299) -- (5.25,1.299) -- (5.75,1.299) -- (6,0.866) --
  (5.75,0.433) -- (5.5,0) -- (4.848,-0.348) -- (4.152,-0.348) -- (3.302,-0.231) -- (2.5,0) --
  (2.25,0.433) -- (2,0.866) -- (1.75,1.299) -- (1.5,1.732) -- (1.75,2.165) -- (2,2.598) --
  (1.75,3.031) -- (1.451,2.975) -- (1.035,2.11) -- (0.578,1.238) -- (0.25,0.433) -- (0.5,0) --
  (0.368,-0.094);
\draw[Lstar] (0.368,-0.094) .. controls (-0.9,-2.6) and (9.9,-2.6) .. (8.632,-0.094);
\node[fill=white,inner sep=1pt] at (4.5,7.217) {$D$};
\begin{scope}[x={(1,0)},y={(0.5,0.866)}]
\foreach \p in {(0,0),(3,0),(4,0),(5,0),(9,0),(2,1),(4,1),(5,1),(7,1),(1,2),(6,2),
  (7,2),(1,3),(2,3),(3,3),(6,3),(0,4),(2,4),(0,5),(1,5),(2,5),(3,5),(1,6),(0,7),(1,7),
  (2,7),(0,8),(1,8),(0,9)} {\plus{\p}}
\foreach \p in {(1,0),(2,0),(6,0),(7,0),(8,0),(0,1),(1,1),(3,1),(6,1),(8,1),(0,2),
  (2,2),(3,2),(4,2),(5,2),(0,3),(4,3),(5,3),(1,4),(3,4),(4,4),(5,4),(4,5),(0,6),(2,6),
  (3,6)} {\minus{\p}}
\end{scope}
\fill (7.746,4.472) circle (2.3pt);
\fill (1.254,4.472) circle (2.3pt);
\filldraw[fill=white,line width=0.6pt] (4.5,-1.15) circle (2.3pt);
\node[below left=0pt] at (0,0) {$x_3$};
\node[below right=0pt] at (9,0) {$x_1$};
\node[above=1pt] at (4.5,7.794) {$x_2$};
\node[below=2pt] at (4.5,-1.15) {$s_3$};
\node[right=2pt] at (7.746,4.472) {$s_1$};
\node[left=2pt] at (1.254,4.472) {$s_2$};
\node at (4.5,-2.45) {$\sface$};
\node at (6.013,6.174) {$A_1$};
\node at (2.987,6.174) {$A_2$};
\node at (1.8,-0.55) {$A_3$};
\node at (1,-2.05) {\color{red!80!black}$\gamma$};
\end{tikzpicture}
\caption{The interface $\gamma$ for a configuration of spins on a triangular domain (black: plus;
white: minus; open bonds are not shown). The arcs $A_1,A_2,A_3$ join the marked points
$x_1,x_2,x_3$, the ghost vertices $s_1,s_2,s_3$ are joined to every site of $A_1$, $A_2$, and $A_3$,
respectively (lines), and the face ${\sface}=s_1s_2s_3$ is the outer face. Here
$\sigma_{s_1}=\sigma_{s_2}=+$ and $\sigma_{s_3}=-$, and $\gamma$ (thick red) is the only interface
through ${\sface}$. It is a closed curve: outside ${\sface}$, it starts by crossing an edge at $x_3$
and ends by crossing an edge at $x_1$, and between these crossings it runs through $\Omega_\delta$
and, in places, along the arcs, between an arc and its ghost vertex. It separates the plus spins
joined to $A_1$ and $A_2$ by plus paths from the minus spins joined to $A_3$ by minus paths. The
thin red curves are other interfaces; they do not enter the definition of $D$, the component
attached to $x_2$.}\label{fig:interface}
\end{figure}

\begin{figure}[tp]
\centering
\begin{tikzpicture}[scale=1,
  bd/.style={line width=0.45pt},
  il/.style={line width=1pt,red!80!black,line join=round}]
\def\domain{plot[smooth cycle,tension=0.62] coordinates {(4,0.3) (4.15,1.55) (4.5,2.35) (3.85,2.4)
  (3.1,1.95) (2.35,2.4) (1.45,2.8) (0.45,2.35) (0.05,1.1) (0.35,-0.25) (1.25,-1) (2.55,-1.3)
  (3.6,-0.85)}}
\draw[bd] \domain;
\draw[il] plot[smooth,tension=0.55] coordinates {(1.25,-1) (0.95,-0.45) (0.52,0.3) (0.07,1.1)
  (0.7,1.02) (1.4,0.6) (2.1,0.82) (2.75,1.42) (3.1,1.93) (3.45,1.35) (3.75,0.8) (4,0.3)};
\fill (0.06,1.1) circle (1.3pt);
\fill (3.1,1.94) circle (1.3pt);
\fill (4,0.3) circle (1.7pt) node[right=1pt] {$x_1$};
\fill (0.45,2.35) circle (1.7pt) node[above left=-1pt] {$x_2$};
\fill (1.25,-1) circle (1.7pt) node[below left=-1pt] {$x_3$};
\node at (2.4,2.75) {$A_1$};
\node at (-0.3,0.3) {$A_2$};
\node at (2.6,-1.62) {$A_3$};
\node at (2.6,0.72) {\color{red!80!black}$\gamma$};
\node at (1.25,1.8) {$D$};
\fill (2.1,1.95) circle (1.5pt) node[right=1pt] {$z$};
\end{tikzpicture}
\caption{A large-scale picture of an interface $\gamma$ (red) as in Figure~\ref{fig:interface}, when
$\sigma_{s_1}=\sigma_{s_2}=+$ and $\sigma_{s_3}=-$: its passage through ${\sface}$ is not shown, and
its parts along the arcs lie on the boundary, so that it runs in $\Omega$ from $x_3$ to $x_1$. Its
contacts with the arcs $A_2$ and $A_1$ (dots) cut the region between $\gamma$ and these arcs into
several pieces. The component $D$ is the piece attached to $x_2$. Roughly speaking, $z$ lies in $D$
if a simple plus path joining $A_1$ to $A_2$ separates $z$ from $A_3$. Here $W=w_2$, and $z$
receives the contribution $4w_2$ in~\eqref{eq:Fdelta}.}\label{fig:observable}
\end{figure}
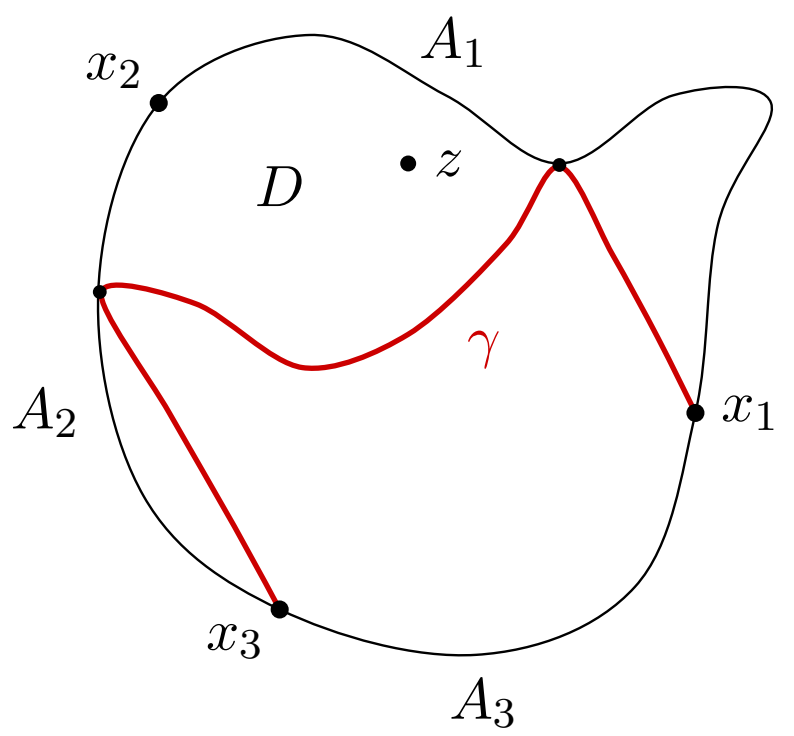

If $W=w_k$, a single interface $\gamma$ passes through ${\sface}$; it separates $s_{k+1}$ from
$s_{k-1}$ and $s_k$. No boundary condition is needed for this, since $\gamma$ exists as soon as the
three ghost spins are not all equal, or equivalently $W\neq0$, and the case $W=0$ does not
contribute to the observable below, which contains the factor $W$. The interface $\gamma$ is a
closed curve, which passes through ${\sface}$ outside the domain (Figure~\ref{fig:interface}).
Outside ${\sface}$, it starts by crossing an edge at $x_{k+1}$ and ends by crossing an edge at
$x_{k-1}$, the endpoints of $A_{k+1}$, and between these crossings it runs through $\Omega_\delta$
and, in places, along the arcs, between an arc and its ghost vertex. Let $G_\omega$ be the graph
obtained from $\delta\cH$, which now also has a vertex for each of the triangles formed with the
ghost vertices (Section~\ref{sec:averaged}), by removing the edges $e^*$ of the open bonds $e$ and
keeping all its vertices; since contracting an edge $e$ deletes $e^*$, the graph $G_\omega$ plays
the role of the dual of the contracted graph. Let $D$ be the component of $G_\omega\setminus\gamma$
that contains the edge $(s_{k-1}s_k)^*$, where we remove the points of $\gamma$, its vertices
included, but not the other edges at these vertices. It is adjacent to $x_k$
(Figures~\ref{fig:interface} and~\ref{fig:observable}). If $W=0$, we set $D=\varnothing$. For a
point $z\in\Omega$, we write $z\in D$ if the face of $\delta\T$ that contains $z$, viewed as a
vertex of $\delta\cH$, is a vertex of $D$, or if $z$ lies in a hole of $D$, that is, in a region
that is enclosed by a circuit of open bonds and such that the faces outside this circuit that have a
side on it are all vertices of $D$. We also regard $D$ as the planar region formed by these points,
namely the union of the closed faces of $\delta\T$ that are vertices of $D$, together with the holes
of $D$. For a point $z\in \Omega$ and an edge $e$, let
\begin{align}
F^\omega(e)&=4\,E^\omega\bigl[W\,\one_{\{e^*\in D\}}\bigr],\label{eq:observable}\\
F_\delta(z)&=4\,\Eb\bigl[W\,\one_{\{z\in D\}}\bigr]=\sum_{k=1}^3w_k\,\Pb(z\in D\mid W=w_k).\label{eq:Fdelta}
\end{align}
If $e$ is not an open bond, it is an edge of the contracted graph, and $e^*$ is an edge of its dual
$G_\omega$; if $e$ is an open bond, $e^*$ is not in $G_\omega$, and $F^\omega(e)=0$. The factor $4$
compensates for the probability $\frac14$ of each event $\{W=w_k\}$. Thus $F_\delta$ is the analogue
of~\eqref{eq:H}, to which it is equal at $\beta=0$, and $F^\omega$ is its quenched version, in which
the probabilities are computed given $\omega$. More precisely,
$F_\delta(z)=\Eb[F^\omega_\delta(z)]$, where $F^\omega_\delta(z)=4\,E^\omega[W\,\one_{\{z\in D\}}]$,
and $F^\omega_\delta(z)$ is close to the values of $F^\omega$ at the edges near $z$. To see this,
fix $\omega$, and let $z$ lie in a face $\fce$ of $\delta\T$ that is not enclosed by a circuit of
open bonds, so that $z\in D$ if and only if $\fce$ is a vertex of $D$, and some side of $\fce$ is
not an open bond. If $\gamma$ does not cross $\fce$, every edge $e_j(\fce)^*$ of $G_\omega$ lies in
$D$ if and only if $\fce$ does, and, if $\gamma$ crosses $\fce$, then $\fce$ is not a vertex of $D$.
Hence, for every side $e_j(\fce)$ that is not an open bond,
\begin{equation*}
F^\omega_\delta(z)=F^\omega\bigl(e_j(\fce)\bigr)-4\,E^\omega\bigl[W\,\one_{\{\gamma\text{ crosses }\fce,\ e_j(\fce)^*\in
D\}}\bigr],
\end{equation*}
where the last term is bounded by $4P^\omega(\gamma\text{ crosses }\fce)$. If $\fce$ is enclosed by
a circuit of open bonds, then, similarly, $F^\omega_\delta(z)$ differs from $F^\omega(e)$, for every
side $e$, not an open bond, of one of the faces that lie outside the outermost such circuit and have
a side on it, by at most four times the probability that $\gamma$ crosses a face with a vertex in
the cluster of this circuit. Since the clusters are small, these error terms tend to zero in
expectation as $\delta\to0$ for the points $z$ at a positive distance from $\partial\Omega$, so that
$F_\delta(z)$ is, up to a small error, the average over $\omega$ of the values of $F^\omega$ at the
edges near $z$ that are not open bonds. This error plays no role in the holomorphicity of the limit;
near the boundary, it is controlled by the strong RSW estimates (Section~\ref{sec:cardy}). The ghost
vertices play a role similar to the disorders at the marked points in the loop representation of
percolation of Khristoforov and Smirnov~\cite{KS21}, where the interfaces starting at the marked
points are the paths of a uniformly random loop configuration.

Color switching on clusters gives an exact analogue of Smirnov's relation~\eqref{eq:smirnov} for
$F^\omega$. For every face $\fce\neq{\sface}$ with no open side, we have (Section~\ref{sec:face} and
Corollary~\ref{cor:face})
\begin{equation}\label{eq:face}
\sum_{j=0}^2\tau^j\,F^\omega\bigl(e_j(\fce)\bigr)=0,
\end{equation}
so that the discrete contour integral of $F^\omega$ around $\fce$, in $\delta\T$, vanishes. This
relation does not involve the map $\cB$, and it holds even if the image of $\fce$ under $\cB$ is
degenerate, as happens when two vertices of $\fce$ lie in one cluster although the side that joins
them is closed. At a face $\fce$ with exactly one open side, the relation degenerates. We label the
sides so that $e_0=e_0(\fce)$ is the open one, and we write $e_1=e_1(\fce)$ and $e_2=e_2(\fce)$.
First, $F^\omega(e_0)=0$, since $e_0^*$ has been removed from $G_\omega$ and is never in $D$.
Second, $F^\omega(e_1)=F^\omega(e_2)$, because the vertex $\fce$ of $\delta\cH$ has only two edges
in $G_\omega$, namely $e_1^*$ and $e_2^*$, so that either $\gamma$ passes through $\fce$ and
contains both of them, or $e_1^*$ and $e_2^*$ are joined through $\fce$ in
$G_\omega\setminus\gamma$; in both cases, either both lie in $D$ or both lie outside $D$. This
agrees with the picture on the contracted graph, on which such a face has two sides with the same
endpoints. Since $\cB$ maps $e_0$ to a point, $d_1(\fce)=-d_2(\fce)$, so that the image of $\fce$ is
a segment, or a point if its three vertices lie in one cluster, and the contour integral
$F^\omega(e_1)\,d_1(\fce)+F^\omega(e_2)\,d_2(\fce)$ of $F^\omega$ around this image vanishes. Thus
$F^\omega$ satisfies Smirnov's relations exactly, but on the contracted graph, not on $\delta\T$.

Around the image of a face under $\cB$, however, the discrete contour integral of $F^\omega$ does
not vanish. If no side of $\fce$ is open, the face relation~\eqref{eq:face} gives
\begin{equation}\label{eq:defect}
\sum_{j=0}^2F^\omega\bigl(e_j(\fce)\bigr)\,d_j(\fce)=q_{\fce}\sum_{j=0}^2\bar\tau^j\,F^\omega\bigl(e_j(\fce)\bigr),
\end{equation}
that is, rotations and dilations preserve Smirnov's relation, as the maps $z\mapsto Az$ preserve the
Cauchy--Riemann equations, and only shears create defects. Local defects of analyticity created by
the embedding were studied by Beffara~\cite{Bef08}. Here the defects cannot be bounded face by face.
To count them, we say that a face borders $D$ if $\gamma$ crosses it and the dual edge of the side
that $\gamma$ does not cross lies in $D$, and we assume, for the sake of exposition, that
$r^{d+o(1)}$ faces border $D$ along a piece of $\gamma$ of diameter $r\delta$, for some $d\in(1,2)$.
A face borders $D$ at the scale $r\delta$ when three disjoint paths join it to distance $r\delta$,
two of them of the color on the side of $D$ and one of the other color. For percolation, this
three-arm event has probability $r^{-2/3+o(1)}$~\cite{SW01}, so that $d=\frac43$. The interface
$\gamma$ itself crosses more faces, $r^{7/4+o(1)}$ for percolation, as determined by the two-arm
event, which has probability $r^{-1/4+o(1)}$, since it also crosses the faces of the pockets cut off
from $D$. RSW estimates give upper bounds of the form $r^{2-c}$, for some $c>0$, which suffice for
all the upper bounds in our strategy below. In~\eqref{eq:defect}, the sum on the right-hand side is
bounded by $4P^\omega(\fce\text{ borders }D)$, which is of order $\delta^{2-d+o(1)}$, and
$|q_{\fce}|\asymp\delta$ at a positive fraction of the faces. Summing these bounds over the
$\asymp\delta^{-2}$ faces of a macroscopic region gives
$\delta^{-2}\cdot\delta\cdot\delta^{2-d+o(1)}=\delta^{1-d+o(1)}$, which tends to infinity since
$d>1$, so that the sum of the defects can be shown to be small only through cancellations.

\subsection{Reduction to a single interface}\label{sec:ideas-interface}

By Weyl's lemma~\cite{Wey40}, a locally integrable function $F$ is almost everywhere equal to a
holomorphic function if and only if $\iint\dbar f\,F\,dA=0$ for every smooth test function $f$ with
compact support in $\Omega$. Unlike Morera's theorem, used in Section~\ref{sec:smirnov-beta0}, it
requires no continuity of $F$, and it is formulated with smooth test functions, which suit the
homogenization of Section~\ref{sec:ideas-selfavg}. As a function of the point $z\in\Omega$, the
observable $F_\delta(z)=4\Eb[W\one_{\{z\in D\}}]$ is an expectation of indicator functions of $D$,
with its holes filled. Applying Green's formula to $D$ in each configuration and taking
expectations, we obtain
\begin{equation}\label{eq:green}
\tfrac{\ii}{2}\iint\dbar f\,F_\delta\,dA=\Eb\Bigl[W\oint_{\partial D}f\,dz\Bigr].
\end{equation}
The observable enters the right-hand side through the random region $D$ and the weight $W$, whereas
$f$ is a fixed test function. Thus the limits of $F_\delta$ are holomorphic if and only if these
expected contour integrals tend to zero for every test function $f$. Since $f$ vanishes near
$\partial\Omega$, only the part of $\partial D$ along $\gamma$ matters, together with the open bonds
on $\partial D$, which $\cB$ maps to points. A face $\fce$ crossed by $\gamma$ has exactly one side
$e_{j(\fce)}(\fce)$ that $\gamma$ does not cross, and, as in Section~\ref{sec:ideas-switching},
$\fce$ borders $D$ if $e_{j(\fce)}(\fce)^*$ lies in $D$; this never happens if the side is open.
These sides form the boundary of $D$ along $\gamma$, and the clockwise orientation around $\fce$ is
the positive orientation of $\partial D$. Replacing them by their images under $\cB$ changes the
contour integral by $O(N\delta\ell)=\delta^{2-d+o(1)}$, which tends to zero, where
$N=\delta^{-d+o(1)}$ is the number of faces that border $D$. Hence
$\Eb[W\,V(f)]\approx(\ii/2)\iint\dbar f\,F_\delta\,dA$, where
\begin{equation}\label{eq:fields}
\begin{aligned}
V(f)&=\sum_{\fce}f(z_{\fce})\,d_{j(\fce)}(\fce)=U(f)+Q(f),\\
U(f)&=\sum_{\fce}f(z_{\fce})\,p_{\fce}\,\tau^{j(\fce)},\\
Q(f)&=\sum_{\fce}f(z_{\fce})\,q_{\fce}\,\bar\tau^{j(\fce)},
\end{aligned}
\end{equation}
the sums run over the faces $\fce$ that border $D$, and $z_{\fce}$ is the image of the center of
$\fce$. The rotation and dilation part does not contribute at all, for two reasons. First, for a
face $\fce$ with no open side,
\begin{equation*}
\sum_{j=0}^2\tau^j\,\one_{\{e_j(\fce)^*\in D\}}=\tau^{j(\fce)}\,\one_{\{\fce\text{ borders }D\}},
\end{equation*}
since the three indicators are equal if $\gamma$ does not cross $\fce$, and $\sum_j\tau^j=0$,
whereas only $e_{j(\fce)}(\fce)^*$ can lie in $D$ if $\gamma$ crosses $\fce$. A face with an open
side, on the other hand, never borders $D$. If $\gamma$ crosses such a face, it leaves the open side
$e$ uncrossed, since interfaces never cross open bonds, and $e^*$ is not in $G_\omega$. Such faces
therefore do not enter $U(f)$, $V(f)$, and $Q(f)$ at all, and the failure of the face
relation~\eqref{eq:face} at them does not matter. No information is lost, since $\cB$ maps their
open side to a point, which does not contribute to the contour integral along the image of
$\partial D$. Second, given $\omega$, the numbers $f(z_{\fce})$ and $p_{\fce}$ are fixed. The
coefficients $p_{\fce}$ depend on $\omega$, so that they can be taken out of the expectation, and
the face relation applied, only for a fixed $\omega$, which is why we work with the quenched
observable $F^\omega$. Therefore, with sums over the faces with no open side,
\begin{equation}\label{eq:VQ}
\begin{aligned}
U(f)&=\sum_{\fce}f(z_{\fce})\,p_{\fce}\sum_{j=0}^2\tau^j\,\one_{\{e_j(\fce)^*\in D\}},\\
E^\omega\bigl[W\,U(f)\bigr]&=\frac14\sum_{\fce}f(z_{\fce})\,p_{\fce}\sum_{j=0}^2\tau^j\,F^\omega\bigl(e_j(\fce)\bigr)=0,
\end{aligned}
\end{equation}
by the definition~\eqref{eq:observable} of $F^\omega$ and the face relation~\eqref{eq:face}.
Consequently, $\Eb[W\,V(f)]=\Eb[W\,Q(f)]$, so that only the shears matter. At $\beta=0$, there are
no shears, $Q=0$, and we recover Smirnov's argument.

\subsection{Homogenization of the shears}\label{sec:ideas-selfavg}
The key heuristic step is that, along macroscopic interfaces,
\begin{equation*}
\begin{aligned}
Q&\approx\eta\,V,\\
U&=V-Q\approx(1-\eta)\,V,
\end{aligned}
\end{equation*}
for a real constant $\eta$ that depends only on $\beta$. In other words, the shears carry the part
$\eta V$ of the contour integral, and the rotations and dilations carry the part $(1-\eta)V$. By our
convention, this means that $Q(f)-\eta V(f)$ should tend to zero in probability, for every test
function $f$; for the held interfaces below, it will tend to zero in mean square. We will explain
this relation, and define $\eta$, in the following steps. Throughout, $\delta$ tends to zero first,
and the scales of the pieces below tend to zero only at the end. Many relations below hold only in
mean square, in the sense that the expectation of the square of the error is small; we will state
this for each of them, together with the law under which the expectation is taken. The other
relations either hold for every configuration or relate deterministic quantities, such as $m(R)$,
$g$, and $\lambda$ below.

\emph{Pieces.} Along $\gamma$, the sides $e_{j(\fce)}(\fce)$ of the faces that border $D$ form a
path in the contracted graph, from cluster to cluster, and their images under $\cB$, the vectors
$d_{j(\fce)}(\fce)$, follow each other. Cut this path at some of its vertices into pieces on which
$f$ is nearly constant, and let $\zeta_0,\zeta_1,\dots$ be the images of these vertices under $\cB$,
namely the barycenters of the corresponding clusters. The increments of $V(1)$ and $Q(1)$, with
$f=1$, along the $i$-th piece are
\begin{equation*}
\begin{aligned}
\Delta V_i&=\zeta_{i+1}-\zeta_i,\\
\Delta Q_i&=\sum_{\fce\text{ in the $i$-th piece}}q_{\fce}\,\bar\tau^{j(\fce)} .
\end{aligned}
\end{equation*}
The first identity holds for every configuration. The second increment is a sum of many small random
terms, which does not telescope.

\emph{Closed interfaces.} Recall that the interfaces are circuits of $\delta\cH$
(Section~\ref{sec:ES}); $\gamma$ is the interface through ${\sface}$, and we call the other ones
closed interfaces. For every interface $L$ and every component $C$ of $G_\omega\setminus L$, we
define $V_{L,C}$, $U_{L,C}$, and $Q_{L,C}$ as $V$, $U$, and $Q$ in~\eqref{eq:fields}, with
$(\gamma,D)$ replaced by $(L,C)$. Around a closed interface $L$, $V_{L,C}(1)=0$ for every
configuration, since the image of a closed curve is closed. The shear part $Q_{L,C}(1)$ should be
small as well, in mean square under $\Pb$, compared with the diameter $R$ of $L$; for comparison,
the contour integral along a piece of $L$ of diameter $R$ is of order $R$. Here and below, lengths
such as $R$ are Euclidean lengths in the plane, in which the lattice $\delta\T$ has spacing
$\delta$, and $o(R)$ denotes a quantity whose ratio to $R$ tends to zero as $\delta/R\to0$, for
instance as $\delta\to0$ with $R$ fixed. More precisely, let $m(R)$ be the expected sum of
$|Q_{L,C}(1)|^2$ over the pairs $(L,C)$ in which $L$ is a closed interface of diameter between $R$
and $2R$ in a fixed region away from $\partial\Omega$, divided by the expected number of these
interfaces $L$. Then
\begin{equation*}
m(R)=o(R^2)\qquad(\delta/R\to0).
\end{equation*}
This will follow from the Bessel identities of Section~\ref{sec:ideas-energy}.

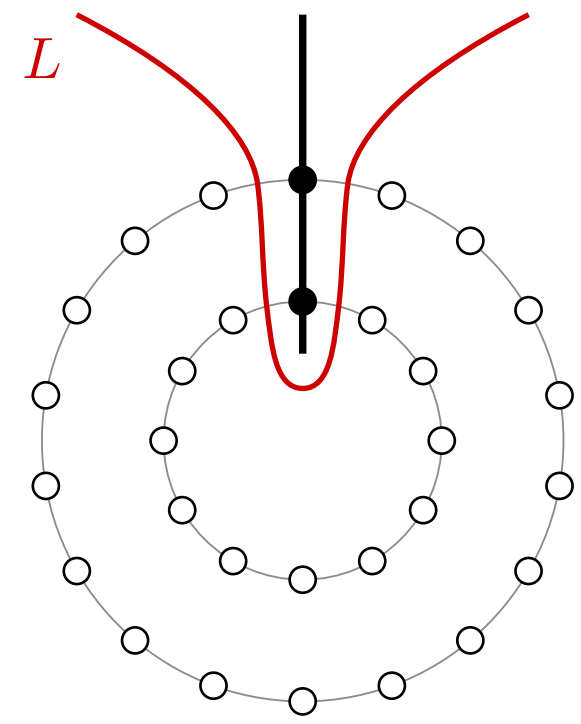
\begin{figure}[t]
\centering
\begin{tikzpicture}[scale=1.15,
  circ/.style={line width=0.35pt,black!45},
  site/.style={line width=0.5pt},
  pl/.style={line width=1.4pt},
  intf/.style={line width=1pt,red!80!black,line join=round}]
\draw[circ] (0,0) circle (1.5);
\draw[circ] (0,0) circle (0.8);
\draw[pl] (0,2.45) -- (0,0.5);
\draw[intf] (-1.3,2.45) .. controls (-0.62,2.1) and (-0.3,1.755) .. (-0.26,1.48)
  .. controls (-0.229,1.255) and (-0.2355,1.0105) .. (-0.207,0.773)
  .. controls (-0.18,0.548) and (-0.16,0.3) .. (0,0.3)
  .. controls (0.16,0.3) and (0.18,0.548) .. (0.207,0.773)
  .. controls (0.2355,1.0105) and (0.229,1.255) .. (0.26,1.48)
  .. controls (0.3,1.755) and (0.62,2.1) .. (1.3,2.45);
\foreach \t in {10,30,...,350} {
  \ifnum\t=90
    \filldraw[site,fill=black] (\t:1.5) circle (0.075);
  \else
    \filldraw[site,fill=white] (\t:1.5) circle (0.075);
  \fi}
\foreach \t in {0,30,...,330} {
  \ifnum\t=90
    \filldraw[site,fill=black] (\t:0.8) circle (0.075);
  \else
    \filldraw[site,fill=white] (\t:0.8) circle (0.075);
  \fi}
\node[red!80!black] at (-1.5,2.2) {$L$};
\end{tikzpicture}
\caption{Holding an interface in a small disk. The configuration in the disk contains two circuits
of sites around the center (thin lines), each consisting of minus sites (white) except for one plus
site (black), and a plus path (thick) that reaches the center through these plus sites. An interface
can cross a circuit only at an edge that joins a plus site to a minus site. Hence the interface $L$
(red) that separates the plus path from the minus sites crosses each circuit next to its plus site,
at places prescribed by the configuration in the disk.}\label{fig:holding}
\end{figure}

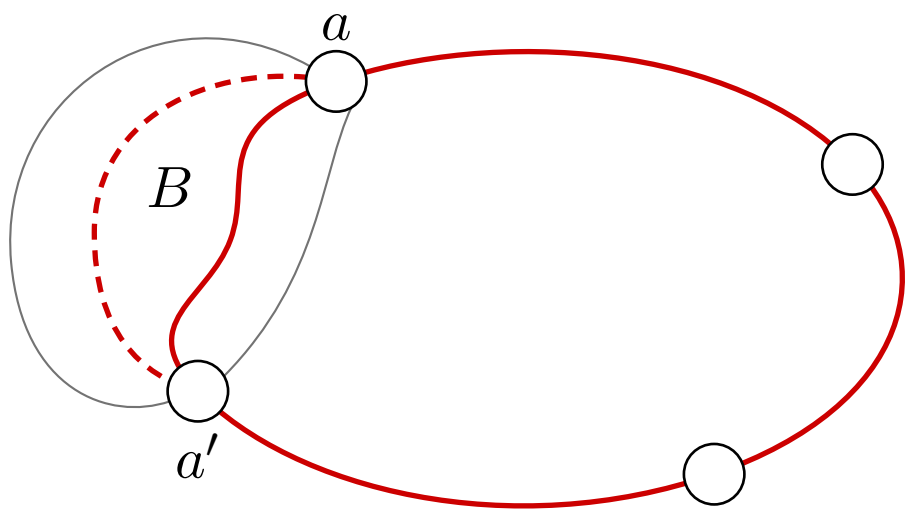
\begin{figure}[t]
\centering
\begin{tikzpicture}[scale=1,
  intf/.style={line width=1pt,red!80!black,line join=round},
  alt/.style={line width=1pt,red!80!black,dash pattern=on 3pt off 2pt,line join=round},
  bnd/.style={line width=0.4pt,black!55}]
\draw[bnd] (-1.15,1.126) .. controls (-1.35,0.7) and (-1.35,0) .. (-1.992,-0.65);
\draw[bnd] (-1.423,1.399) .. controls (-2.4,1.95) and (-3.5,1.2) .. (-3.4,0.1)
  .. controls (-3.32,-0.75) and (-2.75,-0.95) .. (-2.353,-0.818);
\draw[intf] (2.165,0.75) .. controls (1.39,1.555) and (-0.211,1.659) .. (-1.25,1.299)
  .. controls (-2.116,0.999) and (-1.8,0.65) .. (-1.95,0.25)
  .. controls (-2.1,-0.15) and (-2.6,-0.295) .. (-2.165,-0.75)
  .. controls (-1.54,-1.4) and (-0.092,-1.764) .. (1.25,-1.299)
  .. controls (2.549,-0.849) and (2.79,0.1) .. cycle;
\draw[alt] (-1.25,1.299) .. controls (-2,1.45) and (-2.85,1.1) .. (-2.85,0.3)
  .. controls (-2.85,-0.35) and (-2.55,-0.65) .. (-2.165,-0.75);
\foreach \p in {(2.165,0.75),(-1.25,1.299),(-2.165,-0.75),(1.25,-1.299)}
  \filldraw[fill=white,line width=0.5pt] \p circle (0.2);
\node[above=4pt] at (-1.25,1.299) {$a$};
\node[below=4pt] at (-2.165,-0.75) {$a'$};
\node at (-2.35,0.6) {$B$};
\end{tikzpicture}
\caption{A closed interface (red) held at four small disks, in which the configuration forces it to
pass at prescribed places; $a$ and $a'$ are the places where it passes through two consecutive
disks. The band $B$ between these two disks is the region between its two boundary curves (thin
lines); it contains a plus path and a minus path (not shown), which keep the interface away from
these curves, except near the disks. Resampling the configuration in the band replaces the piece of
the interface between $a$ and $a'$ by another piece (dashed), and leaves the rest unchanged, except
near $a$ and $a'$. Since the shear part around the closed interface is small before and after the
resampling, the two pieces carry nearly the same shear increment.}\label{fig:routes}
\end{figure}

\emph{Path independence.} We now hold a closed interface at a few places, as in
Figure~\ref{fig:routes}. To hold an interface in a small disk, we require that the configuration in
the disk contain two circuits of sites around the center, each consisting of minus sites except for
one or two adjacent plus sites, and a plus path that reaches the center through these plus sites
(Figure~\ref{fig:holding}). An interface can cross only edges that join sites of different spins,
and, along each circuit, there are only two such edges, at the ends of its plus sites. Hence the
interface that separates the plus path from the minus sites crosses each circuit through these two
edges, and it passes through the disk at places determined by the configuration in the disk
(Sections~\ref{sec:ideas-difficulty} and~\ref{sec:constrained}). Let $L$ be a closed interface of
diameter $R$ that is held in this way at a few small disks, of radius $h$, with $\delta\ll h\ll R$,
and let $a$ and $a'$ be the places where it passes through two consecutive disks; in
Figure~\ref{fig:routes}, $L$ is held at four disks, and two of them are marked $a$ and $a'$. Between
$a$ and $a'$, we also keep $L$ in a band $B$, a region shaped like a curved strip, which runs from
the disk at $a$ to the disk at $a'$ and is fixed in advance, independently of the configuration. The
boundary of $B$ consists of two short arcs, where $B$ meets the two disks, and of two curves that
run from the disk at $a$ to the disk at $a'$, one on each side of $B$; we call these two curves the
boundary curves of $B$ (thin lines in Figure~\ref{fig:routes}). Away from the disks, the width of
$B$, defined as the distance between its boundary curves, is comparable to the distance between $a$
and $a'$. We require that the band contain a plus path and a minus path from one disk to the other,
both at a distance comparable to the width of the band from the boundary curves, except near the
disks. The piece of $L$ between $a$ and $a'$, which separates these two paths, then also stays away
from the boundary curves. All the expectations in this step are conditional on these requirements.
Let $\Pb^{\mathrm{held}}$ be the law of the configuration conditioned on them, and
$\Eb^{\mathrm{held}}$ its expectation. Since the number of requirements is bounded, they hold with
probability bounded below, by the strong RSW estimates, so that
\begin{equation}\label{eq:held}
\frac{d\Pb^{\mathrm{held}}}{d\Pb}\le C .
\end{equation}
Now resample the configuration in the band, including its two paths, according to its conditional
law under $\Pb^{\mathrm{held}}$ given the outside; this preserves $\Pb^{\mathrm{held}}$. Let
$(L',C')$ be the new pair. In Figure~\ref{fig:routes}, the piece of $L$ between $a$ and $a'$ (solid)
is replaced by the dashed piece, and the rest of $L$ does not move. We will use three facts:
\begin{itemize}[leftmargin=2em]
\item Splitting. $Q_{L,C}(1)$ is a sum over the faces of $L$ that border $C$. We split it into the
contribution $\Delta Q$ of the faces in the band, which form the piece of $L$ between $a$ and $a'$,
and the contribution $Q_{\mathrm{rest}}$ of the other faces.
\item The rest does not change. The configuration outside the band is not resampled, and the new
interface still passes through $a$ and $a'$. Thus $L'$ coincides with $L$ outside the band, and only
the piece between $a$ and $a'$ is replaced, with a new contribution $\Delta Q'$; $Q_{\mathrm{rest}}$
changes only near $a$ and $a'$, through the clusters that meet the band, and we neglect this change.
\item Both shear parts, $Q_{L,C}(1)$ and $Q_{L',C'}(1)$, are $o(R)$ in mean square under
$\Pb^{\mathrm{held}}$. Under $\Pb$, the bound $m(R)=o(R^2)$ above says that $\sum_C|Q_{L,C}(1)|^2$,
and hence $|Q_{L,C}(1)|^2$ for each component $C$, has expectation $o(R^2)$ per closed interface of
diameter $R$. By~\eqref{eq:held}, since $|Q_{L,C}(1)|^2\ge0$, and since, as we will see in
Section~\ref{sec:ideas-energy}, the expected number of closed interfaces of diameter between $R$ and
$2R$ that meet a given disk of radius $R$ is bounded, we get
$\Eb^{\mathrm{held}}|Q_{L,C}(1)|^2=o(R^2)$ for the held interface. Since the resampling preserves
$\Pb^{\mathrm{held}}$, the pair $(L',C')$ has the same law as $(L,C)$, and
$\Eb^{\mathrm{held}}|Q_{L',C'}(1)|^2=o(R^2)$.
\end{itemize}
Thus, for every configuration, up to the neglected change near $a$ and $a'$,
\begin{equation*}
\begin{aligned}
Q_{L,C}(1)&=\Delta Q+Q_{\mathrm{rest}},\\
Q_{L',C'}(1)&=\Delta Q'+Q_{\mathrm{rest}},
\end{aligned}
\end{equation*}
where both left-hand sides are $o(R)$ in mean square under $\Pb^{\mathrm{held}}$. Subtracting the
two lines, we obtain
\begin{equation*}
\Eb^{\mathrm{held}}|\Delta Q-\Delta Q'|^2\le2\,\Eb^{\mathrm{held}}|Q_{L,C}(1)|^2+2\,\Eb^{\mathrm{held}}|Q_{L',C'}(1)|^2=o(R^2).
\end{equation*}
Now let $\langle\Delta Q\rangle=\Eb^{\mathrm{held}}[\Delta Q\mid\text{outside}]$ be the conditional
expectation of $\Delta Q$ given the configuration outside the band. Given the outside, $\Delta Q$
and $\Delta Q'$ are independent and have the same law, since $\Delta Q'$ comes from a fresh
configuration in the band, drawn from the same conditional law. Therefore the cross term has zero
expectation when
$\Delta Q-\Delta Q'=(\Delta Q-\langle\Delta Q\rangle)-(\Delta Q'-\langle\Delta Q\rangle)$ is
squared, and
\begin{equation}\label{eq:concentration}
\Eb^{\mathrm{held}}\bigl|\Delta Q-\langle\Delta Q\rangle\bigr|^2=\tfrac12\,\Eb^{\mathrm{held}}|\Delta Q-\Delta Q'|^2=o(R^2).
\end{equation}
Thus $\Delta Q$ concentrates around $\langle\Delta Q\rangle$, a function of the configuration
outside the band, with fluctuations $o(R)$ in mean square under $\Pb^{\mathrm{held}}$. Apart from
the changes near $a$ and $a'$ neglected above, the argument so far uses only exact identities, the
strong RSW estimates, and the bound $m(R)=o(R^2)$.

It remains to see that $\langle\Delta Q\rangle$ depends only on the configurations near the ends of
the piece. Two further steps are required. Locality uses spatial mixing, and the independence of the
band uses once more the smallness of the shear part around closed interfaces.

\emph{Locality.} The configuration outside the band acts on the configuration in the band only
across the boundary of the band, namely across the two arcs, through the configurations $\chi_a$ and
$\chi_{a'}$ near $a$ and $a'$, and across the boundary curves. Away from the disks, the two paths,
the piece of $L$ between them, and the clusters that meet it stay at a distance comparable to the
width of the band, hence to $R$, from the boundary curves. Spatial mixing, which expresses that the
model has a finite correlation length for $\beta<\beta_c$, makes the influence of the configuration
across the boundary curves negligible at such distances. More precisely, changing the configuration
on a set $\Lambda$ changes the probability of every event that depends only on the configuration in
a set $\Lambda'$ by at most
\begin{equation}\label{eq:mixing}
C\,|\Lambda'|\,e^{-\mathrm{dist}(\Lambda',\Lambda)/(C\delta)},
\end{equation}
where $|\Lambda'|$ is the number of sites of $\Lambda'$. The bound is uniform in the domain and in
the configuration elsewhere, even next to $\Lambda'$, and it applies to connection events. This
follows from the exponential decay of the Fortuin--Kasteleyn connectivities~\cite{ABF87,DRT19},
which holds even with wired boundary conditions (Section~\ref{sec:ES}). Through the Edwards--Sokal
coupling, the configuration on $\Lambda$ acts on that in $\Lambda'$ only through the clusters that
join the two sets (see also~\cite{MOS94} for general two-dimensional spin systems). Let $\Lambda$ be
the outside of the band along the boundary curves, away from the disks, and $\Lambda'$ the part of
the band at distance comparable to $R$ from the boundary curves, which contains the piece of $L$
away from the disks. Then $\Delta Q$, a sum of $O((R/\delta)^2)$ terms of order at most $\ell$,
depends only on the configuration in $\Lambda'$ and near the disks, and the bound~\eqref{eq:mixing}
is $O((R/\delta)^2e^{-cR/\delta})$ for some $c>0$. Since the requirements on the band have
probability bounded below, changing the configuration in $\Lambda$ therefore changes
$\langle\Delta Q\rangle$ by $O(Re^{-cR/\delta})$, for a smaller $c>0$, which is negligible. Hence
$\langle\Delta Q\rangle$ depends on the outside of the band only through $\chi_a$ and $\chi_{a'}$,
and it is close to the conditional expectation of $\Delta Q$ given $\chi_a$ and $\chi_{a'}$ only.
More precisely, for every configuration outside the band,
\begin{equation}\label{eq:locality}
\langle\Delta Q\rangle=\varphi_B(\chi_a,\chi_{a'})+O(Re^{-cR/\delta}),\qquad
\varphi_B(\chi_a,\chi_{a'})=\Eb^{\mathrm{held}}\bigl[\Delta Q\bigm|\chi_a,\chi_{a'}\bigr].
\end{equation}
For the same reason, the deterministic function $\varphi_B$ is the same, up to $O(Re^{-cR/\delta})$,
for all the held laws in which an interface is held at $a$ and $a'$ and kept in $B$, whatever the
requirements outside $B$; it may depend on the band $B$.

\begin{figure}[t]
\centering
\begin{tikzpicture}[scale=0.9,
  intf/.style={line width=1pt,red!80!black,line join=round},
  bnd/.style={line width=0.4pt,black!55}]
\def\rightband{%
  \draw[bnd] (0.13,1.375) .. controls (0.55,1) and (0.62,0.45) .. (0.62,0)
    .. controls (0.62,-0.45) and (0.55,-1) .. (0.13,-1.375);
  \draw[bnd] (0.17,1.56) .. controls (1.2,1.5) and (1.75,0.75) .. (1.75,0)
    .. controls (1.75,-0.75) and (1.2,-1.5) .. (0.17,-1.56);}
\def\disks{%
  \foreach \p in {(0,1.5),(0,-1.5)} \filldraw[fill=white,line width=0.5pt] \p circle (0.18);
  \node[above=3pt] at (0,1.5) {$a$};
  \node[below=3pt] at (0,-1.5) {$a'$};
  \node at (2.05,0) {$B_0$};}
\begin{scope}
\rightband
\draw[bnd] (-0.13,1.375) .. controls (-0.55,1) and (-0.62,0.45) .. (-0.62,0)
  .. controls (-0.62,-0.45) and (-0.55,-1) .. (-0.13,-1.375);
\draw[bnd] (-0.17,1.56) .. controls (-1.2,1.5) and (-1.75,0.75) .. (-1.75,0)
  .. controls (-1.75,-0.75) and (-1.2,-1.5) .. (-0.17,-1.56);
\draw[intf] (0,1.5) .. controls (0.6,1.3) and (1.05,0.85) .. (1.1,0.3)
  .. controls (1.14,-0.15) and (1.3,-0.65) .. (0.95,-1.05)
  .. controls (0.7,-1.32) and (0.4,-1.35) .. (0,-1.5);
\draw[intf] (0,1.5) .. controls (-0.55,1.3) and (-1,1) .. (-1.12,0.45)
  .. controls (-1.22,-0.05) and (-1.1,-0.6) .. (-0.85,-1)
  .. controls (-0.6,-1.35) and (-0.35,-1.37) .. (0,-1.5);
\disks
\node at (-2.05,0) {$B_1$};
\node at (0,-2.65) {(a)};
\end{scope}
\begin{scope}[shift={(6.5,0)}]
\rightband
\draw[bnd] (-0.13,1.375) .. controls (-0.75,1.15) and (-1.3,0.75) .. (-1.3,0.1)
  .. controls (-1.3,-0.55) and (-0.75,-1.1) .. (-0.13,-1.375);
\draw[bnd] (-0.17,1.56) .. controls (-1.3,2) and (-2.65,1.35) .. (-2.65,0.2)
  .. controls (-2.65,-0.85) and (-1.4,-1.6) .. (-0.17,-1.56);
\draw[intf] (0,1.5) .. controls (0.65,1.3) and (0.95,0.95) .. (1.05,0.45)
  .. controls (1.15,-0.05) and (1.35,-0.55) .. (1,-1)
  .. controls (0.75,-1.3) and (0.4,-1.36) .. (0,-1.5);
\draw[intf] (0,1.5) .. controls (-0.9,1.5) and (-1.75,1.15) .. (-1.9,0.45)
  .. controls (-2,-0.05) and (-1.95,-0.55) .. (-1.55,-0.95)
  .. controls (-1.15,-1.32) and (-0.55,-1.38) .. (0,-1.5);
\disks
\node at (-2.35,1.65) {$B_2$};
\node at (0,-2.65) {(b)};
\end{scope}
\end{tikzpicture}
\caption{Independence of the band. (a) A closed interface $L_1$ (red), held at the disks at $a$ and
$a'$, and kept in the bands $B_1$ and $B_0$, away from their boundary curves (thin lines). (b) The
same construction, with a band $B_2$ in place of $B_1$. In (a) and (b), the configuration in $B_0$
has the same conditional law given its outside, so that the shear increment along $B_0$ is close to
the same function of the configurations in the two disks. Since the shear part around a closed
interface is small, the increments along $B_1$ and $B_2$ are then close to the same function as
well.}\label{fig:bands}
\end{figure}
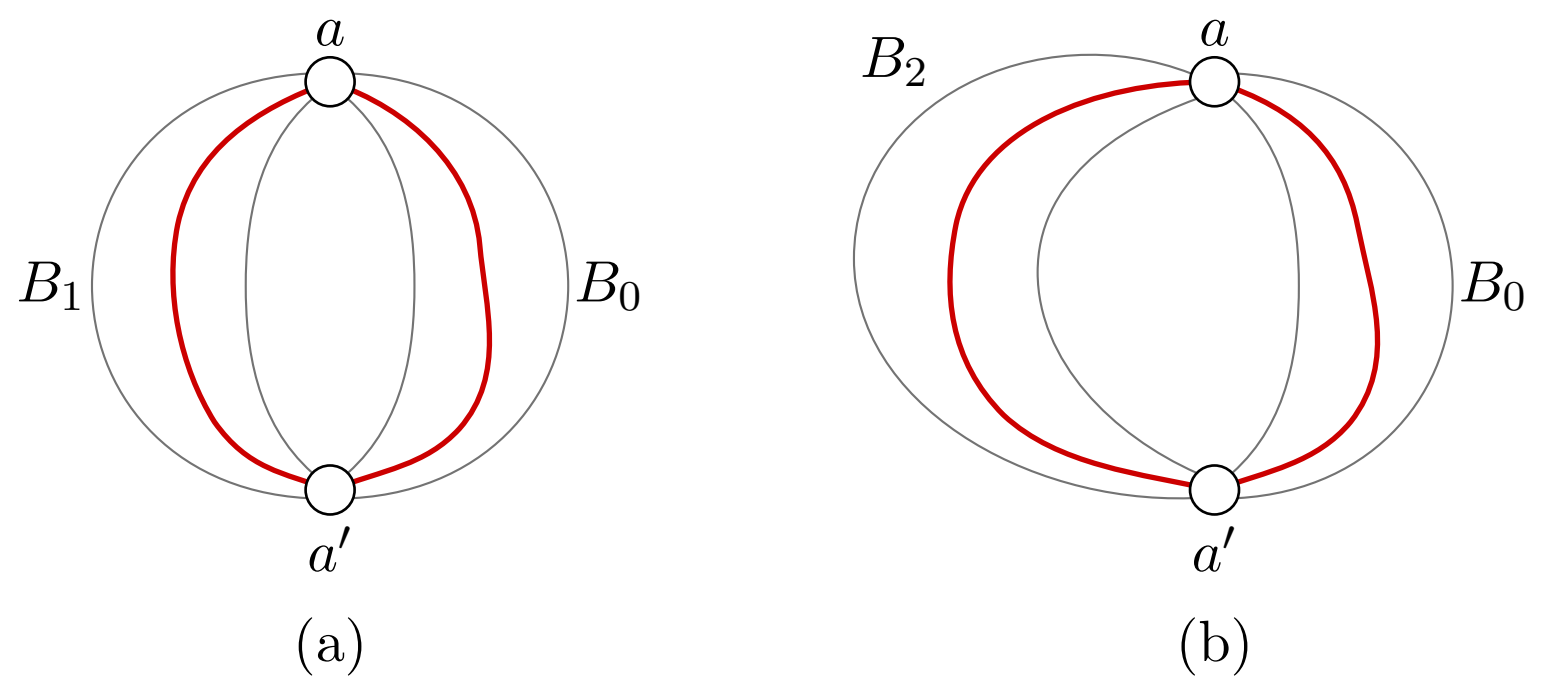

\emph{Independence of the band.} Let $B_1$ and $B_2$ be two bands from the disk at $a$ to the disk
at $a'$, and $B_0$ a band from $a'$ back to $a$ that meets them only at the two disks
(Figure~\ref{fig:bands}). For $i=1,2$, holding a closed interface $L_i$ at the two disks and keeping
it in $B_i$ and $B_0$, as above, gives a held law $\Pb^{\mathrm{held}}_i$. Let $\Delta Q_{B_i}$ and
$\Delta Q_{B_0}$ be the shear increments along the pieces of $L_i$ in $B_i$ and in $B_0$. The two
held laws differ only outside $B_0$, so that, under both of them, the configuration in $B_0$ has the
same conditional law given its outside. Hence~\eqref{eq:concentration} and~\eqref{eq:locality},
applied to each band, give, in mean square under $\Pb^{\mathrm{held}}_i$,
\begin{equation*}
\begin{aligned}
\Delta Q_{B_i}&=\varphi_{B_i}(\chi_a,\chi_{a'})+o(R),\\
\Delta Q_{B_0}&=\varphi_{B_0}(\chi_{a'},\chi_a)+o(R),
\end{aligned}
\end{equation*}
with the same function $\varphi_{B_0}$ for $i=1$ and $i=2$. Up to the parts of $L_i$ in the two
disks, which we neglect as above, the two increments add up to the shear part around $L_i$, for
every configuration, and this shear part is small in mean square under $\Pb^{\mathrm{held}}_i$:
\begin{equation*}
\Delta Q_{B_i}+\Delta Q_{B_0}=Q_{L_i,C_i}(1)=o(R).
\end{equation*}
The laws of $(\chi_a,\chi_{a'})$ under $\Pb^{\mathrm{held}}_1$ and $\Pb^{\mathrm{held}}_2$ have
bounded densities with respect to each other, by the strong RSW estimates. Therefore, in mean square
under either of the two held laws,
\begin{equation*}
\varphi_{B_1}(\chi_a,\chi_{a'})=-\varphi_{B_0}(\chi_{a'},\chi_a)+o(R)=\varphi_{B_2}(\chi_a,\chi_{a'})+o(R).
\end{equation*}
Thus the shear increment from $a$ to $a'$ is, up to $o(R)$ in mean square, the same deterministic
function of $\chi_a$ and $\chi_{a'}$, whichever band the piece follows; it depends on the ends of
the piece, and not on its route. This is how the shears average out. The increment is a sum of
$(R/\delta)^{d+o(1)}$ random terms of order $\delta$, one for each face of the piece that borders
$C$, and yet the configurations at the two ends of the piece determine it up to $o(R)$ in mean
square.

\emph{Linearity and the definition of $\eta$.} Two kinds of bounds appear in this step. Bounds on
random quantities, such as the increment $\Delta Q$ of a piece, hold only in mean square under the
held law, and we will state this each time. Bounds on deterministic quantities, such as the expected
increment $g$ and the function $\lambda$ below, are ordinary bounds; they come from the mean-square
ones, by $|\Eb X|\le(\Eb|X|^2)^{1/2}$.

We now compare the random increment $\Delta Q$ of a held piece with a deterministic function of its
displacement. To use translations, we work with the infinite-volume law. The law of the model in
$\Omega_\delta$ is not invariant under translations, rotations, or reflections, because of the
boundary of $\Omega_\delta$, but the infinite-volume law $\Pb_\beta$ of the configuration
$(\sigma,\omega)$ in the plane is invariant under the symmetries of $\delta\T$, since it is unique
for $\beta<\beta_c$. We may use it in place of the law in $\Omega_\delta$, because the pieces that
matter lie in the support of $f$, at a positive distance $c$ from $\partial\Omega$, and,
by~\eqref{eq:mixing}, the held laws in $\Omega_\delta$ and in the plane differ near such a piece by
$O((R/\delta)^2e^{-c/(C\delta)})$ in total variation, which is negligible.

For $v\in\C$ large compared with $h$, consider a held piece whose disks are centered at the site $0$
and at the point $v$, with some band from the first disk to the second one and some requirements
outside it; we will see that their choice does not matter. Let $\Eb^{\mathrm{held}}_{0,v}$ be the
expectation under $\Pb_\beta$ conditioned on these requirements, and let
\begin{equation*}
g(v)=\Eb^{\mathrm{held}}_{0,v}\bigl[\Delta Q\bigr],\qquad\Delta Q=\sum_{\fce}q_{\fce}\,\bar\tau^{j(\fce)},
\end{equation*}
where $C$ is the component of $G_\omega\setminus L$, for the held interface $L$, such that the
positive orientation of $\partial C$ runs along the piece from the disk at $0$ to the disk at $v$,
and the sum runs over the faces of the piece that border $C$, so that $g(v)$ is the expected shear
increment of the piece. Thus $g$ is deterministic, whereas $\Delta Q$ is random. For a held piece
whose disks are centered at $a$ and $a'$, they are related by
\begin{equation}\label{eq:g}
\begin{aligned}
\Delta Q&=\varphi_B(\chi_a,\chi_{a'})+o(R)&&\text{(concentration and locality)}\\
&=\Eb^{\mathrm{held}}[\Delta Q]+O(h)+o(R)&&\text{(the disks)}\\
&=g(a'-a)+O(h)+o(R)&&\text{(translation invariance)},
\end{aligned}
\end{equation}
where the error terms $O(h)$ and $o(R)$ are random, and their bounds hold in mean square under the
held law. Put differently, $\Delta Q$ is close to its conditional expectation given the
configurations in the two disks, this conditional expectation is close to the expectation of
$\Delta Q$, and this expectation is $g(a'-a)$, up to $O(h)+o(R)$. The three steps are the following.
\begin{itemize}[leftmargin=2em]
\item Concentration and locality. By~\eqref{eq:concentration} and~\eqref{eq:locality}, $\Delta Q$ is
close to its conditional expectation given the outside of the band, which depends on the outside
only through $\chi_a$ and $\chi_{a'}$.
\item The disks. This step says that the conditional expectation
$\varphi_B(\chi_a,\chi_{a'})=\Eb^{\mathrm{held}}[\Delta Q\mid\chi_a,\chi_{a'}]$, which is random, is
close to its average
$\Eb^{\mathrm{held}}[\varphi_B(\chi_a,\chi_{a'})]=\Eb^{\mathrm{held}}[\Delta Q]$, which is
deterministic. In other words, $\varphi_B$ changes only by $O(h)$, in mean square, when $\chi_a$ and
$\chi_{a'}$ change. To see this, split $\Delta Q=\Delta Q_{\mathrm{near}}+\Delta Q_{\mathrm{far}}$,
where $\Delta Q_{\mathrm{near}}$ is the contribution of the faces of the piece in the two disks, of
radius $h$, and $\Delta Q_{\mathrm{far}}$ that of the other faces of the piece. Then
\begin{equation*}
\begin{aligned}
\varphi_B(\chi_a,\chi_{a'})-\Eb^{\mathrm{held}}[\Delta Q]
&=\bigl(\Eb^{\mathrm{held}}[\Delta Q_{\mathrm{near}}\mid\chi_a,\chi_{a'}]-\Eb^{\mathrm{held}}[\Delta Q_{\mathrm{near}}]\bigr)\\
&\quad+\bigl(\Eb^{\mathrm{held}}[\Delta Q_{\mathrm{far}}\mid\chi_a,\chi_{a'}]-\Eb^{\mathrm{held}}[\Delta Q_{\mathrm{far}}]\bigr).
\end{aligned}
\end{equation*}
We expect the second term to be negligible, that is, $\chi_a$ and $\chi_{a'}$ should not change the
far part on average, since their influence on the piece should decay with the distance to the disks.
This is an additional heuristic input. The independence of the band only shows that the increment is
the same function of $\chi_a$ and $\chi_{a'}$ for all routes, and not that this function depends
little on them. The first term is $O(h)$ in mean square, since conditioning does not increase the
variance, so that its mean square is at most $\Eb^{\mathrm{held}}|\Delta Q_{\mathrm{near}}|^2$, and
since the parts of the piece in the disks have shear increments at most of order $h$ in mean square.
To see the last point, let $\one_{\mathrm{disk}}$ be the indicator function of one of the disks.
Then $Q_{L,C}(\one_{\mathrm{disk}})$ is the shear increment of the parts of $L$ in this disk, and
the Bessel identity for $Q$, which will be stated in~\eqref{eq:bessel}, gives
\begin{equation*}
\Eb\sum_{(L,C)}\bigl|Q_{L,C}(\one_{\mathrm{disk}})\bigr|^2=\Eb\sum_{\fce\text{ in the disk}}\kappa_{\fce}\,|q_{\fce}|^2
\asymp\Bigl(\frac h\delta\Bigr)^2\delta^2=h^2,
\end{equation*}
since the disk contains $\asymp(h/\delta)^2$ faces, $|q_{\fce}|$ is of order $\delta$, and
$\kappa_{\fce}$, the conditional probability given $\omega$ that $\fce$ is crossed by an interface,
is of order one; the held law changes this expectation by a bounded factor only, by~\eqref{eq:held}.
The summation step below, in which $h$ tends to zero, requires this factor to be bounded uniformly
in $h$, which we expect, since, by the strong RSW estimates, the requirements outside the disk
should change the law of the configuration near the disk by a bounded factor only. Together with the
error $o(R)$ of the first line, this gives the second line of~\eqref{eq:g}.
\item Translation invariance. By the translation invariance of $\Pb_\beta$, the expectation is that
of the translated piece, whose disks are centered at $0$ and $a'-a$, up to a deterministic error
$O(h)$, which accounts for the position of the disks with respect to the lattice, since translation
invariance holds only for translations by lattice vectors. The band and the requirements outside it
do not matter. By the independence of the band and by locality, the function $\varphi_B$ depends on
neither of them, up to $o(R)$, and they change only the law of $\chi_a$ and $\chi_{a'}$, which, by
the second step, affects the expectation by $O(h)+o(R)$ only. The expectation is therefore
$g(a'-a)$, up to $O(h)+o(R)$.
\end{itemize}

\emph{Additivity.} Concatenating a held piece with displacement $v$ and a held piece with
displacement $v'$ gives a piece with displacement $v+v'$, which is also held at the junction. By
path independence, its increment is, up to $O(h)+o(|v|+|v'|)$ in mean square, that of a piece held
at its two ends only. Taking expectations shows that $g$ is additive, up to deterministic errors,
each of which is at most the mean-square size of the corresponding random error. These errors come
from the disks, which contribute $O(h)$, from the fluctuations of the increments around their
conditional expectations, and from the shear part around the closed interfaces used to compare the
routes; the last two are at most of order $m(R)^{1/2}$ in mean square, with $R\asymp|v|+|v'|$,
by~\eqref{eq:concentration}. The bound $m(R)=o(R^2)$ only makes them $o(|v|+|v'|)$, which is not
enough. With errors of order $|v|/\log(|v|/\delta)$, for instance, $g(v)$ could be equal to
$v\cos\log\log(|v|/\delta)$, and the limit below would not exist. We expect, however, the following
additivity estimate, which is deterministic, since $g$ is:
\begin{equation*}
g(v+v')=g(v)+g(v')+O\Bigl(h+\delta\,\Bigl(\frac{|v|+|v'|}\delta\Bigr)^{d/2+o(1)}\Bigr).
\end{equation*}
The reason is that, along a piece, the shear increment should consist of a systematic part, which
grows with the displacement, and of fluctuations. Around a closed interface of diameter $R$, the
displacement $V_{L,C}(1)$ vanishes, and only the fluctuations remain. Now
$Q_{L,C}(1)=\sum_{\fce}q_{\fce}\bar\tau^{j(\fce)}$, the sum running over the faces that border $C$,
and
\begin{itemize}[leftmargin=2em]
\item there are $(R/\delta)^{d+o(1)}$ such faces, by the count of Section~\ref{sec:ideas-switching}
at the scale $R$;
\item each term is of order $\delta$, since $q_{\fce}=A'_{\fce}\,\overline{p^0_{\fce}}$ and
$|p^0_{\fce}|=\delta$, the length of a side of $\fce$, by Section~\ref{sec:ideas-distorted}, whereas
$|A'_{\fce}|$, which measures the shear of $\fce$, is of order one at a typical face whose vertices
lie in different clusters.
\end{itemize}
Thus $Q_{L,C}(1)$ should behave like a sum of $(R/\delta)^{d+o(1)}$ weakly dependent centered terms
of order $\delta$, that is, be of order $\delta(R/\delta)^{d/2+o(1)}$ in mean square, the size of a
term times the square root of the number of terms. This gives $m(R)=\delta^2(R/\delta)^{d+o(1)}$, in
agreement with $m(R)=o(R^2)$, since $d<2$. Any exponent smaller than $1$ would do in place of $d/2$.

\emph{The linear function $\lambda$.} An approximately additive function of this kind is close to a
linear one, by an argument that goes back to Hyers~\cite{Hye41}. Since $g$ is deterministic, the
bounds that follow are ordinary bounds, except where we say otherwise. With $v'=v$, the additivity
estimate gives
\begin{equation*}
\bigl|2^{-n-1}g(2^{n+1}v)-2^{-n}g(2^nv)\bigr|\lesssim2^{-n}h+\delta\,(|v|/\delta)^{d/2+o(1)}\,2^{-n(1-d/2)+o(n)},
\end{equation*}
which is summable in $n$, since $d/2<1$. Hence, for every $v\in\C$, the limit
\begin{equation*}
\lambda(v)=\lim_{n\to\infty}2^{-n}g(2^nv)
\end{equation*}
exists, and, summing these corrections, for $v$ large compared with $h$,
\begin{equation*}
|g(v)-\lambda(v)|\lesssim h+\delta\,(|v|/\delta)^{d/2+o(1)} .
\end{equation*}
Like $g$, the function $\lambda$ is deterministic, and it is real-linear:
\begin{itemize}[leftmargin=2em]
\item it is additive, since the additivity estimate at $2^nv$ and $2^nv'$, divided by $2^n$, gives
$\lambda(v+v')=\lambda(v)+\lambda(v')$ as $n\to\infty$;
\item it is Lipschitz, and hence continuous, because $g$ grows at most linearly, with a constant $C$
that does not depend on $\delta$, in the ordinary sense:
\begin{equation*}
|g(v)|\le C\,(|v|+h) .
\end{equation*}
This follows from the Bessel identity, as for the disks above. Let $B$ be the band of a held piece
with displacement $v$. Then $Q_{L,C}(\one_B)$ is the shear increment of the piece in $B$, since $L$
meets $B$ only along this piece. As $B$ has area $\asymp|v|^2$, it contains $\asymp(|v|/\delta)^2$
faces, and
\begin{equation*}
\Eb\sum_{(L,C)}\bigl|Q_{L,C}(\one_B)\bigr|^2=\Eb\sum_{\fce\text{ in }B}\kappa_{\fce}\,|q_{\fce}|^2
\asymp\Bigl(\frac{|v|}\delta\Bigr)^2\delta^2=|v|^2 .
\end{equation*}
Therefore the shear increment of the piece in $B$ is $O(|v|)$ in mean square under the held law, and
its expectation is $O(|v|)$ in the ordinary sense; the parts in the disks add $O(h)$ to the
expectation, so that $|g(v)|\le C(|v|+h)$. Consequently,
\begin{equation*}
|\lambda(v)|=\lim_{n\to\infty}2^{-n}\,|g(2^nv)|\le\lim_{n\to\infty}2^{-n}\,C\,(2^n|v|+h)=C\,|v|,
\end{equation*}
and, by additivity, $|\lambda(v)-\lambda(w)|=|\lambda(v-w)|\le C\,|v-w|$;
\item an additive function is linear over the rationals, and, if it is continuous, it is linear over
the reals.
\end{itemize}
Thus $\lambda(v)=\eta\,v+\eta'\,\bar v$ for two constants $\eta$ and $\eta'$. The function $\lambda$
does not depend on $\delta$, since $\Pb_\beta$ at spacing $\delta$ is $\Pb_\beta$ at spacing $1$
scaled by $\delta$, so that $g$ at spacing $\delta$ is, up to $O(h)$, $\delta$ times $g$ at spacing
$1$, evaluated at $v/\delta$, and $\lambda$, which is linear, is not affected by bounded errors.
This defines $\eta$, which, for the barycenter map, depends only on $\beta$, and not on $\Omega$, on
$f$, or on the local configurations that hold the interfaces.

\emph{The law of large numbers.} The increment $\Delta V=\zeta_{a'}-\zeta_a$ of $V(1)$ along a piece
differs from $a'-a$ by $O(h)$, in mean square, since the ends of the piece lie in the disks, and
$\cB$ moves each site by at most the diameter of its cluster, which has exponential tails at the
scale $\delta$. Hence, by~\eqref{eq:g}, the bound on $|g-\lambda|$, and
$\lambda(v)=\eta v+\eta'\bar v$, along a held piece,
\begin{equation}\label{eq:lln-piece}
\Delta Q=\eta\,\Delta V+\eta'\,\overline{\Delta V}+O(h)+o(|\Delta V|),
\end{equation}
where the error terms are random, and their bounds hold in mean square under the held law only. This
is a law of large numbers, of the kind satisfied by the passage times of first-passage
percolation~\cite{HW65,Kin68}. It is not of the usual form, in which a sum of $N$ terms is close to
$N$ times their expectation, since the expected contribution of a face to $V$ and to $Q$ vanishes by
symmetry, and a piece of $N$ faces has a displacement of order only $N^{1/d+o(1)}\delta$, much
smaller than $N\delta$ since $d>1$.

\emph{Summation.} Cut the part of a macroscopic interface that borders a component into
$r^{-d+o(1)}$ pieces of size $r$, with $h\ll r$. By Minkowski's inequality, as $\delta\to0$ with $r$
and $h$ fixed, the errors $o(|\Delta V|)$ in~\eqref{eq:lln-piece}, summed over the boundedly many
pieces, tend to zero in mean square, and the errors $O(h)$ add up to $O(r^{-d+o(1)}h)$, which is
small if, for instance, $r=h^{1/2}$, since $d<2$. The variation of $f$ inside the pieces is handled
in the same way at smaller scales. Hence, in mean square under the held law,
\begin{equation}\label{eq:lln}
Q(f)\approx\sum_if(\zeta_i)\bigl[\eta\,\Delta V_i+\eta'\,\overline{\Delta V_i}\bigr]
\approx\eta\oint f\,dz+\eta'\oint f\,d\bar z .
\end{equation}
An arc of $\gamma$ between two good places, which will be defined in
Section~\ref{sec:ideas-difficulty}, becomes a piece of a held interface after a resampling of the
configuration near it. Since the resampling changes probabilities by a bounded factor
only,~\eqref{eq:lln} also holds along $\gamma$, in probability (Section~\ref{sec:transfer}).

\emph{Symmetry.} Consider the rotation $z\mapsto\tau z$ by $2\pi/3$ about the site $0$. It maps
$\delta\T$ onto itself, and the Ising model has the same coupling on all the edges, so that the
rotation preserves the infinite-volume law $\Pb_\beta$ in the plane. This is not true in
$\Omega_\delta$, whose boundary is not invariant under the rotation, and this is why $g$ is defined
under $\Pb_\beta$. The rotation maps every cluster of $\omega$ to a cluster of the rotated
configuration, and the barycenter of a set of points is a linear function of these points. Hence the
rotation commutes with $\cB$, in the sense that the image of $\tau z$ in the rotated configuration
is $\tau$ times the image of $z$. If we label the sides of the rotated face $\tau\fce$ by the images
of the sides of $\fce$, it follows that
\begin{equation*}
d_j(\tau\fce)=\tau\,d_j(\fce),\qquad\text{so that}\qquad p_{\tau\fce}=\tau\,p_{\fce},\quad q_{\tau\fce}=\tau\,q_{\fce},
\end{equation*}
by the uniqueness of the decomposition~\eqref{eq:pq}. The rotation also maps the interfaces to the
interfaces of the rotated configuration, with the same uncrossed sides, and held pieces to held
pieces. It therefore multiplies the displacement $\Delta V$ and the shear increment
$\Delta Q=\sum_{\fce}q_{\fce}\bar\tau^{j(\fce)}$ of every piece by $\tau$. Since it maps a held
piece with displacement $v$ to a held piece with displacement $\tau v$, and since it preserves
$\Pb_\beta$, $g(\tau v)=\tau g(v)$, up to $O(h)$, so that $\lambda(\tau v)=\tau\lambda(v)$. With
$\lambda(v)=\eta v+\eta'\bar v$, this reads
\begin{equation*}
\eta\,\tau v+\eta'\,\bar\tau\bar v=\tau\eta\,v+\tau\eta'\,\bar v\qquad(v\in\C),
\end{equation*}
so that $\eta'\bar\tau=\eta'\tau$, and therefore $\eta'=0$, since $\bar\tau\neq\tau$. Thus the term
with $\eta'$, which would make the limit anisotropic, vanishes. The reflection $z\mapsto\bar z$,
which reverses orientations, gives $\lambda(-\bar v)=-\overline{\lambda(v)}$, so that $\eta$ is
real. The term with $\eta$ is compatible with both symmetries, and it need not vanish. Therefore
\begin{equation*}
\lambda(v)=\eta\,v\qquad(v\in\C),
\end{equation*}
so that $\eta$ is the slope of the shear increments along an interface, and $Q\approx\eta V$ in the
sense of~\eqref{eq:lln}.

\emph{The meaning of $\eta$.} Along an interface, every face that borders the component contributes
a step $d_{j(\fce)}(\fce)$ to the image of the interface under $\cB$, and these steps add up to the
displacement $\Delta V$. Each step is the sum of the part $p_{\fce}\tau^{j(\fce)}$, which comes from
the rotation and the dilation of the face, and the part $q_{\fce}\bar\tau^{j(\fce)}$, which comes
from its shear. These parts are of order $\delta$ and point in all directions, so that most of them
cancel along a long piece, and only a net displacement remains. On average, the shears carry the
part $\eta\,\Delta V$ of this displacement, and the rotations and dilations carry the part
$(1-\eta)\,\Delta V$. Thus $\eta$ is a signed share, which may be negative if, on average, the shear
parts of the steps point opposite to the displacement. At $\beta=0$, there are no shears, and
$\eta=0$. For $\beta>0$, $\eta$ is an effective constant, which summarizes at large scales how the
random distortion of the faces acts on the interfaces, like the effective coefficients in
homogenization. Its value is not known and is not needed, since the argument only requires
$\eta\neq1$.

\subsection{The energy identity and the conclusion}\label{sec:ideas-energy}

By Section~\ref{sec:ideas-interface}, $\Eb[W\,V(f)]=\Eb[W\,Q(f)]$, and, by~\eqref{eq:lln}, with
$\eta'=0$, $Q(f)\approx\eta V(f)$ along $\gamma$. Hence, once this approximation is known to hold in
$L^1$, which we will check at the end of this section,
\begin{equation}\label{eq:conclusion}
\begin{aligned}
\Eb\bigl[W\,V(f)\bigr]&=\Eb\bigl[W\,Q(f)\bigr]\approx\eta\,\Eb\bigl[W\,V(f)\bigr],\\
\text{that is,}\qquad(1-\eta)\,\Eb\bigl[W\,V(f)\bigr]&\approx0 .
\end{aligned}
\end{equation}
This forces $\Eb[W\,V(f)]\approx0$, provided that $\eta\neq1$. We will derive the bound
$\eta<\frac12$ from an energy identity, which rests on the Bessel identities below. These identities
also give the bound $m(R)=o(R^2)$ used in Section~\ref{sec:ideas-selfavg}, and the bounds on second
moments needed to take expectations in~\eqref{eq:conclusion}.

\emph{The Bessel identities.} Color switching shows more than $E^\omega[W\,U(f)]=0$. For an
interface $L$, a component $C$ of $G_\omega\setminus L$, and a face $\fce$ with no open side, let
\begin{equation*}
Z_{\fce}(L,C)=\sum_{j=0}^2\tau^j\,\one_{\{e_j(\fce)^*\in C\}},
\end{equation*}
which, as in Section~\ref{sec:ideas-interface}, is equal to $\tau^{j(\fce)}$ if $L$ crosses $\fce$
and $\fce$ borders $C$, and to $0$ otherwise. The contour integrals of
Section~\ref{sec:ideas-selfavg} then read
\begin{equation*}
\begin{aligned}
U_{L,C}(f)&=\sum_{\fce}f(z_{\fce})\,p_{\fce}\,Z_{\fce}(L,C),\\
Q_{L,C}(f)&=\sum_{\fce}f(z_{\fce})\,q_{\fce}\,\overline{Z_{\fce}(L,C)},
\end{aligned}
\end{equation*}
the sums running over the faces with no open side. Given $\omega$, distinct faces are orthogonal
when one sums over all the pairs $(L,C)$. More precisely, for all faces $\fce$ and $\fce'$ with no
open side,
\begin{equation*}
E^\omega\sum_{(L,C)}Z_{\fce}(L,C)\,\overline{Z_{\fce'}(L,C)}=\kappa_{\fce}\,\one_{\{\fce=\fce'\}},
\end{equation*}
where $\kappa_{\fce}$ is the conditional probability that $\fce$ is crossed by an interface
(Proposition~\ref{prop:orth}). For $\fce=\fce'$, the sum is the indicator that $\fce$ is crossed by
an interface $L$, since such an interface is unique, and exactly one component of
$G_\omega\setminus L$ contains $e_{j(\fce)}(\fce)^*$. For $\fce\neq\fce'$, the orthogonality follows
from color switching, as does the face relation~\eqref{eq:face}, which is the case in which $\fce'$
is the added face ${\sface}$ (Corollary~\ref{cor:face}). Expanding the squares and using the
orthogonality gives the Bessel identities:
\begin{equation}\label{eq:bessel}
\begin{aligned}
E^\omega\sum_{(L,C)}|U_{L,C}(f)|^2&=\sum_{\fce}\kappa_{\fce}\,|f(z_{\fce})|^2\,|p_{\fce}|^2,\\
E^\omega\sum_{(L,C)}|Q_{L,C}(f)|^2&=\sum_{\fce}\kappa_{\fce}\,|f(z_{\fce})|^2\,|q_{\fce}|^2 .
\end{aligned}
\end{equation}
A face is crossed by an interface if and only if its three vertices do not all have the same spin.
Given $\omega$, the clusters receive independent fair spins, so that, if the vertices of $\fce$ lie
in $k_{\fce}$ distinct clusters, they all have the same spin with probability $2\cdot2^{-k_{\fce}}$,
and
\begin{equation*}
\kappa_{\fce}=1-2^{1-k_{\fce}}=
\begin{cases}
\frac34&\text{if }k_{\fce}=3,\\
\frac12&\text{if }k_{\fce}=2,\\
0&\text{if }k_{\fce}=1.
\end{cases}
\end{equation*}
If $\kappa_{\fce}<\frac34$, two vertices of $\fce$ lie in the same cluster, and hence have the same
image under $\cB$; the image of $\fce$ is then degenerate, and $|p_{\fce}|=|q_{\fce}|$.

\emph{The shear part around closed interfaces.} The second identity in~\eqref{eq:bessel} gives the
bound $m(R)=o(R^2)$ used in Section~\ref{sec:ideas-selfavg}. Let $f$ be the indicator function of a
region of area $A$. The region contains $\asymp A\delta^{-2}$ faces, and
$\Eb[\kappa_{\fce}|q_{\fce}|^2]\asymp\delta^2$, since $|q_{\fce}|$ is of order $\delta$. For a
closed interface $L$ in the region, $Q_{L,C}(f)=Q_{L,C}(1)$, and the other pairs only add
nonnegative terms, so that, by~\eqref{eq:bessel},
\begin{equation*}
\Eb\sum_{(L,C):\,L\text{ in the region}}|Q_{L,C}(1)|^2\lesssim A\delta^{-2}\cdot\delta^2=A .
\end{equation*}
We split this sum according to the diameter of $L$. By the definition of $m(R)$, the pairs in which
$L$ has diameter between $R$ and $2R$ contribute $m(R)$ times the expected number of these
interfaces, which is $\asymp AR^{-2}$ for $\delta\le R\le1$, as we will explain below. Therefore
\begin{equation*}
\sum_{R=\delta,2\delta,4\delta,\dots}AR^{-2}\,m(R)\lesssim A,\qquad\text{that is,}\qquad
\sum_{R=\delta,2\delta,4\delta,\dots}\frac{m(R)}{R^2}\lesssim1 .
\end{equation*}
Since the lattice $\delta\T$ is the lattice $\T$ scaled by $\delta$, and since,
by~\eqref{eq:mixing}, the law of the model away from $\partial\Omega$ is close to the
infinite-volume law $\Pb_\beta$, the ratio $m(R)/R^2$ depends only on $R/\delta$, up to a negligible
error, and the last sum is a convergent series in $R/\delta=1,2,4,\dots$. Its terms tend to zero, so
that $m(R)/R^2\to0$ as $\delta/R\to0$, and $Q_{L,C}(1)$ is small compared with $R$ in mean square.

The count $\asymp AR^{-2}$ comes from the RSW estimates, which show that the expected number of
interfaces of diameter between $R$ and $2R$ that meet a disk of radius $R$ is bounded above and
below.
\begin{itemize}[leftmargin=2em]
\item Lower bound. With probability bounded below, the disk contains a circuit of plus sites at
distance between $\frac34R$ and $R$ from its center, and a circuit of minus sites at distance
between $\frac12R$ and $\frac34R$; the outer boundary of the minus cluster of the inner circuit is
then an interface in the disk that surrounds the disk of radius $\frac12R$, so that its diameter is
between $R$ and $2R$.
\item Upper bound. An interface of diameter at least $R$ that meets a disk of radius $\frac18R$
crosses the annulus of radii $\frac18R$ and $\frac14R$ around it twice, and, given the interfaces
already explored, each new crossing has probability at most $1-c$, by the strong RSW estimates,
which are uniform in the boundary conditions; the number of these interfaces therefore has an
exponential tail, and hence a bounded expectation.
\end{itemize}
A region of area $A$ contains $\asymp AR^{-2}$ disjoint disks of radius $R$, and it is covered by
$\asymp AR^{-2}$ disks of radius $\frac18R$, which gives the count; the exponent $-2$ comes from the
dimension of the plane. We use the lower bound here. The upper bound was used in
Section~\ref{sec:ideas-selfavg}, and we will use it again in the proof that $\eta<\frac12$ below.

\emph{The energy identity.} By Section~\ref{sec:ideas-distorted},
$\frac{\sqrt3}4(|p_{\fce}|^2-|q_{\fce}|^2)$ is the signed area of the image of $\fce$ under $\cB$.
In the difference of the two identities in~\eqref{eq:bessel}, we may replace $\kappa_{\fce}$ by
$\frac34$, since $|p_{\fce}|=|q_{\fce}|$ if $\kappa_{\fce}<\frac34$, and add the faces with an open
side, whose images are degenerate as well. The signed areas then add up to the total area, even
though some images are degenerate or inverted. To see this, let $z$, $z'$, and $z''$ be the vertices
of $\fce$, in counterclockwise order, and $\zeta$, $\zeta'$, and $\zeta''$ their images under $\cB$.
The signed area of the image of $\fce$ is
\begin{equation*}
\tfrac12\,\mathrm{Im}\bigl(\bar\zeta\zeta'+\bar\zeta'\zeta''+\bar\zeta''\zeta\bigr),
\end{equation*}
a sum of terms $\frac12\mathrm{Im}(\bar\zeta\zeta')$ over the sides $zz'$ of $\fce$, oriented
counterclockwise. In the sum over the faces of a region $\Lambda$, every side inside $\Lambda$
appears twice, with opposite orientations, and its two terms cancel, since
$\mathrm{Im}(\bar\zeta\zeta')+\mathrm{Im}(\bar\zeta'\zeta)=0$. The remaining sum, over the sides on
the boundary of $\Lambda$, is the signed area enclosed by the image of $\partial\Lambda$, for every
configuration. Since $\cB$ moves every site by at most the diameter of its cluster, which is at most
$\ell$ except on an event of negligible probability (Section~\ref{sec:ES}), this area differs from
that of $\Lambda$ by $O(\ell\,|\partial \Lambda|)$. Applying this to small regions $\Lambda$, on
which $f$ is almost constant, and using $\frac34\cdot\frac4{\sqrt3}=\sqrt3$, we obtain
\begin{equation}\label{eq:energy}
\Eb\sum_{(L,C)}\bigl(|U_{L,C}(f)|^2-|Q_{L,C}(f)|^2\bigr)
=\frac34\,\Eb\sum_{\fce}\bigl(|p_{\fce}|^2-|q_{\fce}|^2\bigr)|f(z_{\fce})|^2\approx\sqrt3\iint|f|^2\,dA>0 .
\end{equation}
We call its left-hand side the energy. Since the distortion preserves the total area, the rotations
and dilations dominate the shears in the energy.

\emph{The bound $\eta<\frac12$.} For the pairs $(L,C)$ in which $L$ and $C$ are macroscopic, the
homogenization gives $Q_{L,C}(f)\approx\eta V_{L,C}(f)$, and hence
$U_{L,C}(f)=V_{L,C}(f)-Q_{L,C}(f)\approx(1-\eta)V_{L,C}(f)$, so that
\begin{equation*}
|U_{L,C}(f)|^2-|Q_{L,C}(f)|^2\approx\bigl((1-\eta)^2-\eta^2\bigr)|V_{L,C}(f)|^2=(1-2\eta)\,|V_{L,C}(f)|^2
\end{equation*}
(this uses~\eqref{eq:lln} in mean square, uniformly over the macroscopic arcs of interfaces in a
bounded region, since the expectations involve the squares $|Q_{L,C}(f)|^2$;
Section~\ref{sec:transfer} gives~\eqref{eq:lln} along actual interfaces only in probability, and the
passage to mean square, which requires the uniform integrability of these squares, is an additional
heuristic step). The other pairs, around single clusters or pockets, where the homogenization does
not apply, contribute to both sums in~\eqref{eq:energy}, but little to their difference. Here
$|U|^2-|Q|^2=|V|^2-2\,\mathrm{Re}(V\bar Q)$, and $V_{L,C}(f)$ is the contour integral of $f$ along
the part of $L$ that borders $C$ and that becomes a closed curve under $\cB$. If this curve has
diameter between $R$ and $2R$, then $|V_{L,C}(f)|\lesssim R^2$, and there are $O(R^{-2})$ such pairs
per unit area, by the upper bound above, since each crossing of an annulus by an interface borders
at most two such components. The pairs in which this curve has diameter less than $\rho$ therefore
contribute $O(\rho^2)$ to the expected sum of $|V|^2$, and, by the Cauchy--Schwarz inequality and
the Bessel identity for $Q$, $O(\rho)$ to that of $\mathrm{Re}(V\bar Q)$. Consequently,
\begin{equation*}
(1-2\eta)\,\Eb\sum_{(L,C)}|V_{L,C}(f)|^2\approx\sqrt3\iint|f|^2\,dA>0,
\end{equation*}
where the sum runs over the pairs in which this curve has diameter at least $\rho$, and $\delta\to0$
before $\rho\to0$. The sum on the left is nonnegative, so that $\eta<\frac12$, and in particular
$\eta\neq1$. The identity does not determine $\eta$, since this sum also depends on how the
interfaces pinch the domain into components at the scale of a single cluster, which is not
determined by the scaling limit of the interfaces.

\emph{The conclusion.} The relation~\eqref{eq:conclusion} involves expectations,
whereas~\eqref{eq:lln} says only that, along $\gamma$, the difference $Q(f)-\eta V(f)$ tends to zero
in probability as $\delta\to0$; almost sure convergence is not defined here, since the
configurations for different values of $\delta$ are not coupled. By~\eqref{eq:bessel}, however,
$Q(f)$ and $V(f)$ are bounded in mean square, uniformly in $\delta$, and $|W|\le1$. Since a sequence
of random variables that tends to zero in probability and is bounded in mean square tends to zero in
$L^1$, and since $\eta\neq1$,~\eqref{eq:conclusion} gives $\Eb[W\,V(f)]\to0$. By
Section~\ref{sec:ideas-interface}, it follows that $\iint\dbar f\,F_\delta\,dA\to0$ for every test
function $f$.

The functions $F_\delta$ take values in the triangle $\Delta$ (Section~\ref{sec:averaged}), and they
should be precompact, by the RSW estimates. Their subsequential limits $F$ satisfy
$\iint\dbar f\,F\,dA=0$ for all test functions $f$, so that they are holomorphic, by Weyl's lemma.
Their boundary values come from the boundary identity, the analogue of~\eqref{eq:H-boundary}: for a
site $u$ of $A_k$, and the edge $s_ku$ that joins $u$ to the ghost vertex $s_k$,
\begin{equation}\label{eq:fan}
F^\omega(s_ku)=w_k+(w_{k+1}-w_k)\,P^\omega\bigl([x_k,u]\pc A_{k+1}\bigr);
\end{equation}
taking expectations replaces $P^\omega$ by $\Pb$. As $u$ runs along $A_k$ from $x_k$ to $x_{k+1}$,
this probability increases from $0$ to $1$ in the limit. Hence, using the strong RSW estimates near
the boundary (Section~\ref{sec:cardy}), $F$ maps $A_k$ monotonically onto the side $[w_k,w_{k+1}]$
of $\Delta$. By the argument principle, $F$ is then a conformal map from $\Omega$ onto $\Delta$ with
$F(x_k)=w_k$, so that $F=\Phi_\Omega$. The identity~\eqref{eq:fan} at $u=b$ then gives Cardy's
formula (Section~\ref{sec:cardy}).

\subsection{The main difficulty}\label{sec:ideas-difficulty}

As in stochastic homogenization~\cite{Koz79,PV81}, the random distortion plays the role of rapidly
varying coefficients, and $\eta$ and $\eta'$ are effective coefficients. If $\eta'$ did not vanish,
the relation $\Eb[W\,V]=\Eb[W\,Q]$ would give
\begin{equation*}
(1-\eta)\,\dbar F+\eta'\,\partial F=0,
\end{equation*}
a Beltrami equation with the constant coefficient $\mu=-\eta'/(1-\eta)$, whose solutions, if
$|\mu|\neq1$, are $F(z)=G(z+\mu\bar z)$ with $G$ holomorphic; since $\eta'=0$, it is the
Cauchy--Riemann equation. Effective complex structures of this kind also arise in the homogenization
of random quasiconformal maps~\cite{ARST22,IM23} and in Beffara's extension of Smirnov's argument to
triangulations made of large triangular regions~\cite{Bef13}, where the faces neither collapse nor
reverse their orientation. Unlike in classical homogenization, however, the curve along which the
medium is probed here depends on the medium. This dependence produces $\eta$, and it is the main
difficulty of the homogenization (Section~\ref{sec:difficulty}).

We overcome it by holding interfaces at prescribed places, with the holding configurations and the
bands of Section~\ref{sec:ideas-selfavg} (Figures~\ref{fig:holding}, \ref{fig:routes},
and~\ref{fig:bands}). An actual interface should meet holding configurations at a dense set of
places, which we call good (Section~\ref{sec:transfer}). Two kinds of rare events must be excluded
(Section~\ref{sec:returns}). Let $c$ be the color on one side of the interface, and $c'$ the other
color; in Figure~\ref{fig:holding}, $c$ is plus.
\begin{itemize}[leftmargin=2em]
\item Near a point of the interface, no annulus contains a circuit with at most two sites of color
$c$, like the circuits of Figure~\ref{fig:holding}. By Menger's theorem, an annulus contains no such
circuit if and only if it is crossed by three disjoint paths of color $c$; together with the
interface, which gives an arm of color $c'$, this makes four arms. The bound~\eqref{eq:A4} says that
four such arms cross an annulus of radii $r<R$ with probability at most $C(r/R)^{\alpha_4}$, for
some $\alpha_4>1$.
\item The interface returns close to itself after a macroscopic excursion, and a short path of open
bonds joins the two passages. This path cuts the excursion off from $D$, so that $D$ is not
determined by the configuration near the interface. Around such a point, six arms with the colors
$(c,c,c',c,c,c')$ cross the annuli, and the bound~\eqref{eq:A6} says that this has probability at
most $C(r/R)^{\alpha_6}$, for some $\alpha_6>2$.
\end{itemize}
These two exponents must be larger than $1$ and $2$, respectively. To see why, consider a square of
side length $k$, measured in units of the mesh of a grid. An interface that crosses it passes near
$\gtrsim k$ of its $\asymp k^2$ grid points, one on each grid line, and it misses all good places in
the square only if all of these grid points are bad, which requires $\gtrsim k$ bad points; by
Markov's inequality, this has probability $\lesssim k^2\cdot k^{-\alpha_4}/k$. A return, on the
other hand, may occur at any of the $\asymp k^2$ grid points. Hence the bounds $\alpha_4>1$ and
$\alpha_6>2$ give
\begin{equation*}
k\cdot k^{-\alpha_4}\to0\quad \text{and}\quad k^2\cdot k^{-\alpha_6}\to0 \quad\text{as }k\to\infty.
\end{equation*}

\subsection{What is proved, and the plan of the rest of the paper}\label{sec:ideas-status}

As mentioned after Conjecture~\ref{conj:cardy}, except for the exact identities, which we indicate
as such, and for the results that we cite, the arguments of this paper are heuristic. The exact
identities all come from color switching on clusters, and they enter the strategy as follows:
\begin{itemize}[leftmargin=2em]
\item the uniform law of the interface configurations given $\omega$ (Section~\ref{sec:face}) is
color switching itself, and the other identities follow from it;
\item the face relation~\eqref{eq:face} gives $E^\omega[W\,U(f)]=0$, so that the rotations and
dilations do not contribute (Sections~\ref{sec:ideas-interface} and~\ref{sec:defect});
\item the orthogonality of distinct faces gives the Bessel identities, which make the shear part
small around closed interfaces (Sections~\ref{sec:ideas-selfavg} and~\ref{sec:ideas-energy}) and
give the energy identity, hence $\eta<\frac12$ (Section~\ref{sec:ideas-energy});
\item the boundary identity~\eqref{eq:fan}, and the fact that the values of $F^\omega$ lie in
$\Delta$, identify the limit as $\Phi_\Omega$ and give Cardy's formula (Section~\ref{sec:cardy}).
\end{itemize}
The uniform law follows from the classical correspondence between spin configurations and their
interfaces (Section~\ref{sec:face}). The orthogonality, the face relation~\eqref{eq:face}, and the
values in $\Delta$ will be proved in Appendix~\ref{app:exact} (Proposition~\ref{prop:orth},
Corollary~\ref{cor:face}, and Proposition~\ref{prop:convex}), and the identity~\eqref{eq:fan}
follows readily from the definition of the observable. Color switching on clusters is also used in
the proof of~\eqref{eq:A4}, to compare three arms of one color and one of the other with four
alternating arms (Section~\ref{sec:returns}). The barycenter map, the resampling and mixing
arguments, and the symmetry argument do not use it. The remaining inputs are of two kinds:
\begin{itemize}[leftmargin=2em]
\item results that are proved in the literature: the exponential decay of the Fortuin--Kasteleyn
connectivities~\cite{ABF87,DRT19}, and the spatial mixing~\eqref{eq:mixing}, which follows from it;
\item results that will appear in~\cite{SRSW}: the strong RSW estimates, the five-arm estimate, and
the arm bounds~\eqref{eq:A4} and~\eqref{eq:A6} of Section~\ref{sec:returns}.
\end{itemize}
The steps that remain heuristic are the homogenization of Section~\ref{sec:selfavg} and the energy
argument of Section~\ref{sec:ideas-energy}; the passage from holomorphicity to Cardy's formula,
through the precompactness of $F_\delta$ and its boundary values, should follow as in Smirnov's
argument, given the strong RSW estimates (Section~\ref{sec:cardy}).

The rest of the paper will give more details on these steps. In Section~\ref{sec:switching}, we will
explain color switching on clusters and the two parts of the contour integrals.
Section~\ref{sec:selfavg} will be devoted to the homogenization, and, in Section~\ref{sec:cardy}, we
will explain how Cardy's formula follows and conclude.

\section{Color switching on clusters}\label{sec:switching}

\subsection{The observable}\label{sec:averaged}

The faces of $\Omega_\delta$, the faces $s_kuu'$ for consecutive sites $u,u'$ of $A_k$, the faces
$s_ks_{k+1}x_{k+1}$, and the added face ${\sface}$ form a triangulation of the sphere, which we call
the completed sphere; we still write $\delta\cH$ for its dual. Every face has zero or two sides that
separate different spins, so that the interfaces are disjoint circuits of $\delta\cH$. If $W=w_k$,
the sides $s_ks_{k+1}$ and $s_{k+1}s_{k-1}$ of ${\sface}$ separate different spins, and $\gamma$ is
the interface through ${\sface}$. It starts at $x_{k+1}$ and ends at $x_{k-1}$, and between these
points it runs through $\Omega_\delta$ and, in places, along the arcs, so that at large scales it is
a curve in $\Omega$ that may touch the arcs (Figure~\ref{fig:observable}).

We define $D$ in $G_\omega$, and not in $\delta\cH$, because the interfaces, which never cross open
bonds, are circuits of $G_\omega$, on which color switching on clusters acts
(Section~\ref{sec:face}). Where $\gamma$ touches the arcs, or comes back within one lattice step of
itself, it passes on both sides of a site or of a cluster, and it cuts off a piece of the region
between itself and the arcs. Since a cluster is a single vertex of the contracted graph, these cuts
are the single-site cuts of Smirnov's argument. At $\beta=0$, if $W=w_k$ and $\sigma_{s_k}=+$, a
face lies in $D$ if and only if $E_k$ holds for it, so that $F_\delta=H$. Components of the
complement of an interface also enter the observable of~\cite{KSS25}. Finally, if $e^*$ is joined to
${\sface}$ in $G_\omega$, then $P^\omega(e^*\in D)=\frac14$ (Proposition~\ref{prop:convex}), so that
$F^\omega(e)=\sum_kw_k\,P^\omega(e^*\in D\mid W=w_k)$ is a convex combination of $w_1,w_2,w_3$;
otherwise $F^\omega(e)=0$. Hence the values of $F^\omega$ lie in $\Delta$. The values of $F_\delta$
lie in $\Delta$ as well, since $\Delta$ is convex and contains $0$, and since
$F^\omega_\delta(z)=4E^\omega[W\,\one_{\{z\in D\}}]$ is a combination of $w_1,w_2,w_3$ with
nonnegative weights whose sum $4P^\omega(z\in D)$ is at most $4P^\omega(e^*\in D)\le1$; here $e^*$
is an edge of $G_\omega$ at the face that contains $z$, or, if $z$ is enclosed by a circuit of open
bonds, at a face that lies outside the outermost such circuit and has a side on it.

\subsection{Smirnov's relation on the contracted graph}\label{sec:face}

An even subgraph of $G_\omega$ is a set of edges such that every vertex is the endpoint of an even
number of them, that is, a disjoint union of circuits. By the classical correspondence between spin
configurations and their domain walls on the sphere (see, e.g.,~\cite[Sections~2 and~5]{GJ09}), the
map that sends $\sigma$ to the union $I(\sigma)$ of the interfaces is two-to-one from the
$2^{k(\omega)}$ spin configurations that are constant on the clusters onto the even subgraphs of
$G_\omega$, where $k(\omega)$ now counts the clusters of the completed sphere, the ghost vertices
included. Hence, by~\eqref{eq:ES-intro},
\begin{equation}\label{eq:uniform}
P^\omega\bigl(I(\sigma)=\xi\bigr)=2^{1-k(\omega)}\qquad\text{for every even subgraph $\xi$ of $G_\omega$,}
\end{equation}
that is, given $\omega$, all interface configurations compatible with $\omega$ are equally likely.
Reversing the colors of the clusters on one side of a circuit toggles this circuit in $I(\sigma)$;
this is color switching on clusters. Here, toggling a set $S$ of edges in a subgraph $\xi$ means
replacing $\xi$ by the symmetric difference $\xi\triangle S$, which removes the edges of $S$ that
lie in $\xi$ and adds the other edges of $S$.

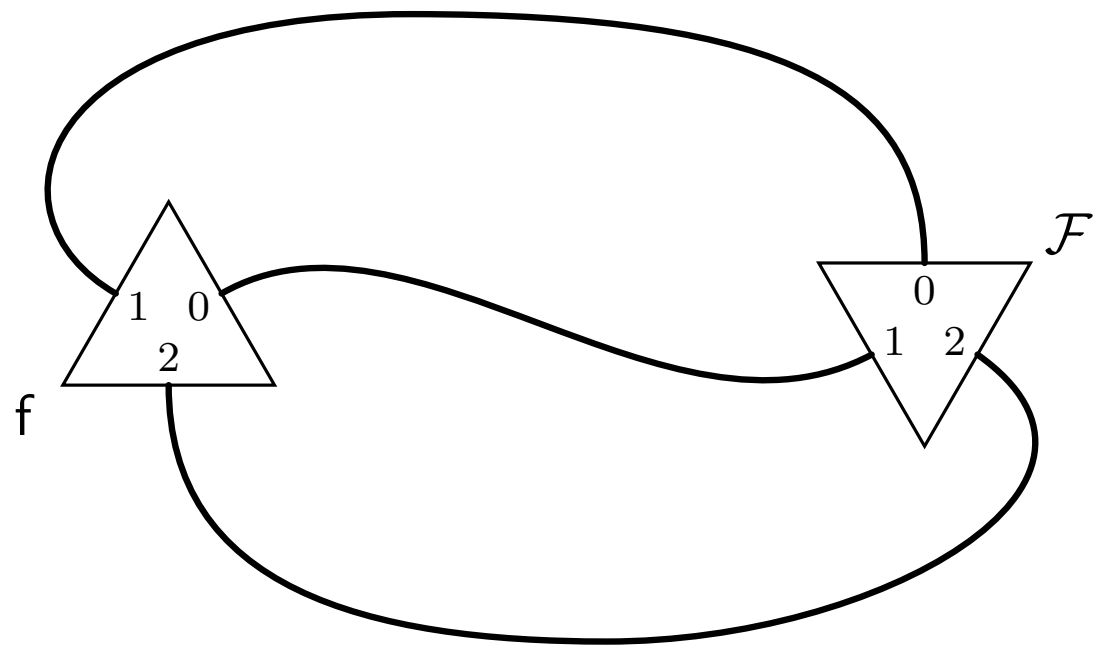
\begin{figure}[t]
\centering
\begin{tikzpicture}[scale=1,
  tri/.style={line width=0.6pt},
  str/.style={line width=1.2pt,line cap=round}]
\draw[tri] (0,0.808) -- (-0.7,-0.404) -- (0.7,-0.404) -- cycle;
\draw[tri] (5,-0.808) -- (4.3,0.404) -- (5.7,0.404) -- cycle;
\draw[str] (-0.35,0.202) .. controls (-1.25,0.72) and (-0.9,2.05) .. (1.6,2.05)
           .. controls (4.1,2.05) and (5,1.55) .. (5,0.404);
\draw[str] (0.35,0.202) .. controls (1.6,0.92) and (3.3,-0.92) .. (4.65,-0.202);
\draw[str] (0,-0.404) .. controls (0,-1.75) and (1.3,-2.1) .. (2.9,-2.1)
           .. controls (4.7,-2.1) and (6.55,-1.05) .. (5.35,-0.202);
\node at (0.2,0.115) {\scriptsize $0$};
\node at (-0.2,0.115) {\scriptsize $1$};
\node at (0,-0.22) {\scriptsize $2$};
\node at (5,0.22) {\scriptsize $0$};
\node at (4.8,-0.115) {\scriptsize $1$};
\node at (5.2,-0.115) {\scriptsize $2$};
\node at (-0.95,-0.6) {$\fce$};
\node at (5.95,0.6) {$\sface$};
\end{tikzpicture}
\caption{Three disjoint strands joining the sides of a face $\fce$ to those of the face ${\sface}$;
the sides of each face are numbered $0,1,2$ counterclockwise. Planarity reverses the cyclic order,
and the side $j$ of $\fce$ is joined to the side $1-j$ of ${\sface}$, modulo
$3$.}\label{fig:strands}
\end{figure}

\emph{The origin of the face relation.} Let $\fce\neq{\sface}$ be a face with no open side. If
$\gamma$ does not cross $\fce$, the three edges $e_j(\fce)^*$ all lie in $D$ or all lie outside $D$,
so that, for all $c_0,c_1,c_2\in\C$ with $c_0+c_1+c_2=0$,
\begin{equation*}
\sum_{j=0}^2c_j\,\one_{\{e_j(\fce)^*\in D\}}=c_{j(\fce)}\,\one_{\{\fce\text{ borders }D\}}.
\end{equation*}
With $c_j=\tau^j$, the left-hand side of~\eqref{eq:face} is therefore
$4E^\omega[W\,\tau^{j(\fce)}\one_{\{\fce\text{ borders }D\}}]$. Since $W$ records the side of
${\sface}$ that $\gamma$ does not cross,~\eqref{eq:face} says that the sides of $\fce$ and of
${\sface}$ that $\gamma$ does not cross are uncorrelated on the event that $\fce$ borders $D$.
By~\eqref{eq:uniform}, this is a statement about the subgraphs of $G_\omega$ whose only vertices of
odd degree are $\fce$ and ${\sface}$. Unless three disjoint strands of such a subgraph join the
three sides of $\fce$ to those of ${\sface}$, its contribution is the sum of a term that depends
only on the side of $\fce$ and a term that depends only on the side of ${\sface}$, and both cancel
since $\sum_j\tau^j=0$. If three strands join the side $j$ of $\fce$ to the side $\pi(j)$ of
${\sface}$, planarity reverses the cyclic order (Figure~\ref{fig:strands}), so that $\pi(j)=j_0-j$
for some $j_0$, and $\sum_j\tau^j\bar\tau^{\pi(j)}=\bar\tau^{j_0}\sum_j\tau^{2j}=0$. This
cancellation replaces Smirnov's color switching. With ${\sface}$ replaced by any face
$\fce'\neq\fce$, it gives the orthogonality relation of Section~\ref{sec:ideas-energy}; the proofs
are in Appendix~\ref{app:exact}.

\subsection{Boundary values}\label{sec:boundary-values}

The boundary identity~\eqref{eq:fan} holds for every $\omega$ and every $\beta$. For a site $u$ of
$A_k$, the edge $s_ku$ is never open, and $(s_ku)^*$ is joined to ${\sface}$ in $G_\omega$. If
$W=w_{k-1}$, then $\gamma$ separates $s_k$ from the other ghost vertices, so that
$(s_ku)^*\notin D$, and, by Proposition~\ref{prop:convex}, $F^\omega(s_ku)$ is a convex combination
of $w_k$ and $w_{k+1}$. Cutting the sphere along $s_ku$ identifies the weight of $w_{k+1}$ with the
crossing probability in~\eqref{eq:fan}; this follows readily from the definition of the observable.

\subsection{The defect is carried by one interface}\label{sec:defect}

By~\eqref{eq:pq} and~\eqref{eq:face}, for a face $\fce$ with no open side,
\begin{equation*}
\begin{aligned}
\sum_{j=0}^2F^\omega\bigl(e_j(\fce)\bigr)\,d_j(\fce)
&=p_{\fce}\sum_{j=0}^2\tau^j\,F^\omega\bigl(e_j(\fce)\bigr)+q_{\fce}\sum_{j=0}^2\bar\tau^j\,F^\omega\bigl(e_j(\fce)\bigr)\\
&=q_{\fce}\sum_{j=0}^2\bar\tau^j\,F^\omega\bigl(e_j(\fce)\bigr),
\end{aligned}
\end{equation*}
which is~\eqref{eq:defect}; if a side of $\fce$ is open, this contour integral vanishes
(Section~\ref{sec:ideas-switching}). If we regard $D$ as the union of the faces of $\delta\T$ that
are vertices of $D$ or lie in its holes, then Green's formula
$\oint_{\partial D}f\,dz=2\ii\iint_D\dbar f\,dA$ and~\eqref{eq:Fdelta} give~\eqref{eq:green}.
Finally, the identity of Section~\ref{sec:face}, with $c_j=d_j(\fce)$ and with $c_j=\bar\tau^j$,
together with~\eqref{eq:observable} and~\eqref{eq:defect}, gives
\begin{equation*}
\begin{aligned}
4E^\omega\bigl[W\,V(f)\bigr]&=\sum_{\fce}f(z_{\fce})\sum_{j=0}^2F^\omega\bigl(e_j(\fce)\bigr)\,d_j(\fce)\\
&=\sum_{\fce}f(z_{\fce})\,q_{\fce}\sum_{j=0}^2\bar\tau^j\,F^\omega\bigl(e_j(\fce)\bigr)=4E^\omega\bigl[W\,Q(f)\bigr],
\end{aligned}
\end{equation*}
the sums running over the faces with no open side, in agreement with~\eqref{eq:VQ}.

\section{Homogenization along the interfaces}\label{sec:selfavg}

This section gives more details on the homogenization of Section~\ref{sec:ideas-selfavg}, namely on
the role of symmetry (Section~\ref{sec:symmetry}) and on the derivation of~\eqref{eq:lln} through
held interfaces (Sections~\ref{sec:difficulty} to~\ref{sec:returns}).

\subsection{The role of symmetry}\label{sec:symmetry}

In terms of the differentials in~\eqref{eq:lln}, the rotation argument of the symmetry step of
Section~\ref{sec:ideas-selfavg} reads
\begin{equation*}
\begin{aligned}
dz&\mapsto\tau\,dz,\\
Q&\mapsto\tau\,Q,\\
d\bar z&\mapsto\bar\tau\,d\bar z .
\end{aligned}
\end{equation*}
Both sides of~\eqref{eq:lln} must transform in the same way, so that $\eta'=0$. The rotation by
$\pi$, which is a symmetry of every lattice, would not suffice, since it multiplies $dz$ and
$d\bar z$ by the same factor $-1$. In the language of field theory, the Beltrami coefficient
$\mu_{\fce}$ couples to the stress tensor, and $\eta'$ is a quantity of conformal spin $2$, which
the rotation by $2\pi/3$ forces to vanish.

\emph{Anisotropic interactions.} If the interactions on $\T$ are anisotropic, the rotation by
$2\pi/3$ is not a symmetry, and only $Q\approx\eta\oint f\,dz+\eta'\oint f\,d\bar z$ remains. The
limit then satisfies the Beltrami equation of Section~\ref{sec:ideas-difficulty}, in agreement with
the prediction that Cardy's formula then holds after a linear change of coordinates.

\subsection{The difficulty: interfaces depend on the medium}\label{sec:difficulty}

The interface $\gamma$ cannot be cut directly into pieces that can be resampled one at a time, with
$\gamma$ held at their ends and the rest of the configuration unchanged:
\begin{itemize}[leftmargin=2em]
\item an interface is determined by the whole configuration;
\item its shape is correlated with the shears along it;
\item the configurations in disjoint regions are not independent.
\end{itemize}
This is why we work with held interfaces in Section~\ref{sec:ideas-selfavg}. We will give more
details on them (Section~\ref{sec:constrained}), on the comparison of actual interfaces with held
ones (Section~\ref{sec:transfer}), and on the rare events in which this comparison fails
(Section~\ref{sec:returns}).

\subsection{Interfaces held at prescribed places}\label{sec:constrained}

We use the scales $\delta\ll h\ll r_h\ll1$, with $r_h=h^{1/2}$. In this subsection and the next one,
$h$ and $r_h$ are fixed while $\delta\to0$, and they tend to zero only at the end; constants that
depend on $h$ but not on $\delta$ do not affect the errors that tend to zero as $\delta\to0$, but
the errors of order $h$ must be bounded uniformly in $h$, as in the step on the disks in
Section~\ref{sec:ideas-selfavg}.

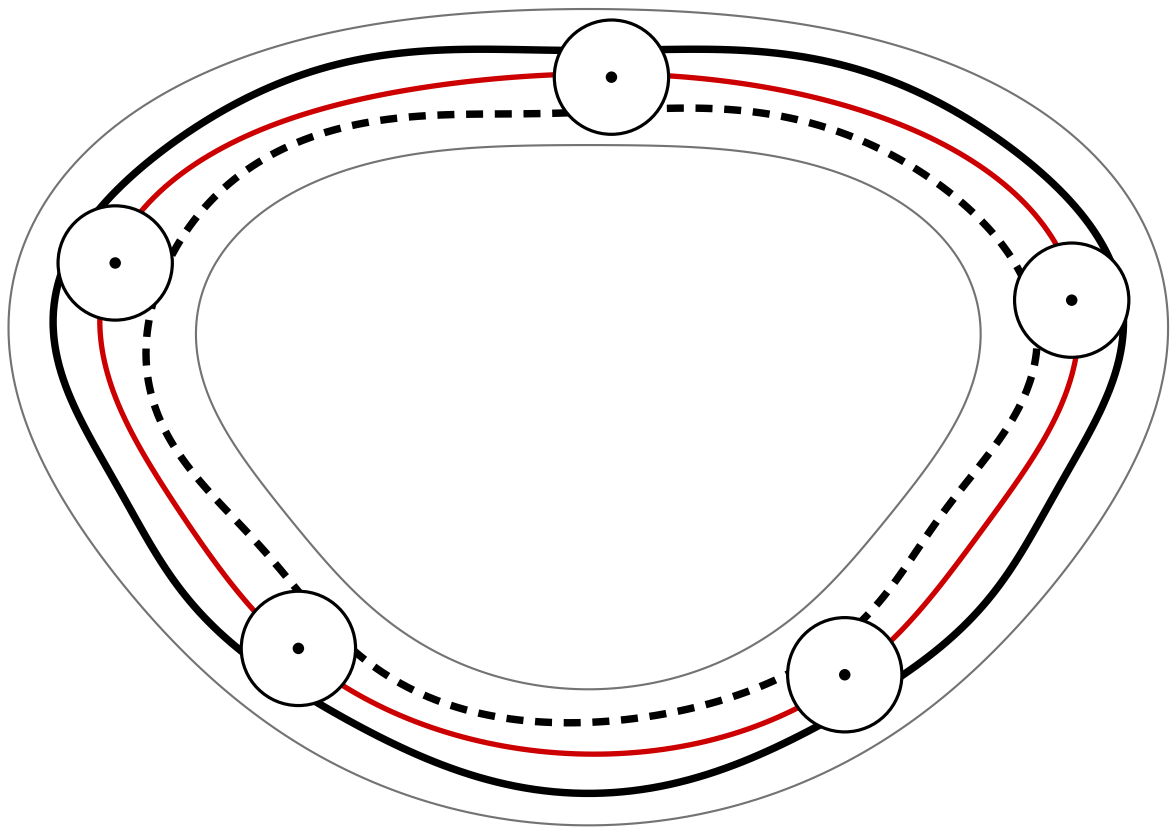
\begin{figure}[t]
\centering
\begin{tikzpicture}[scale=0.9,
  pl/.style={line width=1.4pt},
  mi/.style={line width=1.4pt,dash pattern=on 3.2pt off 2.2pt},
  intf/.style={line width=1pt,red!80!black},
  bnd/.style={line width=0.4pt,black!55}]
\def\rad#1{(2.5+0.18*sin(3*#1))}
\def\loopat#1#2{plot[domain=0:360,samples=180,smooth,variable=\t]
   ({1.4*(2.5+0.18*sin(3*\t)+#1+#2)*cos(\t)},{(2.5+0.18*sin(3*\t)+#1+#2)*sin(\t)})}
\draw[bnd] \loopat{0.5}{0};
\draw[bnd] \loopat{-0.5}{0};
\draw[pl] \loopat{0.23}{0.035*sin(7*\t)};
\draw[intf] \loopat{0}{0.025*sin(5*\t+20)};
\draw[mi] \loopat{-0.23}{0.035*sin(6*\t+50)};
\foreach \a in {15,87,159,231,303} {
  \pgfmathsetmacro{\rr}{2.5+0.18*sin(3*\a)}
  \filldraw[fill=white,line width=0.6pt] ({1.4*\rr*cos(\a)},{\rr*sin(\a)}) circle (0.42);
  \fill ({1.4*\rr*cos(\a)},{\rr*sin(\a)}) circle (1.2pt);
}
\end{tikzpicture}
\caption{An interface (red) held along a ring. Each small disk contains a holding configuration,
with circuits that the interface crosses exactly twice. Between consecutive disks, a plus path
(solid) and a minus path (dashed) keep the interface in the band.}\label{fig:cycle}
\end{figure}

\emph{Held interfaces.} The held interfaces are those of Section~\ref{sec:ideas-selfavg}. They are
held at a bounded number of disks of radius comparable to $h$, by holding configurations as in
Figure~\ref{fig:holding}, with the colors possibly exchanged, and kept in a band between consecutive
disks. The disks and the bands form a closed band, which we call a ring (Figure~\ref{fig:cycle}),
and a unique interface follows the ring and passes through every disk at the prescribed places. By
the strong RSW estimates~\cite{SRSW}, which are uniform in the boundary conditions, each requirement
holds with probability bounded below, regardless of the configuration outside, so
that~\eqref{eq:held} holds, with a constant $C_h$ that does not depend on $\delta$.

\emph{Almost independent pieces.} Two facts make the pieces of a held interface between consecutive
disks almost independent, and they underlie the locality step of Section~\ref{sec:ideas-selfavg}.
\begin{itemize}[leftmargin=2em]
\item The Markov property of the Edwards--Sokal coupling: given the configurations in the disks and
outside the ring, the configurations in the bands are independent.
\item Spatial mixing~\eqref{eq:mixing}, which follows from the exponential decay of the
Fortuin--Kasteleyn connectivities (Section~\ref{sec:ideas-selfavg}); see also~\cite{MOS94} for
general two-dimensional spin systems, and the correlation inequality of Ding, Song, and
Sun~\cite{DSS23}.
\end{itemize}
Hence~\eqref{eq:concentration} and~\eqref{eq:locality} hold for every piece, with $R$ comparable to
the distance between the disks at its ends, and the linearity and symmetry steps of
Section~\ref{sec:ideas-selfavg} give~\eqref{eq:lln-piece}, with $\eta'=0$.

\subsection{From held interfaces to actual interfaces}\label{sec:transfer}

A place $a$ on an interface is good if the disk of radius comparable to $h$ around $a$ contains a
holding configuration $\chi_a$, which holds the interface at $a$
(Section~\ref{sec:ideas-difficulty}).

\emph{An arc of $\gamma$.} Let $a_0,\dots,a_J$ be good places along a macroscopic arc of $\gamma$
that borders $D$, with consecutive places at distance comparable to $r_h$, so that
$J=r_h^{-d+o(1)}$; an upper bound $J\le r_h^{c-2}$, which suffices below, follows from bounds on
multiple crossings of annuli, as in~\cite{AB99}. As in Section~\ref{sec:ideas-selfavg}, let
$\zeta_j$ be the image of $a_j$, and write $\Delta V_j=\zeta_{j+1}-\zeta_j$ and $\Delta Q_j$ for the
increments of $V(1)$ and $Q(1)$, respectively, along the arc from $a_j$ to $a_{j+1}$.

\emph{Comparison with a held interface.} Resampling the configuration near the arc from $a_j$ to
$a_{j+1}$ creates a band along it, with a plus path and a minus path, and turns the arc into a piece
of a held interface, whose disks are those at $a_j$ and $a_{j+1}$ and a bounded number of disks that
close the band. Since a bounded number of requirements are imposed, the resulting law has a density
with respect to $\Pb^{\mathrm{held}}$ that is bounded independently of $\delta$, and a bounded
density changes the probabilities of events and the expectations of nonnegative quantities by only a
bounded factor. Hence~\eqref{eq:lln-piece}, with $\eta'=0$, carries over to the arc:
\begin{equation*}
\Delta Q_j=\eta\,\Delta V_j+O(h)+o(|\Delta V_j|),
\end{equation*}
where the error terms are random, and their bounds hold in mean square under the resulting law. At
the scale $h$, the arcs fall into a finite family of shapes, which depends on $h$ but not on
$\delta$, so that the errors $o(|\Delta V_j|)$ tend to zero uniformly over the arcs.

\emph{Summation.} The summation step of Section~\ref{sec:ideas-selfavg}, with $r=r_h$, then gives
$Q(f)\approx\eta V(f)$ along $\gamma$, in probability.

\subsection{Rare events and the arm bounds}\label{sec:returns}

For a sequence $\mathbf s$ of colors, let $\Arm_{\mathbf s}(z;r,R)$ be the event that disjoint paths
with the colors of $\mathbf s$, in counterclockwise order, cross the annulus $\{y\in\C:r<|y-z|<R\}$,
and recall that $\Pb_\beta$ is the infinite-volume law of Section~\ref{sec:ideas-selfavg}. We use
the following bounds, which will appear in~\cite{SRSW}: for some $\alpha_4>1$ and $\alpha_6>2$,
\begin{align}
\Pb_\beta\bigl(\Arm_{(+,+,+,-)}(z;r,R)\bigr)&\le C\Bigl(\frac rR\Bigr)^{\alpha_4},\tag{A4}\label{eq:A4}\\
\Pb_\beta\bigl(\Arm_{(+,+,-,+,+,-)}(z;r,R)\bigr)&\le C\Bigl(\frac rR\Bigr)^{\alpha_6},\tag{A6}\label{eq:A6}
\end{align}
and the same bounds for the color sequences obtained by exchanging plus and minus.

\begin{figure}[t]
\centering
\begin{tikzpicture}[scale=0.92,
  pl/.style={line width=1.4pt},
  mi/.style={line width=1.4pt,dash pattern=on 3.2pt off 2.2pt},
  cc/.style={line width=0.45pt}]
\begin{scope}
\draw[cc] (0,0) circle (0.4); \draw[cc] (0,0) circle (2.2);
\draw[pl] (90:0.4) .. controls (80:1.1) and (105:1.4) .. (98:2.2);
\draw[pl] (210:0.4) .. controls (225:1) and (195:1.5) .. (212:2.2);
\draw[pl] (-30:0.4) .. controls (-15:1) and (-55:1.5) .. (-38:2.2);
\draw[mi] (20:0.4) .. controls (40:1.1) and (10:1.6) .. (28:2.2);
\node at (98:2.5) {$+$}; \node at (212:2.5) {$+$}; \node at (-38:2.5) {$+$};
\node at (28:2.5) {$-$};
\node at (0,-2.9) {(a)};
\end{scope}
\begin{scope}[xshift=6.4cm]
\draw[cc] (0,0) circle (0.5); \draw[cc] (0,0) circle (2.2);
\draw[pl] (135:2.2) .. controls (145:1.2) and (-0.55,0.2) .. (0,0.2) .. controls (0.55,0.2) and
  (35:1.2) .. (45:2.2);
\draw[pl] (225:2.2) .. controls (215:1.2) and (-0.55,-0.2) .. (0,-0.2) .. controls (0.55,-0.2) and
  (-35:1.2) .. (-45:2.2);
\draw[mi] (-0.32,0) .. controls (-1,0) and (-1.5,0.12) .. (180:2.2);
\draw[mi] (0.32,0) .. controls (1,0) and (1.5,-0.12) .. (0:2.2);
\draw[line width=1.6pt,line cap=round,black!55] (0,0.2) -- (0,-0.2);
\node at (45:2.55) {$c$}; \node at (135:2.55) {$c$}; \node at (180:2.5) {$c'$};
\node at (225:2.55) {$c$}; \node at (-45:2.55) {$c$}; \node at (0:2.5) {$c'$};
\node at (0,-2.9) {(b)};
\end{scope}
\end{tikzpicture}
\caption{(a) Four arms with the colors $(+,+,+,-)$: no good place at this scale. (b) An interface
whose side of color $c$ passes twice near a point: the paths of color $c$ along the two passages
(solid) and the paths of color $c'$ between them (dashed) form six arms with the colors
$(c,c,c',c,c,c')$. A short path of open bonds (gray) between the two passages encloses a large
region.}\label{fig:arms}
\end{figure}
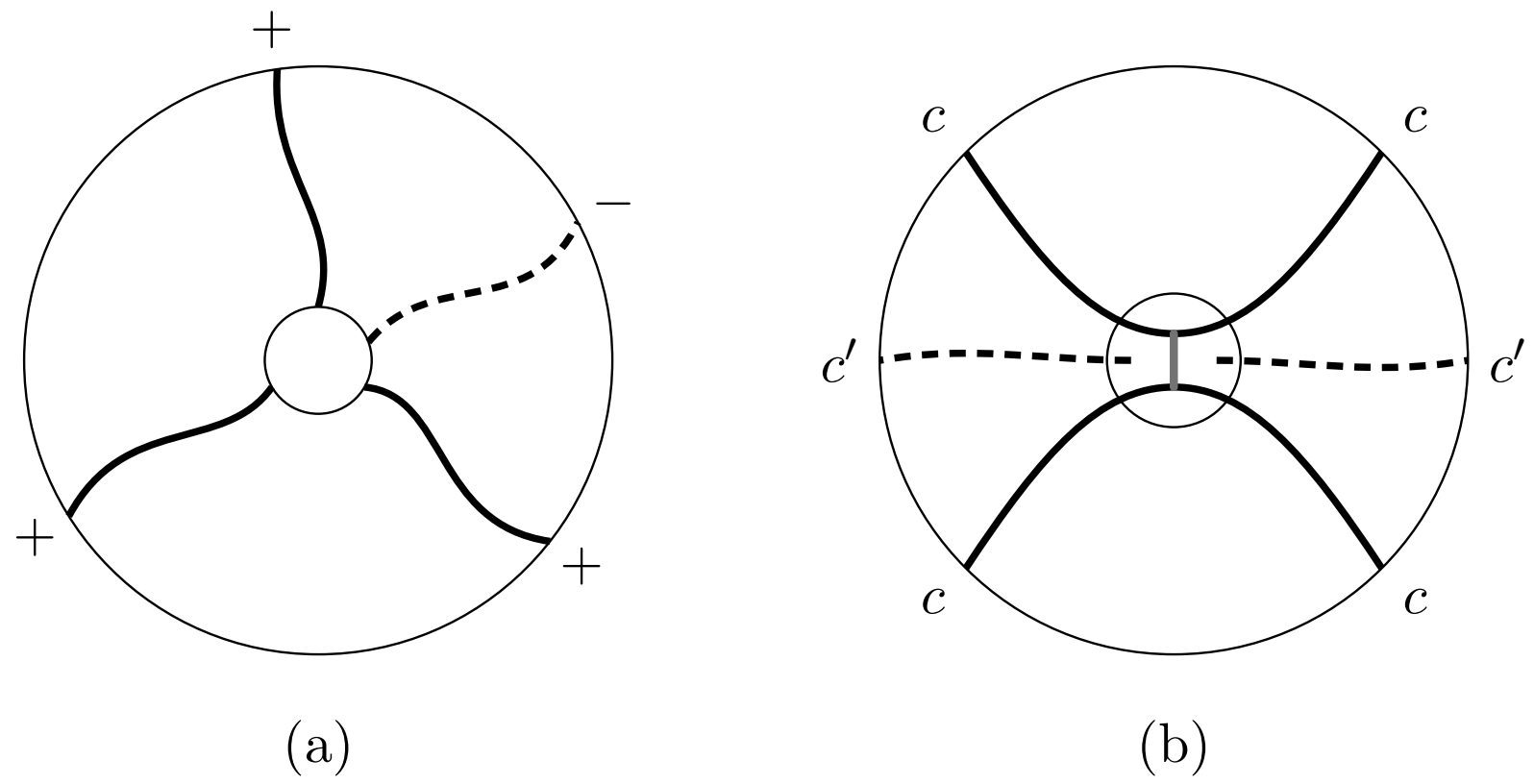

The two rare events of Section~\ref{sec:ideas-difficulty}, with the colors $c$ and $c'$ as there,
give arm events. The absence of a good place gives four arms with the colors $(c,c,c,c')$
(Figure~\ref{fig:arms}(a)), and, at a return, the sites of color $c$ along the two passages and the
paths of color $c'$ between them form six arms with the colors $(c,c,c',c,c,c')$
(Figure~\ref{fig:arms}(b)).

\emph{Returns.} A return may occur anywhere in a disk of radius $R$. A union bound over the
$\asymp(R/r)^2$ points of a grid of mesh $r$, with six arms between the scales $r$ and $r'$ around
each of them, where $r'$ is a lower bound on the diameter of the excursion, gives the bound
\begin{equation*}
\Bigl(\frac Rr\Bigr)^2\Bigl(\frac r{r'}\Bigr)^{\alpha_6}=\frac{R^2}{r'^{\alpha_6}}\,r^{\alpha_6-2},
\end{equation*}
which tends to zero as $r\to0$, since $\alpha_6>2$.

\emph{Good places.} Here the count is one-dimensional. Call a point of a grid of mesh $h$ bad if an
interface passes near it and crosses $n$ nested annuli around it, of aspect ratio $M$, and fewer
than three of these annuli contain a circuit with at most two sites of color $c$. In each of these
annuli, the interface gives an arm of color $c'$, so that, by Menger's theorem and~\eqref{eq:A4}, as
for the first rare event of Section~\ref{sec:ideas-difficulty}, each annulus fails with probability
at most $CM^{-\alpha_4}$. Failures in distinct annuli are almost independent (see~\cite{Kes87}, and
also~\cite{Nol08}, for such arguments in Bernoulli percolation), so that a grid point is bad with
probability at most $C_nM^{-\alpha_4(n-2)}$. An interface that crosses a square of side $kh$, with
$k\asymp M^n$ so that the annuli fit, passes near at least $k/2$ grid points, one on each grid line,
and it meets no good place only if these grid points are all bad, so that at least $k/2$ of the
$k^2$ grid points of the square are bad. By Markov's inequality, this has probability at most
\begin{equation*}
\frac{k^2}{k/2}\,C_nM^{-\alpha_4(n-2)}\lesssim M^{2\alpha_4-(\alpha_4-1)n}\longrightarrow0\quad \text{as } M\to\infty,
\end{equation*}
if $n>2\alpha_4/(\alpha_4-1)$. If the exponent in~\eqref{eq:A4} were equal to $1$, the bound would
be of order $M^2$ for every $n$, which is why the strategy requires $\alpha_4>1$.

\emph{The values of the exponents.} At $\beta=0$, the probabilities of four and of six arms, not all
of the same color, decay like $(r/R)^{5/4+o(1)}$ and $(r/R)^{35/12+o(1)}$~\cite{SW01}, and the same
exponents are expected for $\beta<\beta_c$. For every $\beta<\beta_c$, the five-arm event with the
colors $(+,-,+,-,+)$, which is obtained from the event in~\eqref{eq:A6} by removing one plus arm,
has probability comparable to $(r/R)^2$~\cite{SRSW}, so that~\eqref{eq:A6} requires a strictly
larger exponent. Both bounds rest on the strong RSW estimates. In~\eqref{eq:A6}, the sixth arm must
cross a region left free by the five others, which is unlikely given these five arms;
in~\eqref{eq:A4}, three plus arms and a minus arm can be compared with four alternating arms by
switching the colors of clusters of $\omega$, as in Section~\ref{sec:face}.

\section{From holomorphicity to Cardy's formula, and conclusions}\label{sec:cardy}

Once the subsequential limits $F$ of $F_\delta$ are known to be holomorphic, Cardy's formula should
follow as in Smirnov's argument, with one modification. The boundary of $\Omega_\delta$ may be rough
at all scales, and, unlike in percolation, the spins at distinct sites are not independent. We
therefore need the strong RSW estimates of~\cite{SRSW}, which are uniform in the boundary conditions
(Section~\ref{sec:background}), to control the behavior of the observable near the boundary. Then
$F=\Phi_\Omega$, by Proposition~\ref{prop:convex} and~\eqref{eq:fan}, as explained at the end of
Section~\ref{sec:ideas-energy}, and~\eqref{eq:fan} at $u=b^\delta$ gives
\begin{equation*}
\Pb_{\Omega_\delta,\beta}\bigl([x^\delta_1,b^\delta]\pc[x^\delta_2,x^\delta_3]\bigr)
=\frac{|\Eb F^\omega(s_1b^\delta)-w_1|}{|w_2-w_1|}\longrightarrow\frac{|\Phi_\Omega(b)-w_1|}{|w_2-w_1|},
\end{equation*}
which is~\eqref{eq:cardy}. The equality is exact; the convergence is what the strategy is designed
to show.

\emph{Conclusions.} We have explained a heuristic strategy to derive Cardy's formula for the spin
clusters of the Ising model on the triangular lattice at every fixed $\beta\in(0,\beta_c)$, and we
will give the details in~\cite{CardyLong}. Unlike in critical percolation, there is no exact
color-switching lemma for Smirnov's observable. In this strategy, holomorphicity should emerge at
large scales from the following mechanism:
\begin{itemize}[leftmargin=2em]
\item color switching on clusters: given $\omega$, Smirnov's relation
$\sum_j\tau^jF^\omega(e_j(\fce))=0$ holds exactly at every face of $\delta\T$ without an open side,
that is, on the graph in which every cluster is contracted to a point, which is not necessarily a
triangulation;
\item a random distortion: placing every cluster at its barycenter rotates, dilates, and shears the
faces, $d_j(\fce)=p_{\fce}\tau^j+q_{\fce}\bar\tau^j$; the rotations and dilations do not contribute,
$E^\omega[W\,U(f)]=0$, and the shears create a defect, $\Eb[W\,V(f)]=\Eb[W\,Q(f)]$, which is carried
by a single interface;
\item homogenization: along long interfaces, $Q\approx\eta V$, with an effective constant $\eta$,
and the rotation by $2\pi/3$ rules out a term with $d\bar z$;
\item area: since the distortion preserves the total area,
$(1-2\eta)\,\Eb\sum|V_{L,C}(f)|^2\approx\sqrt3\iint|f|^2\,dA$, so that $\eta<\frac12$, and
$(1-\eta)\,\Eb[W\,V(f)]\approx\Eb[W\,U(f)]=0$ gives $\Eb[W\,V(f)]\to0$.
\end{itemize}
With Cardy's formula and the strong RSW estimates at hand, one can prove the conformal invariance of
the interfaces of the spin clusters, namely their convergence to SLE$_6$, as conjectured by
Smirnov~\cite[Section~2.2, Conjecture~3]{Smi06}, with additional work that is by now relatively
standard (see~\cite{CN07,KS17}).

At $\beta_c$, the Fortuin--Kasteleyn clusters are macroscopic, and the map $\cB$ is not close to the
identity. Moreover, the strong RSW estimates, which hold for $\beta<\beta_c$~\cite{SRSW}, fail, and
spatial mixing holds only at a polynomial rate, instead of the exponential rate
of~\eqref{eq:mixing}. The strategy therefore breaks down, as it should, since the spin clusters of
the critical Ising model are not described by percolation. The structure on which the exact
identities rest, namely independent fair colors on the clusters of a random graph with small
clusters, is shared by the divide-and-color models of H\"aggstr\"om~\cite{Hag01} built on
subcritical Bernoulli percolation, which are critical on the triangular lattice~\cite{BCM09}, and we
expect the same mechanism to apply to them.

\appendix
\section{Proof of the exact identities}\label{app:exact}

Most arguments of this paper are heuristic, but the strategy rests on a few exact identities, which
we prove here. The proofs are short and elementary, and they use only color switching on clusters
and the counting of even subgraphs. Uniform even subgraphs and their toggling are classical (see,
e.g.,~\cite[Sections~2 and~5]{GJ09}), and the subgraphs with two odd vertices used below are the
loop configurations with disorders of~\cite{KS21}; see also~\cite{KSS25}.

We fix the configuration $\omega$ of open bonds on the completed sphere of
Section~\ref{sec:averaged}. We orient the completed sphere so that its orientation agrees with that
of the plane on $\Omega_\delta$, and we label the sides $e_0(\fce),e_1(\fce),e_2(\fce)$ of every
face $\fce$ counterclockwise for this orientation; for the outer face ${\sface}$, this is the
clockwise order in the plane. Starting from the uniform law~\eqref{eq:uniform}, we prove the
orthogonality relation behind the Bessel identities~\eqref{eq:bessel}, the face
relation~\eqref{eq:face}, and the fact that the values of $F^\omega$ lie in $\Delta$. The arguments
use only that the completed sphere is a finite triangulation of the sphere; they apply in the same
way to the faces of the plane, by completing a large region to a sphere.

Recall that $G_\omega$ is obtained from $\delta\cH$ by removing the edges $e^*$ of the open bonds
$e$, while all the vertices of $\delta\cH$ are kept. We identify subgraphs with their sets of edges,
and we recall from Section~\ref{sec:face} that toggling a set $S$ of edges in a subgraph $\xi$ means
replacing $\xi$ by the symmetric difference $\xi\triangle S$. For an interface $L$, the pairs
$(L,C)$ of Section~\ref{sec:ideas-selfavg} are formed with the connected components $C$ of the space
$G_\omega\setminus L$ obtained by removing all the points of $L$, its vertices included; an edge of
$G_\omega$ that is not on $L$ but has both endpoints on $L$ is thus a component by itself. For
points $m$ and $n$ of $G_\omega$ that are not on $L$, we write $m\sim_Ln$ if they lie in one
component of $G_\omega\setminus L$, and we write $\one_C(m)=\one_{\{m\in C\}}$. For two points $m$
and $n$ in the interiors of edges of $G_\omega$, let
\[
\mathcal K_\omega(m,n)=E^\omega\sum_{(L,C)}\one_C(m)\,\one_C(n)=E^\omega\sum_L\one_{\{m\sim_Ln\}},
\]
the expected number of interfaces $L$, through neither $m$ nor $n$, such that $m\sim_Ln$. For a face
$\fce$ and $j\in\{0,1,2\}$, let $m_j(\fce)$ be a point in the interior of the half of $e_j(\fce)^*$
that is adjacent to $\fce$. For a face $\fce$ with no open side, let
$Z_{\fce}(L,C)=\sum_j\tau^j\one_C(m_j(\fce))$, as in Section~\ref{sec:ideas-energy}, with
$e_j(\fce)^*$ replaced by the point $m_j(\fce)$. Finally, for $i\in\{1,2,3\}$, let $n_i$ be a point
in the interior of the half of $(s_{i-1}s_i)^*$ that is adjacent to ${\sface}$. If $W=w_k$, the
points $n_i$ with $i\neq k$ lie on $\gamma$, and $n_k$ lies in $D$; if an interface does not pass
through ${\sface}$, the three points $n_i$ lie in the component of its complement that contains
${\sface}$.

\begin{lemma}\label{lem:defects}
Let $m\neq n$ be points of one component of $G_\omega$, each a vertex or a point in the interior of
an edge, and let $P_{m,n}$ and $E_{m,n}$ be the probability and the expectation under the uniform
measure on the subgraphs $\xi$ of the subdivided graph whose vertices of odd degree are exactly $m$
and $n$, where the subdivided graph is obtained from $G_\omega$ by subdividing the edges that
contain $m$ or $n$ in their interiors. Then, for every circuit $L$ of $G_\omega$ that passes through
neither $m$ nor $n$,
\[
P_{m,n}(L\subset\xi)=P^\omega(L\text{ is an interface})\,\one_{\{m\sim_Ln\}}.
\]
In particular, if $m$ and $n$ lie in the interiors of edges, then
$\mathcal K_\omega(m,n)=E_{m,n}[N]$, where $N$ is the number of circuits of $\xi$.
\end{lemma}

\begin{proof}
The vertices of $L$ have even degree in $\xi$, hence degree $2$ if $L\subset\xi$. Thus the subgraphs
$\xi$ that contain $L$ are the unions of $L$ with the subgraphs of $G^L_\omega$ with odd vertices
$m$ and $n$, where $G^L_\omega$ is obtained from the subdivided graph by removing the vertices of
$L$ and the edges incident to them. Such subgraphs exist if and only if $m\sim_Ln$, and then, by
toggling a fixed path from $m$ to $n$, there are as many of them as there are even subgraphs of
$G^L_\omega$. The same toggling shows that $G_\omega$ has as many subgraphs with odd vertices $m$
and $n$ as even subgraphs. Therefore $P_{m,n}(L\subset\xi)$ is $\one_{\{m\sim_Ln\}}$ times the ratio
of the numbers of even subgraphs of $G^L_\omega$ and of $G_\omega$. By~\eqref{eq:uniform}, this
ratio is $P^\omega(L\text{ is an interface})$, since the even subgraphs of $G_\omega$ that contain
$L$ are the unions of $L$ with the even subgraphs of $G^L_\omega$. If $m$ and $n$ lie in the
interiors of edges, they have degree $1$ in $\xi$, so that no circuit of $\xi$ passes through them,
and summing the identity over $L$ gives $E_{m,n}[N]=\mathcal K_\omega(m,n)$.
\end{proof}

\begin{proposition}\label{prop:orth}
Let $\fce$ be a face with no open side. Then $Z_{\fce}(L,C)=\tau^{j(\fce)}$ if $L$ crosses $\fce$
and $\fce$ borders $C$, and $Z_{\fce}(L,C)=0$ otherwise; in particular,
$\sum_{(L,C)}|Z_{\fce}(L,C)|^2$ is the indicator that $\fce$ is crossed by an interface. Moreover,
for every face $\fce'\neq \fce$ with no open side,
\[
E^\omega\sum_{(L,C)}Z_{\fce}(L,C)\,\overline{Z_{\fce'}(L,C)}=0 .
\]
\end{proposition}

\begin{proof}
If $L$ does not pass through $\fce$, the three points $m_j(\fce)$ lie in one component of
$G_\omega\setminus L$, and $Z_{\fce}(L,C)$ is a multiple of $\sum_j\tau^j=0$. If $L$ passes through
$\fce$, the points on its two crossed sides lie on $L$, and only the third one may lie in $C$. This
proves the first statement.

For the second one, note that
\[E^\omega\sum_{(L,C)}Z_{\fce}(L,C)\overline{Z_{\fce'}(L,C)}=\sum_{j,k}\tau^j\bar\tau^k
\mathcal K_\omega(m_j(\fce),m_k(\fce')).\]
If $\fce$ and $\fce'$ lie in different components of $G_\omega$, all the terms vanish. Otherwise,
let $E_{\fce,\fce'}$ be the expectation under the uniform measure on the subgraphs $\xi$ whose
vertices of odd degree are the vertices $\fce$ and $\fce'$ of $\delta\cH$, in the graph obtained
from $G_\omega$ by subdividing it at the six points $m_j(\fce)$ and $m_k(\fce')$; these vertices
have degree $3$ in $G_\omega$, since none of their sides is open. Toggling the edges
$\fce m_j(\fce)$ and $\fce'm_k(\fce')$ of the subdivided graph is a bijection onto the subgraphs
whose odd vertices are $m_j(\fce)$ and $m_k(\fce')$, so that, by Lemma~\ref{lem:defects},
\[
\mathcal K_\omega\bigl(m_j(\fce),m_k(\fce')\bigr)=E_{\fce,\fce'}\bigl[N\bigl(\xi\triangle\{\fce m_j(\fce),\fce'm_k(\fce')\}\bigr)\bigr].
\]
Fix $\xi$, and remove small neighborhoods of $\fce$ and $\fce'$. The rest of $\xi$ consists of $N_0$
circuits and of strands that join the occupied half-edges at $\fce$ and $\fce'$, one or three at
each of these faces. We now count the circuits after the toggling. If $\fce$ has degree $1$ in
$\xi$, toggling at $\fce$ only moves the end of a path, and the number of circuits depends only on
$k$; similarly if $\fce'$ has degree $1$. If $\fce$ and $\fce'$ have degree $3$ and a single strand
joins them, at the half-edge $j_1$ of $\fce$ and the half-edge $k_1$ of $\fce'$, with a returning
strand at each face, the number of circuits is $N_0+\one_{\{j=j_1\}}+\one_{\{k=k_1\}}$. If three
strands join $\fce$ to $\fce'$, with the half-edge $j$ of $\fce$ joined to the half-edge $\pi(j)$ of
$\fce'$, the number of circuits is $N_0+\one_{\{k=\pi(j)\}}$. The terms that depend on at most one
of $j$ and $k$ vanish after summation against $\tau^j\bar\tau^k$, since $\sum_j\tau^j=0$. In the
last case, the three strands form, together with $\fce$ and $\fce'$, a $\theta$-graph on the sphere
(two vertices joined by three disjoint paths), and the cyclic orders of the three paths at the two
vertices are opposite (Figure~\ref{fig:strands}). Hence $\pi(j)=j_0-j$ for some $j_0$, and
$\sum_j\tau^j\bar\tau^{\pi(j)}=\bar\tau^{j_0}\sum_j\tau^{2j}=0$. Summing over $\xi$ proves the
proposition.
\end{proof}

\begin{corollary}\label{cor:face}
The face relation~\eqref{eq:face} holds for every face $\fce\neq{\sface}$ with no open side.
\end{corollary}

\begin{proof}
By the properties of $n_1,n_2,n_3$, and since $\sum_iw_i=0$, for every point $m$,
\[
\sum_{(L,C)}\one_C(m)\sum_{i=1}^3w_i\,\one_C(n_i)=W\,\one_{\{m\in D\}},
\]
so that, by~\eqref{eq:observable}, $F^\omega(e)=4\sum_iw_i\mathcal K_\omega(m_e,n_i)$, where $m_e$
is a point in the interior of $e^*$ other than $n_1,n_2,n_3$. Since $s_1,s_2,s_3$ appear
counterclockwise as seen from $\Omega_\delta$, they appear clockwise around ${\sface}$ for the
orientation of the sphere, and, for a suitable choice of the side $e_0({\sface})$, the points
$m_j({\sface})$ and $n_{1-j}$, with indices modulo $3$, lie on the same half-edge, and we may take
$m_j({\sface})=n_{1-j}$, so that
$Z_{\sface}(L,C)=\sum_i\tau^{1-i}\one_C(n_i)=\tau\sum_i\overline{w_i}\,\one_C(n_i)$. Therefore
\[
\sum_{j=0}^2\tau^jF^\omega\bigl(e_j(\fce)\bigr)=4\sum_{j,i}\tau^jw_i\,\mathcal K_\omega\bigl(m_j(\fce),n_i\bigr)
=4\tau\,E^\omega\sum_{(L,C)}Z_{\fce}(L,C)\overline{Z_{\sface}(L,C)},
\]
which vanishes by Proposition~\ref{prop:orth}, since no side of ${\sface}$ is open.
\end{proof}

\begin{proposition}\label{prop:convex}
Let $e$ be a bond that is not open and such that $e^*$ is joined to ${\sface}$ in $G_\omega$. Then
$P^\omega(e^*\in D)=\frac14$, so that
\[
F^\omega(e)=\sum_{k=1}^3w_k\,P^\omega\bigl(e^*\in D\bigm|W=w_k\bigr)
\]
is a convex combination of $w_1,w_2,w_3$. If $e^*$ is not joined to ${\sface}$ in $G_\omega$, then
$F^\omega(e)=0$. In all cases, $F^\omega(e)\in\Delta$.
\end{proposition}

\begin{proof}
If $e^*$ is not joined to ${\sface}$, it never lies in $D$, which is contained in the component of
$G_\omega$ containing ${\sface}$. Otherwise, let $m$ be a point in the interior of $e^*$ that does
not lie on the segments ${\sface}n_1$, ${\sface}n_2$, and ${\sface}n_3$, and let $P_{m,{\sface}}$
and $E_{m,{\sface}}$ be as in Lemma~\ref{lem:defects}, with $n={\sface}$, the graph being subdivided
at $n_1,n_2,n_3$ as well. By the properties of $n_1,n_2,n_3$,
\[
\mathcal K_\omega(m,n_i)=P^\omega\bigl(W=w_i,\ m\in D\bigr)+\bar N,
\]
where $\bar N$, which does not depend on $i$, is the expected number of interfaces $L$ that do not
pass through ${\sface}$ and satisfy $m\sim_L{\sface}$. On the other hand, toggling the edge
${\sface}n_i$ of the subdivided graph is a bijection onto the subgraphs whose odd vertices are $m$
and $n_i$, so that, by Lemma~\ref{lem:defects},
$\mathcal K_\omega(m,n_i)=E_{m,{\sface}}[N(\xi\triangle\{{\sface}n_i\})]$. Let $N_0$ be the number
of circuits of $\xi$ that do not pass through ${\sface}$; summing the identity of
Lemma~\ref{lem:defects} over these circuits gives $E_{m,{\sface}}[N_0]=\bar N$. The vertex
${\sface}$ has degree $1$ or $3$ in $\xi$. If it has degree $1$, the toggling only moves the end of
the path that starts at $m$, and the number of circuits is $N_0$. If ${\sface}$ has degree $3$, the
strand that leaves ${\sface}$ through the edge ${\sface}n_{i_1}$, say, ends at $m$, and the two
other half-edges are joined by a strand; after the toggling, the number of circuits is
$N_0+\one_{\{i=i_1\}}$. Hence
\[
P^\omega\bigl(W=w_i,\ m\in D\bigr)=\mathcal K_\omega(m,n_i)-\bar N=P_{m,{\sface}}\bigl(\deg_\xi{\sface}=3,\ i_1=i\bigr),
\]
and, summing over $i$, $P^\omega(m\in D)=P_{m,{\sface}}(\deg_\xi{\sface}=3)$. Finally, the edges at
$s_i$ are never open, so that the edges $e^*$, for the edges $e$ at $s_i$, form a circuit of
$G_\omega$. This circuit contains $(s_{i-1}s_i)^*$ and $(s_is_{i+1})^*$, but not
$(s_{i+1}s_{i-1})^*$. Toggling this circuit preserves $P_{m,{\sface}}$, and it changes the
occupation of exactly these two sides of ${\sface}$. The occupation vector of the three sides of
${\sface}$, which has an odd number of nonzero entries, is therefore uniform over the four such
vectors, and $P_{m,{\sface}}(\deg_\xi{\sface}=3)=\frac14$. The formula for $F^\omega(e)$ follows
from~\eqref{eq:observable}, since each event in~\eqref{eq:wstar} has probability $\frac14$; the
three conditional probabilities are nonnegative, and their sum is $4P^\omega(e^*\in D)=1$.
\end{proof}

\section*{Acknowledgments}
The author is grateful to his Ph.D.\ advisor, Hao Wu, for drawing his attention, when he was a
Ph.D.\ student at Tsinghua University, to the conjecture that the spin clusters of the
high-temperature Ising model have a conformally invariant scaling limit. The author also thanks
Federico Camia for discussions that helped him realize that it may be more convenient to work on the
random graph defined by the Fortuin--Kasteleyn configuration, and Dmitry Chelkak for conversations
that convinced him of the key role of a four-arm exponent strictly greater than~$1$.

\emph{Use of AI.} The author developed the proof strategy. Claude Opus 5.5 (Anthropic) was used to polish the language, improve the exposition, 
and create illustrative figures. Claude Opus 5.5 and Fable 5.2 (Anthropic), as well as GPT-6 Astra (OpenAI), were also used to check the details of 
the proofs on which this paper relies. The author takes full responsibility for the content of the paper.

\renewcommand{\bibliofont}{\small}
\newcommand{\etalchar}[1]{$^{#1}$}

\end{document}